\documentclass{amsart}
\usepackage{amsfonts,amssymb,amsthm,amsmath}
\usepackage[arrow,matrix]{xy}
\theoremstyle{plain}

\newtheorem*{Th*}{Theorem}

\newtheorem*{Cor*}{Corollary}

\theoremstyle{definition}

\theoremstyle{remark}

\numberwithin{equation}{section}

\makeatletter
\def\Set@Scallop[#1]#2#3{{#1}\Parens{#2}{#3}}
\newcommand\DeclareScalableOperator[2]{%
  \expandafter\def\csname#1\endcsname{\@ifnextchar[{{#2}\Set@Scallop}{{#2}\Set@Scallop[{}]}}
}
\makeatother

\DeclareScalableOperator{Ct}{\mathcal{C}} 
\DeclareScalableOperator{Cc}{\mathcal{C}_c} 
\DeclareScalableOperator{Lp}{\mathbf{L}} 
\DeclareScalableOperator{Hom}{\mathrm{Hom}} 
\DeclareScalableOperator{End}{\mathrm{End}} 
\DeclareScalableOperator{Aut}{\mathrm{Aut}} 
\DeclareScalableOperator{Uenv}{\mathfrak{U}} 
\DeclareScalableOperator{Cuenv}{\widehat{\mathfrak{U}}} 
\DeclareScalableOperator{ABer}{\Abs0{\mathrm{Ber}}}

\newcommand{\Fa}{For all }
\newcommand{\fa}{for all }

\newcommand{\fs}{for some }
\newcommand\mathfa[1][{}]{\quad\text{\fa{#1} }}

\newcommand{\scth}{such that }
\newcommand{\AND}{and}

\newcommand\mathtxt[1]{\quad\text{{#1}}\quad}

\newcommand{\nd}{\mathtxt\AND}

\newcommand\vphi{\varphi}
\newcommand\vrho{\varrho}

\newcommand\eps{\varepsilon}
\newcommand\nats{\mathbb{N}}
\newcommand\ints{\mathbb{Z}}

\newcommand\reals{\mathbb{R}}
\newcommand\cplxs{\mathbb{C}}

\newcommand\vvoid{\varnothing}
\newcommand\sle{\leqslant}
\newcommand\sge{\geqslant}

\DeclareMathOperator\Ad{\mathrm{Ad}}
\DeclareMathOperator\ad{\mathrm{ad}}
\DeclareMathOperator\GL{\mathrm{GL}}
\DeclareMathOperator\id{\mathrm{id}}
\DeclareMathOperator\pr{\mathrm{pr}}

\DeclareMathOperator\sdim{\mathrm{sdim}}
\DeclareMathOperator\sgn{\mathrm{sgn}}

\DeclareMathOperator\SMan{\mathrm{SMan}}

\makeatletter
\newcommand\Size[7][1]{
                                 \ifx#20%
                                        \def\r@l{}\def\r@m{}\def\r@r{}%
                                 \else%
                                    \ifx#21%
                                           \def\r@l{\bigl}\def\r@r{\bigr}\def\r@m{\bigm}%
                                    \else%
                                           \ifx#22%
                                                 \def\r@l{\Bigl}\def\r@r{\Bigr}\def\r@m{\Bigm}%
                                            \else%
                                                 \ifx#23%
                                                        \def\r@l{\biggl}\def\r@r{\biggr}\def\r@m{\biggm}%
                                                  \else
                                                        \ifx#24%
                                                        \def\r@l{\Biggl}\def\r@r{\Biggr}\def\r@m{\Biggm}%
                                                        \fi%
                                                  \fi%
                                            \fi%
                                      \fi%
                                 \fi%
                                 \ifx#10%
                                       \def\r@m{}%
                                 \fi%
                                 \r@l#3{#4}\r@m#5{#6}\r@r#7%
}%

\makeatother

\newcommand\Set[3]{
                                 \Size{#1}{\{}{#2}{|}{#3}{\}}%
}%
\newcommand\Dual[3]{
                                 \Size[0]{#1}{\langle}{#2}{,}{#3}{\rangle}%
}%
\newcommand\Parens[2]{
  \Size[0]{#1}{(}{#2}{}{}{)}
}
\newcommand\Bracks[2]{
  \Size[0]{#1}{[}{#2}{}{}{]}
}

\newcommand\Abs[2]{
  \Size[0]{#1}{\lvert}{#2}{}{}{\rvert}
}
\newcommand\Span[2]{
  \Size[0]{#1}{\langle}{#2}{}{}{\rangle}
}
\makeatletter

\newcommand{\IfUpperCase}[1]{\begingroup 
  \protected@edef\@tempa{\expandafter\@firstofone\@firstofone#1.}%
  \expandafter\IfUpperCasE \@tempa\delimiter}

\def\IfUpperCasE #1#2\delimiter{%
  \protected@edef\@tempa{\meaning#1\meaning a}%
  \ifnum \expandafter\IfUppercaSE\@tempa \IfUppercaSE
   \endgroup \expandafter\@firstoftwo
  \else
   \endgroup \expandafter\@secondoftwo
  \fi}

\@ifundefined{strip@prefix}{\def\strip@prefix#1>{}}{}
\def\@tempa{the letter }
\edef\@tempa{\expandafter\strip@prefix\meaning\@tempa}

\expandafter\def\expandafter\IfUppercaSE\expandafter#\expandafter1\@tempa#2#3\IfUppercaSE{\uccode`#2=`#2 }

\newif\ifuc@se
\def\setuc@se#1{\IfUpperCase{#1}{\uc@setrue}{\uc@sefalse}}

\def\theoremn@me#1{\ifuc@se \lowercase{\csname#1name\endcsname}\ignorespaces%
  \else \edef\@temp{\lowercase{\lowercase{\csname#1name\endcsname}}}\@temp\ignorespaces%
  \fi}
\def\theoremn@mes#1{\ifuc@se \lowercase{\csname#1names\endcsname}\ignorespaces%
  \else \edef\@temp{\lowercase{\lowercase{\csname#1names\endcsname}}}\@temp\ignorespaces%
  \fi}
  
\def\thmref#1#2{\setuc@se{#1}\lowercase{{\theoremn@me{#1}\lowercase{\ref{#1:#2}}}}}

\newcommand{\DefTheorem}[2]{\newenvironment{#1}[1][\empty]{\ignorespaces\begin{#2}\ifx##1\empty{}\else\lowercase{\label{#1:##1}}\fi\ignorespaces}{\end{#2}\ignorespacesafterend}}

\DefTheorem{Th}{theorem}
\DefTheorem{Prop}{proposition}
\DefTheorem{Cor}{corollary}
\DefTheorem{Lem}{lemma}
\DefTheorem{Def}{definition}
\DefTheorem{Rem}{remark}
\DefTheorem{Par}{para}

\newenvironment{Par*}{\ignorespaces\noindent\ignorespaces}{\ignorespacesafterend}

\makeatletter

\newif\if@smallmat
\newif\if@none
\newif\if@paren
\newif\if@brack
\newif\if@brace
\newif\if@vline
\newenvironment{Matrix}[2][1]
                                 {\ifx#20%
                                        \@smallmattrue%
                                  \else%
                                         \@smallmatfalse
                                  \fi%
                                  \ifx#11%
                                         \@nonefalse\@parentrue\@brackfalse\@bracefalse\@vlinefalse%
                                  \else%
                                       \ifx#12%
                                            \@nonefalse\@parenfalse\@bracktrue\@bracefalse\@vlinefalse%
                                        \else%
                                            \ifx#13%
                                                 \@nonefalse\@parenfalse\@brackfalse\@bracetrue\@vlinefalse%
                                            \else%
                                                 \ifx#14%
                                                       \@nonefalse\@parenfalse\@brackfalse\@bracefalse\@vlinetrue
                                                 \else%
                                                       \ifx#15%
                                                             \@nonefalse\@parenfalse\@brackfalse\@bracefalse\@vlinefalse%
                                                       \else%
                                                             \@nonetrue\@parenfalse\@brackfalse\@bracefalse\@vlinefalse%
                                                       \fi%
                                                 \fi%
                                            \fi%
                                        \fi%
                                   \fi%
                                   \if@smallmat%
                                        \if@none%
                                             \begin{smallmatrix}%
                                        \else%
                                            \if@paren%
                                                  \bigl(\begin{smallmatrix}%
                                            \else%
                                                  \if@brack%
                                                          \bigl[\begin{smallmatrix}%
                                                  \else%
                                                          \if@brace%
                                                               \bigl\{\begin{smallmatrix}%
                                                          \else%
                                                               \if@vline%
                                                                    \bigl\lvert\begin{smallmatrix}%
                                                                \else%
                                                                    \bigl\lVert\begin{smallmatrix}%
                                                                \fi%
                                                          \fi%
                                                  \fi%
                                            \fi%
                                        \fi%
                                   \else%
                                        \if@none%
                                             \begin{matrix}%
                                        \else%
                                            \if@paren%
                                                  \begin{pmatrix}%
                                            \else%
                                                  \if@brack%
                                                          \begin{bmatrix}%
                                                  \else%
                                                          \if@brace%
                                                               \begin{Bmatrix}%
                                                          \else%
                                                               \if@vline%
                                                                    \begin{vmatrix}%
                                                                \else%
                                                                    \begin{Vmatrix}%
                                                                \fi%
                                                          \fi%
                                                  \fi%
                                            \fi%
                                        \fi%
                                   \fi}%
                                  {\if@smallmat%
                                        \if@none%
                                             \end{smallmatrix}%
                                        \else%
                                            \if@paren%
                                                  \end{smallmatrix}\bigr)%
                                            \else%
                                                  \if@brack%
                                                          \end{smallmatrix}\bigr]%
                                                  \else%
                                                          \if@brace%
                                                               \end{smallmatrix}\bigr\}%
                                                          \else%
                                                               \if@vline%
                                                                    \end{smallmatrix}\bigr\rvert%
                                                                \else%
                                                                    \end{smallmatrix}\bigr\rVert%
                                                                \fi%
                                                          \fi%
                                                  \fi%
                                            \fi%
                                         \fi%
                                   \else%
                                        \if@none%
                                             \end{matrix}%
                                        \else%
                                            \if@paren%
                                                  \end{pmatrix}%
                                            \else%
                                                  \if@brack%
                                                          \end{bmatrix}%
                                                  \else%
                                                          \if@brace%
                                                               \end{Bmatrix}%
                                                          \else%
                                                               \if@vline%
                                                                    \end{vmatrix}%
                                                                \else%
                                                                    \end{Vmatrix}%
                                                                \fi%
                                                          \fi%
                                                  \fi%
                                            \fi%
                                        \fi%
                                   \fi}%
                         
\makeatother

\def\ger{\mathfrak}
\DeclareMathOperator\supp{\mathrm{supp}}
\DeclareMathOperator\str{\mathrm{str}}
\DeclareMathOperator\tr{\mathrm{tr}}
\DeclareMathOperator\gr{\mathrm{gr}}
\newcommand{\lBr}{[\kern-.65ex[}
\newcommand{\rBr}{]\kern-.65ex]}

\DeclareScalableOperator{Ber}{\mathrm{Ber}}
\DeclareScalableOperator{Vol}{\mathrm{Vol}}
\DeclareScalableOperator{Form}{\Omega}
\DeclareScalableOperator{Ind}{\mathrm{Ind}}
\DeclareScalableOperator{Coind}{\mathrm{Coind}}
\DeclareScalableOperator{Der}{\mathrm{Der}}

\begin{document}

\title[Harish-Chandra isomorphism]{The Harish-Chandra isomorphism\\ for reductive symmetric superpairs}
\author[A.~Alldridge]{Alexander Alldridge}

\address{Mathematisches Institut\\ Universit\"at zu K\"oln\\ Weyertal 86--90\\ 50931 K\"oln, Germany}
\email{alldridg@math.uni-koeln.de}
\thanks{Support from the IRTG 1133, SFB/TR 12, Leibniz group, and SPP 1388 grants, funded by Deutsche Forschungsgemeinschaft (DFG), is gratefully acknowledged}

\keywords{Harish-Chandra isomorphism, invariant differential operator, Lie supergroup, Lie superalgebra, symmetric superspace}
\subjclass[2010]{Primary 17B20. Secondary 58A50, 53C35}

\date{}

\begin{abstract}
	We consider symmetric pairs of Lie superalgebras which are strong\-ly reductive and of even type, and introduce a graded Harish-Chandra homomorphism. We prove that its image is a certain explicit filtered subalgebra of the Weyl invariants on a Cartan subspace whose associated graded is the image of Chevalley's restriction map on symmetric invariants. This generalises results of Harish-Chandra and V.~Kac, M.~Gorelik.
\end{abstract}

\maketitle

\section{Introduction}

Supermanifolds were developed in the 1970s by Berezin, Kostant and Leites as a rigorous mathematical framework for the quantum field theory of bosonic and fermionic particles. A particular class which appears naturally in connection with the representation theory of Lie supergroups is formed by the symmetric supermanifolds. In physics, they arise as the target spaces of non-linear SUSY $\sigma$-models, for instance, in the spectral theory of disordered systems \cite{zirnbauer-rsss}, and more recently, in the study of topological insulators \cite{schnyder-ryu-furusaki-ludwig}.

Given a symmetric superspace $X=G/K$, a fundamental object is the algebra $D(X)^G$ of $G$-invariant differential operators on $X$. For instance, the study of its $K$-invariant joint eigenfunctions (the \emph{spherical superfunctions}) should shed light on the regular $G$-representation on superfunctions on $X$. Indeed, in ongoing joint work with J.~Hilgert and M.~R.~Zirnbauer, we are developing the harmonic analysis on $X$ along these lines.

To state our main result, let us fix some notation. Let $(\ger g,\ger k,\theta)$ be a symmetric pair of complex Lie superalgebras, where the decomposition into $\theta$-eigenspaces is $\ger g=\ger k\oplus\ger p$. We assume that it is \emph{strongly reductive} and that there is an \emph{even Cartan subspace} $\ger a\subset\ger p_0$ (see the main text for precise definitions).

Let $I(\ger a)\subset S(\ger a)$ be the image of the restriction map $S(\ger p)^\ger k\toÊS(\ger a)$. We prove the following generalisation of Harish-Chandra's celebrated Isomorphism Theorem:

\begin{Th*}
	There is an explicitly defined filtered subalgebra $J(\ger a)\subset S(\ger a)$ with $\gr J(\ger a)=I(\ger a)$, and a short exact sequence of algebras,
	\[
	\xymatrix{%
		0\ar[r]&(\Uenv0{\ger g}\ger k)^\ger k\ar[r]&\Uenv0{\ger g}^\ger k\ar[r]^-{\Gamma}&J(\ger a)\ar[r]&0\ .%
	}
	\]
	In particular, if $(G,K,\theta)$ is any symmetric pair of Lie supergroups whose complexified infinitesimal pair is $(\ger g,\ger k,\theta)$, then the algebra $D(X)^G$ of complex $G$-invariant differential operators on $X=G/K$ is isomorphic to $J(\ger a)$. 
\end{Th*}

Here, $\Gamma$ is the super version of Harish-Chandra's well-known homomorphism. We stress that below, $J(\ger a)$ is defined in entirely explicit terms---otherwise the statement of the theorem would be somewhat vacuous. The image $I(\ger a)$ of Chevalley's restriction map was determined in joint work with J.~Hilgert and M.R.~Zirnbauer \cite{ahz-chevalley}, and this determination also forms the basis for the description of $J(\ger a)$.  

Our theorem covers the case of $(\ger g,\ger k,\theta)=(\ger k\oplus\ger k,\ger k,\mathrm{flip})$, the so-called `group type'; for this special case, the theorem is due to V.~Kac \cite{kac-laplace}, M.~Gorelik \cite{gorelik-kacconstruction}. Here, $I(\ger a)=J(\ger a)$, but this not true in general. Our method of proof is entirely different from that of Kac, Gorelik.

Rather, we follow Harish-Chandra's original analytic proof in the even case as closely as possible \cite{hc-sph1,helgason-gaga}. His idea is to consider a non-compact real form $(G,K,\theta)$ of $(\ger g,\ger k,\theta)$, and to use the elementary spherical functions defined on $X=G/K$ to derive the Weyl group invariance of the image of $\Gamma$. 

The graded counterpart of the first step in Harish-Chandra's proof appears to be impossible at first sight, due to the non-existence of compact (and hence, purely non-compact) real forms in the realm of real Lie superalgebras \cite{serganova-symmetric}. However, this can be addressed by using the category of \emph{cs} manifolds introduced by J.~Bernstein \cite{bernstein-qft,deligne-morgan}; it is a full subcategory of \emph{complex} graded ringed spaces. We discuss \emph{cs} manifolds at length in Appendix \ref{s:app-a}, and \emph{cs} Lie supergroups in Appendix \ref{s:app-b}. This framework allows us, in Section \ref{s:iwasawa}, to prove the existence of non-compact \emph{cs} forms (\thmref{Prop}{globalnccsformexists}), and of a global Iwasawa decomposition (\thmref{Prop}{iwasawa-global}).

In Appendix \ref{s:app-c}, we recall the basics of Berezin integration in the framework of \emph{cs} manifolds. In particular, we include the definition of the absolute Berezin integral (which is insensitive to changes of orientation in the even variables), and of invariant (absolute) Berezinians on homogeneous \emph{cs} manifolds, as developed by J.~Hilgert and the author for real supermanifolds \cite{ah-berezin}. In Section \ref{s:integral}, we employ these techniques to generalise the integral \emph{formul\ae{}} for the Iwasawa decomposition. In particular, we introduce certain non-zero joint eigen-superfunctions on $\Uenv0{\ger g}^\ger k$, and certain weighted orbital integrals. Compared to the even case, a complication is that $\int_K1\,\Abs0{Dk}=0$ (reflecting the maximal atypicality of the trivial $\ger k$-module), so that the technique of `invariant integration over $K$' cannot be as liberally employed as is customary for compact Lie \emph{groups}. We overcome this difficulty by introducing auxiliary superfunctions in Harish-Chandra's Eisenstein integral. (This amounts to considering joint eigenfunctions in a non-trivial $K\times K$-type of $\Gamma(\mathcal O_G)$.)

These analytic tools permit the proof of the fact that the image of $\Gamma$ is invariant under the \emph{even} Weyl group (\thmref{Prop}{hc-weylinv}). The remainder of the proof, in Section \ref{s:hc}, of the containment of $\mathop{\mathrm{im}}\Gamma$ in $J(\ger a)$ is entirely algebraic. It requires a rank reduction technique pioneered in the even case by Lepowsky \cite{lepowsky-hc}, and an explicit understanding of the relevant special cases of low rank. This is considerably more difficult than in the even case, where rank one is sufficient and the invariants in rank one are generated by the Casimir. (Both is false in general.)

Once the inclusion of the image of $\Gamma$ in $J(\ger a)$ has been established, a simple spectral sequence argument, together with the results from \cite{ahz-chevalley}, readily implies the Main Theorem (\thmref{Th}{hc}).

\medskip\noindent
\emph{Acknowledgements.} The author wishes to thank J.~Hilgert and M.R.~Zirnbauer for lending their invaluable expertise. 

\section{The Iwasawa decomposition}\label{s:iwasawa}

In the following, we will liberally employ standard notions of the theory of Lie supergroups, \emph{cf.}~\cite{leites,deligne-morgan,kostant-supergeom,schmitt-supergeom,scheunert-liesuperalgs,kac-liesuperalgs}. We will also make use of the notions introduced in Appendices \ref{s:app-a} and \ref{s:app-b}, in particular, the theory of \emph{cs} manifolds and groups.

\subsection{Non-compact \emph{cs} forms} 

First, we recall some definitions from \cite{ahz-chevalley}.

\begin{Def}[red]
	Let $\ger g$ be a complex Lie superalgebra and $\theta$ an even involutive automorphism of $\ger g$. We write $\ger g=\ger k\oplus\ger p$ for the $\theta$-eigenspace decomposition and say that $(\ger g,\ger k,\theta)$ is a \emph{symmetric superpair}. Occasionally, we may drop $\theta$ from the notation, although it is fixed as part of the data. We say that it is \emph{reductive} if $\ger g_0$ is reductive in $\ger g$, the centre $\ger z(\ger g)\subset\ger g_0$, and there exists a non-degenerate even $\ger g$- and $\theta$-invariant supersymmetric bilinear form $b$ on $\ger g$ (which is not part of the data). 
	
	In the following, we will only consider reductive symmetric superpairs. Mostly, will in fact need to impose a slightly more restrictive condition on $\ger g$.
	
	Indeed, let $\ger g$ be a Lie superalgebra such that $\ger z(\ger g)\subset\ger g_0$, $\ger g_0$ is reductive in $\ger g$, and $\ger g$ possesses a $\ger g$-invariant non-degenerate even supersymmetric bilinear form $b$. Then by \cite[Proposition 2.1, Lemma 2.3, Theorem 2.1, Corollary 2.1]{benayadi-quadratic}, $\ger g=\ger z(\ger g)\oplus\ger g'$ where $\ger g'=[\ger g,\ger g]$ if and only if $\ger g'$ is the direct sum of $b$-non-degenerate simple graded ideals. In this case, $\ger g$ will be called \emph{strongly reductive}. Any reductive symmetric superpair $(\ger g,\ger k,\theta)$ \scth $\ger g$ is a strongly reductive Lie superalgebra will call be called a \emph{strongly reductive} symmetric superpair. 
	
	We say that $(\ger g,\ger k,\theta)$ is of \emph{even type} if there exists $\ger a\subset\ger g$ \scth $\ger a$ is an \emph{even Cartan subspace} for $(\ger g,\ger k,\theta)$, \emph{i.e.}~$\ger a\subset\ger p_0$, $\ger a$ equals its centraliser $\ger z_{\ger p}(\ger a)$ in $\ger p$, and $\ad\ger a$ consists of semi-simple endomorphisms of $\ger g_0$.
\end{Def}

We shall have to choose particular \emph{cs} supergroup pairs whose underlying Lie superalgebra is part of a reductive symmetric superpair of even type. On the infinitesimal level, the conditions we will need are captured by the following definitions. 

\begin{Def}
	Let $(\ger g,\ger k,\theta)$ be a reductive symmetric superpair. A \emph{\emph{cs} form} of $(\ger g,\ger k)$ is a $\theta$-invariant real form $\ger g_{0,\reals}$ of $\ger g_0$ which is $b$-non-degenerate for some choice of the invariant form $b$. We write $\ger k_{0,\reals}=\ger g_{0,\reals}\cap\ger k$ and $\ger p_{0,\reals}=\ger g_{0,\reals}\cap\ger p$. 
	
	Let $\ger l$ be a real Lie algebra. Recall from \cite[Lemma 4.1, Definition 4.2]{borel-rss} that $\ger l$ is called \emph{compact} if the following equivalent conditions are fulfilled: the set $\ad\ger l\subset\End0{\ger l}$ consists of semi-simple elements with imaginary spectra; and $\ger l$ is the Lie algebra of a compact real Lie group. More generally, if $\vrho$ a is linear representation on a finite-dimensional real vector space $V$, then $\ger l$ is called \emph{$\vrho$-compact} if $\vrho(\ger l)$ generates a compact analytic subgroup of $\GL(V)$. 
	
	Denoting by $\ad_\ger g$ the adjoint action of $\ger g_0\subset\ger g$ on $\ger g$, a \emph{cs} form $\ger g_{0,\reals}$ will be called \emph{non-compact} if $\ger u_0=\ger k_{0,\reals}\oplus i\ger p_{0,\reals}$ is an $\ad_\ger g$-compact real form of $\ger g_0$; here, $\ad_\ger g$ denotes the adjoint action of $\ger g_0$ on $\ger g$. 
	
	Given a \emph{cs} form, a \emph{real even Cartan subspace} is a subspace $\ger a_\reals\subset\ger p_{0,\reals}$ whose complexification $\ger a$ is an even Cartan subspace of $(\ger g,\ger k,\theta)$. 
\end{Def}

We will prove the existence of non-compact \emph{cs} forms. The subtle point is the action of $\ger z(\ger g_0)$ on $\ger g$. It can be tamed by the following lemma. 

\begin{Lem}[centrecpt]
	Let $\ger g$ be a strongly reductive Lie superalgebra. There exists an $\ad_\ger g$-compact real form of $\ger z(\ger g_0)$. 
\end{Lem}

\begin{proof}
	Let $b$ be a non-degenerate invariant form on $\ger g$. As we have noted, $\ger g'$ is the direct sum of $b$-non-degenerate simple ideals $\ger s$, each of which has $\dim\ger z(\ger s_0)\sle1$. If $\ger z(\ger s_0)\neq0$, then $\ger s_1$ is the direct sum $\ger s_1=V_1\oplus V_2$ of simple $\ger s_0$-modules, and by \cite[Chapter II, \S 2.2, Corollary to Theorem 1]{scheunert-liesuperalgs} there is a unique element $C=C_\ger s\in\ger z(\ger s_0)$ \scth $\ad C=(-1)^j$ on $V_j$, $j=1,2$. The analytic subgroup of $\End0{\ger s}$ generated by $i\reals\cdot\ad C$ is isomorphic to $\mathbb T$. 
	
	Take any real form $\ger a$ of $\ger z(\ger g)$ and let $\ger b$ be sum of the real linear spans of the elements $iC_\ger s\in\ger z(\ger s_0)$, for all $b$-non-degenerate simple ideals $\ger s$ of $\ger g'$ with $\ger z(\ger s_0)\neq0$. Define $\ger z_\reals=\ger a\oplus\ger b$. The analytic subgroup of $\GL(\ger g)$ associated with $\ger z_\reals$ is isomorphic to a finite power of $\mathbb T$. 
\end{proof}

\begin{Lem}
	Let $(\ger g,\ger k,\theta)$ be a strongly reductive symmetric superpair. Then there exists a non-compact \emph{cs} form of $(\ger g,\ger k,\theta)$, and for any two such forms, their semi-simple derived algebras are conjugate by an inner automorphism of $\ger g_0$. If $(\ger g,\ger k,\theta)$ is, moreover, of even type, then there exists a real even Cartan subspace for any non-compact \emph{cs} form of $(\ger g,\ger k,\theta)$. 
\end{Lem}

\begin{proof}
	Following \thmref{Lem}{centrecpt}, we may choose an $\ad_\ger g$-compact real form $\ger z_\reals$ of $\ger z(\ger g_0)$. If $\ger s$ is a $b$-non-degenerate simple ideal of $\ger g'$ with $\ger z(\ger s_0)\neq0$, then $\dim\ger z(\ger s_0)=1$. Since $\ger z(\ger s_0)$ is $b$-non-degenerate, it is generated by a $b$-anisotropic vector. Either $\ger s$ is $\theta$-invariant, in which case so is $\ger z(\ger s_0)$, or $\theta(\ger s)$ is a distinct but isomorphic $b$-non-degenerate simple ideal of $\ger g'$. In the latter case, we have $\theta(C_\ger s)=\pm C_{\theta(\ger s)}$ for the elements $C=C_\ger s$ considered in the proof of \thmref{Lem}{centrecpt}, by their mere definition. Hence, by construction, we may assume that $\ger z_\reals$ is $\theta$-invariant and $b$-non-degenerate. 
	
	Let $\ger g_0'=[\ger g_0,\ger g_0]$. By \cite[pp.~154--155]{loos-ss1}, there exists a $\theta$-invariant compact real form $\ger u_0'$ of $\ger g_0'$ (any real form has a Cartan decomposition compatible with $\theta$, and the compact real form associated with it is $\theta$-stable), and it is unique up to inner automorphisms of $\ger g_0'$. We set 
	\[
	\ger g_{0,\reals}=\ger z_\reals\cap\ger k_0\oplus\ger u_0'\cap\ger k_0\oplus i\ger z_\reals\cap\ger p_0\oplus i\ger u_0'\cap\ger p_0\ .
	\]
	Then $\ger g_{0,\reals}$ is a $\theta$-stable and $b$-non-degenerate real form of $\ger g_0$. Defining $\ger k_{0,\reals}=\ger g_{0,\reals}\cap\ger k_0$ and $\ger p_{0,\reals}=\ger g_{0,\reals}\cap\ger p_0$, we see that $\ger u_0=\ger k_{0,\reals}\oplus i\ger p_{0,\reals}$ is indeed $\ad_\ger g$-compact, so that $\ger g_{0,\reals}$ is a non-compact \emph{cs} form of $(\ger g,\ger k,\theta)$.
	
	Assume that $(\ger g,\ger k,\theta)$ is of even type, and that $\ger g_{0,\reals}$ is a non-compact \emph{cs} form. There exists a maximal Abelian subspace $\ger a_\reals\subset\ger p_{0,\reals}$. By assumption, there exists an even Cartan subspace $\ger a\subset\ger p_0$. Now $\ger a_\reals\otimes\cplxs$ is the centraliser in $\ger p_0$ of any regular element of $\ger a_\reals$ (such elements exist), and any semi-simple element of $\ger p_0$ is conjugate under the adjoint group of $\ger k_0$ to an element of $\ger a$ \cite[Chapter III, Proposition 4.16]{helgason-gaga}. Thus, we may assume $\ger a_\reals\otimes\cplxs\subset\ger a$. Then $\ger a=\ger z_{\ger p_0}(\ger a)\subset\ger z_{\ger p_0}(\ger a_\reals)=\ger a_\reals\otimes\cplxs$, and hence the claim. 
\end{proof}

\subsection{Restricted roots and the Iwasawa decomposition}

\begin{Par}[roots]
	Let $(\ger g,\ger k,\theta)$ be a reductive symmetric superpair of even type, assume given a non-compact \emph{cs} form $\ger g_{0,\reals}$ (which always exists if $(\ger g,\ger k,\theta)$ is strongly reductive), let $\ger a_\reals\subset\ger p_{0,\reals}$ be a real even Cartan subspace, and set $\ger a=\ger a_\reals\otimes\cplxs$. We shall fix these data from now on.
	
	By assumption, $\ger g_0$ is reductive in $\ger g$, and $\ger a$ is a commutative subalgebra consisting of semi-simple elements. In particular, $\ger g$ is a semi-simple $\ger a$-module, and we may decompose it as 
	\begin{equation}\label{eq:restrrootdecomp}
		\ger g=\ger m\oplus\ger a\oplus\bigoplus\nolimits_{\lambda\in\Sigma}\ger g_{\ger a}^\lambda\mathtxt{where}\ger m=\ger z_{\ger k}(\ger a)
	\end{equation}
	is the centraliser of $\ger a$ in $\ger k$, and for $\lambda\in\ger a^*$, 
	\[
	\ger g_{\ger a}^\lambda=\Set1{x\in\ger g}{\forall h\in\ger a\,:\,[h,x]=\lambda(h)x}\nd\Sigma=\Set1{\lambda\in\ger a^*\setminus0}{\ger g_{\ger a}^\lambda\neq0}\ .
	\]
	We also define $\ger g_{j,\ger a}^\lambda=\ger g_j\cap\ger g_{\ger a}^\lambda$ and $\Sigma_j=\Set1{\lambda\in\ger a^*\setminus0}{\ger g_{j,\ger a}^\lambda\neq0}$. Then we have $\Sigma=\Sigma_0\cup\Sigma_1$, but the union may not be disjoint. Occasionally, we will write $\Sigma(\ger g:\ger a)=\Sigma$ and $\Sigma(\ger g_j:\ger a)=\Sigma_j$. Since $\ger u_0=\ger k_{0,\reals}\oplus i\ger p_{0,\reals}$ is by assumption a compact Lie algebra, the even restricted roots $\lambda\in\Sigma_0$ are real on $\ger a_\reals$. Let $\ger g_{0,\reals,\ger a}^\lambda=\ger g_{0,\reals}\cap\ger g_{0,\ger a}^\lambda$ \fa $\lambda\in\Sigma_0$ and $\ger m_{0,\reals}=\ger z_{\ger k_{0,\reals}}(\ger a_\reals)$. 
\end{Par}

\begin{Par}[possyst]
	Let $\Sigma^+\subset\Sigma$ be a \emph{positive system}, \emph{i.e.}~a subset \scth $\Sigma=\Sigma^+\,\dot\cup\,{-\Sigma^+}$ and $\Sigma\cap(\Sigma^++\Sigma^+)\subset\Sigma^+$. Let $\Sigma_j^+=\Sigma_j\cap\Sigma^+$. Then $\Sigma_0^+=\Sigma_0\cap\Sigma^+$ is a positive system of the root system $\Sigma_0$. Set 
	\[
		\ger n=\bigoplus\nolimits_{\lambda\in\Sigma^+}\ger g_{\ger a}^\lambda\nd\ger n_j=\ger g_j\cap\ger n\ .
	\]
	By the assumptions on $\Sigma^+$, $\ger n=\ger n_0\oplus\ger n_1$ is an $\ger a$-invariant subsuperalgebra. Moreover, $\ger n_{0,\reals}=\ger g_{0,\reals}\cap\ger n$ (which is a real form of $\ger n_0$) is an $\ger a_\reals$-invariant nilpotent subalgebra of $\ger g_{0,\reals}$. Since the roots in $\Sigma_0$ are real on $\ger a_\reals$, 
\[
	\ger n_{0,\reals}=\bigoplus\nolimits_{\lambda\in\Sigma_0^+}\ger g_{0,\reals,\ger a}^\lambda\ .
\]	

By \cite[Proposition~5]{kostant-rallis}, we have $\dim\ger m_0-\dim\ger a=\dim\ger k_0-\dim\ger p_0$. We also need the corresponding result for the odd part of $\ger g$, and this is the content of the following lemma. 
\end{Par}

\begin{Lem}[centdim-fmla]
	Let $x\in\ger p_0\,$. Then 
	\[
		\dim\ger z_{\ger k_1}(x)-\dim\ger z_{\ger p_1}(x)=\dim\ger k_1-\dim\ger p_1\ .
	\]
\end{Lem}

\begin{proof}
	We proceed exactly as in \cite[proof of Proposition~5]{kostant-rallis}. Choose a $\ger g$- and $\theta$-invariant non-degenerate even supersymmetric form $b$ on $\ger g$. Then, certainly $\ger z_{\ger g_1}(x)=\ger z_{\ger k_1}(x)\oplus\ger z_{\ger p_1}(x)$, and $\ger g_1/\ger z_{\ger g_1}(x)$ carries a \emph{symmetric} form $b_x\,$, induced by 
	\[
		b_x(y,z)=b(x,[z,y])=b([x,z],y)\mathfa y,z\in\ger g_1\ .
	\]
	Clearly, $b_x$ is non-degenerate on $\ger g_1/\ger z_{\ger g_1}(x)$. 

	We have $\ger g_1/\ger z_{\ger g_1}(x)=\ger k_1/\ger z_{\ger k_1}(x)\oplus\ger p_1/\ger z_{\ger p_1}(x)$, and both summands are $b_x$-totally isotropic, seeing that $b(\ger k_1,\ger p_1)=0$. As a well-known consequence of Witt's cancellation theorem, the dimension of $b_x$-totally isotropic subspaces does not exceed $\frac12\dim\ger g_1/\ger z_{\ger g_1}(x)$. Therefore,
	\[
		\dim\ger k_1-\dim\ger z_{\ger k_1}(x)=\dim\ger k_1/\ger z_{\ger k_1}(x)= \dim\ger p_1/\ger z_{\ger p_1}(x)=\dim\ger p_1-\dim\ger z_{\ger p_1}(x)
	\]
	which proves the assertion. 
\end{proof}

\begin{Prop}[iwasawa-infin]
	Let $(\ger g,\ger k,\theta)$ be a reductive symmetric superpair of even type, $\ger g_{0,\reals}$ a non-compact \emph{cs} form, $\ger a_\reals$ a real even Cartan subspace, $\ger a=\ger a_\reals\otimes\cplxs$, and $\ger n$ the nilpotent subalgebra for some positive system of $\Sigma(\ger g:\ger a)$. Then 
	\begin{equation}\label{eq:iwasawa-infin}
		\ger g=\ger k\oplus\ger a\oplus\ger n\nd\ger g_{0,\reals}=\ger k_{0,\reals}\oplus\ger a_\reals\oplus\ger n_{0,\reals}\ .
	\end{equation}
	We call these the \emph{Iwasawa decompositions} of $\ger g$ and $\ger g_{0,\reals}$, respectively.
\end{Prop}

\begin{proof}
	The subspace $\theta\ger n\cap\ger n\subset\ger g$ is $\ger a$-stable; if it were non-zero, it would contain a non-zero joint $\ger a$-eigenvector. Since this is impossible, $\theta\ger n\cap\ger n=0$. Hence, the intersection $\ger k\cap(\ger a\oplus\ger n)=\ger k\cap\ger n=0$, because $\ger k$ is $\theta$-fixed. Since $\ger g=\ger k\oplus\ger p$, the point is to show $\dim\ger a+\dim\ger n=\dim\ger p$. Now,
	\[
	\dim\ger k+\dim\ger p=\dim\ger g=\dim\ger m+\dim\ger a+2\dim\ger n\ .
	\]
	Thus, 
	\[
		2\cdot(\dim\ger a+\dim\ger n)=\dim\ger k-\dim\ger m+\dim\ger a+\dim\ger p=2\cdot\dim\ger p\ ,
	\]
	where the last equation follows from the remark in \thmref{Par}{possyst} and \thmref{Lem}{centdim-fmla}. Hence the assertion. 
\end{proof}

\subsection{The global Iwasawa decomposition}

We need to globalise the Iwasawa decomposition. This requires appropriate \emph{cs} supergroup pairs.

\begin{Def}
	Let $(\ger g,\ger k,\theta)$ be a symmetric superpair. A triple $(G_0,\ger g,\theta)$ where $(G_0,\ger g)$ is a \emph{cs} supergroup pair (\emph{cf.}~Appendix \ref{s:app-b}) where $G_0$ is \emph{connected} is called a \emph{global \emph{cs} form} of $(\ger g,\ger k,\theta)$ if the Lie algebra $\ger g_{0,\reals}$ of $G_0$ is a \emph{cs} form of $(\ger g,\ger k)$, and if $\theta$ is an involutive automorphism of $G_0$ (denoted by the same letter as the given involution on $\ger g$) whose differential is the restriction of $\theta$ to $\ger g_{0,\reals}$, \scth 
	\[
		\Ad(\theta(g))=\theta\circ\Ad(g)\circ\theta\in\End0{\ger g}\mathfa g\in G_0\ .
	\]
	
	A global \emph{cs} form $(G_0,\ger g,\theta)$ of $(\ger g,\ger k,\theta)$ is called \emph{non-compact} if $\ger g_{0,\reals}$ is a non-compact \emph{cs} form of $(\ger g,\ger k,\theta)$, and if $\Ad_\ger g(K_0)$ is compact, where $K_0$ denotes the analytic subgroup of $G_0$ generated by $\ger k_{0,\reals}$, and $\Ad_\ger g$ denotes the adjoint representation of $G_0$ on the Lie superalgebra $\ger g$. 
\end{Def}

\begin{Prop}[globalnccsformexists]
	Let $(\ger g,\ger k,\theta)$ be a strongly reductive symmetric superpair, and $\ger g_{0,\reals}$ a non-compact \emph{cs} form. Define $\ger g'=[\ger g,\ger g]$, $\ger g_{0,\reals}'=[\ger g_{0,\reals},\ger g_{0,\reals}]=\ger g_{0,\reals}\cap[\ger g_0,\ger g_0]$ and $Z_\ger g=\exp_{G_0}(\ger z(\ger g)\cap\ger g_{0,\reals})$. There exists a non-compact global \emph{cs} form $(G_0,\ger g,\theta)$ \scth the following conditions hold:
	\begin{enumerate}
		\item $\ger g_{0,\reals}$ is the Lie algebra of $G_0$; 
		\item the analytic subgroup $G'_0$ with Lie algebra $\ger g'\cap\ger g_{0,\reals}$ is closed;
		\item the analytic subgroup $G''_0$ with Lie algebra $\ger g_{0,\reals}'$ is closed; 
		\item $G'_0=Z(G'_0)\cdot G''_0$;
		\item $G_0=G'_0\times Z_\ger g$ is connected and simply connected;
		\item $Z(G_0)=Z(G'_0)\times Z_\ger g$. 
	\end{enumerate}
	We shall call $(G_0,\ger g,\theta,G'_0,G''_0,Z_\ger g)$ a \emph{standard} global \emph{cs} form of $(\ger g,\ger k)$. 
\end{Prop}

\begin{proof}
	Let $G'_0$ be the subgroup of $\GL(\ger g)$ generated by $e^{\ad(x)}$, $x\in\ger g_{0,\reals}$. By assumption, $\ger g_{0,\reals}$ contains a maximal compactly embedded subalgebra which is $\ad_\ger g$-compact. By \cite[Corollary 13.5.6 (e)]{hnbook-new}, $G'$ is closed in $\GL(\ger g)$. 
	
	The Lie algebra of $G'_0$ is $\ger g_{0,\reals}\cap\ger g'$ where $\ger g'=[\ger g,\ger g]$. Moreover, we may define $\theta$ on $G'_0$ by $\theta(g)=\theta\circ g\circ\theta$. If $G''_0$ is the analytic subgroup of $G'_0$ whose Lie algebra is $\ger g_{0,\reals}'=\ger g_{0,\reals}\cap[\ger g_0,\ger g_0]$, then $G''_0$ is closed in $G'_0$; indeed, it may be identified, under the restriction of even automorphisms on $\ger g$ to $\ger g_0$, with the adjoint group of $\ger g_{0,\reals}$ which is closed in $\GL(\ger g_0)$ because $\ger g_{0,\reals}'$ is semi-simple \cite[Corollary 13.5.7]{hnbook-new}. We remark that $G'_0=Z(G'_0)\cdot G''_0$.
	
	Let $Z_{\ger g}$ be the connected and simply connected real Lie group whose Lie algebra is $\ger z(\ger g)\cap\ger g_{0,\reals}$, \emph{i.e.}~the additive group $\ger z(\ger g)\cap\ger g_{0,\reals}$. Then $\theta$ is defined in an obvious way on $Z_\ger g$. If we set $G_0=G'_0\times Z_\ger g$, then $Z(G_0)=Z(G'_0)\times Z_\ger g$, and $\theta$ extends to $G_0$. Since $G'_0\subset\Aut0{\ger g}$ (Lie superalgebra automorphisms), it is clear how to define $\Ad:G_0\times\ger g\to\ger g$. 
	
	It is straightforward to check that all of the conditions are verified for the sextuple $(G_0,\ger g,\theta,G'_0,G''_0,Z_\ger g)$. To this end we remark only that by assumption, the analytic subgroup of $G'_0$ generated by $\ger g_{0,\reals}'\cap\ger k$ is compact.
\end{proof}

In the following, recall the functor $C$ defined in Appendix \ref{s:app-b}.

\begin{Prop}[iwasawa-global]
	Let $(\ger g,\ger k)$ be a strongly reductive symmetric superpair of even type, $(G_0,\ger g,\theta,G'_0,G''_0,Z_\ger g)$ a standard non-compact global \emph{cs} form, and $\ger a_\reals$ a real even Cartan subspace. Fix a positive system $\Sigma^+$ of $\Sigma$, and let $K_0$, $A=A_0$, $N_0$ be the analytic subgroups of $G_0$ with Lie algebras $\ger k_{0,\reals}$, $\ger a_\reals$ and $\ger n_{0,\reals}$, respectively. Let $K'_0=K_0\cap G'_0$ and $\ger k'=\ger k\cap\ger g'$ where $\ger g'=[\ger g,\ger g]$. Define \emph{cs} Lie supergroups $G=C(G_0,\ger g)$, $G'=C(G'_0,\ger g')$, $K=C(K_0,\ger k)$, $K'=C(K'_0,\ger k')$, and $N=C(N_0,\ger n)$. If $m$ denotes the multiplication of $G$, then the restriction of $m^{(2)}=m\circ(\id\times m)$ defines isomorphisms of \emph{cs} manifolds (\emph{cf.}~Appendix \ref{s:app-a})
	\[
		\phi:K\times A\times N\to G\nd\phi:K'\times(A\cap G'_0)\times N\to G'\ .
	\]
\end{Prop}

\begin{proof}
	Since $\ger n_{0,\reals}=[\ger a,\ger n_{0,\reals}]$, $N_0\subset G'_0$. Since $G'_0$ is connected and linear real reductive, we have the Iwasawa decomposition $G'_0=(K_0\cap G'_0)(A\cap G'_0)N_0$. Let $Z_\ger k=K_0\cap Z_\ger g$ and $Z_\ger p=\exp\ger p_{0,\reals}\cap Z_\ger g=A\cap Z_\ger g$. We have $G_0=G'_0\times Z_\ger g$, so it follows that $(k,a,n)\mapsto kan:K_0\times A_0\times N_0\to G_0$ is a diffeomorphism. 
	
	Let $\ger a=\ger a_\reals\otimes\cplxs$. From \thmref{Prop}{iwasawa-infin} and the Poincar\'e--Birkhoff--Witt theorem \cite[Chapter I, \S~3, Corollary~1 to Theorem~1]{scheunert-liesuperalgs}, it follows that the multiplication of $\Uenv0{\ger g}$ induces an isomorphism of super-vector spaces $\Uenv0{\ger k}\otimes\Uenv0{\ger a}\otimes\Uenv0{\ger n}\to\Uenv0{\ger g}$. Let $\Psi$ be the inverse map. 
	
	To see that the morphism $\phi$ in the statement of the theorem is an isomorphism of \emph{cs} manifolds, it suffices by \thmref{Cor}{globsectmor} to check that $\phi^*$ induces an isomorphism on the level the algebras of global sections. We have 
	\begin{align*}
		(\phi^*f)(b\otimes c\otimes d;k,a,n)&=f(\Ad((an)^{-1})(b)\Ad(n^{-1})(c)d;kan)\\
		&=f(\Ad((an)^{-1})(bc)d;kan)
	\end{align*}
	\fa $f\in\Gamma(\mathcal O_G)$, $b\in\Uenv0{\ger k}$, $c\in\Uenv0{\ger a}$, $d\in\Uenv0{\ger n}$, $k\in K_0$, $a\in A$, $n\in N_0$. We see that $\phi^{*-1}$ is given by the formula
	\[
		(\phi^{*-1}h)(u;kan)=h\Parens1{(\id\otimes\id\otimes\Ad((an)^{-1}))(\Psi(\Ad(an)(u)));k,a,n}\ .
	\]
	\fa $h\in\Gamma(\mathcal O_{K\times A\times N})$, $u\in\Uenv0{\ger g}$, $k\in K_0$, $a\in A$, $n\in N_0$. 
\end{proof}

The following notation will be used repeatedly.

\begin{Def}[iwasawa-proj]
	Let $\phi:K\times A\times N\to G$ be the Iwasawa isomorphism from \thmref{Prop}{iwasawa-global}. Similarly, let $\tilde\phi:N\times A\times K\to G$ be the isomorphism of \emph{cs} manifolds induced by multiplication in the opposite order. We define $A,H:G\to\ger a$ by $\exp\circ A\circ\tilde\phi=p_2$ and $\exp\circ H\circ\phi=p_2$ where $\exp:\ger a\to A$ is the exponential map $\ger a_{0,\reals}\to A$, considered as a morphism of \emph{cs} manifolds. Further, define $k,u:G\to K$ by $k\circ\phi=p_1$ and $u\circ\tilde\phi=p_3$, and $n:G\to N$ by $n\circ\tilde\phi=p_1$. 
	
	This is consistent with Helgason's notation \cite{helgason-gaga}; we feel that this supplies sufficient justification for indulging in the multiple uses of the letter $A$. The intended meaning will always be clear from the context. 
\end{Def}

\section{Integral formul\ae}\label{s:integral}

In what follows, let us assume given a strongly reductive symmetric superpair $(\ger g,\ger k)$ of even type, a standard non-compact \emph{cs} form $(G_0,\ger g,\theta,G'_0,G''_0,Z_\ger g)$, a real even Cartan subspace $\ger a_\reals$ with complexification $\ger a$, and a positive system $\Sigma^+$ of $\Sigma=\Sigma(\ger g:\ger a)$. We fix the notation from \thmref{Prop}{iwasawa-global}. Moreover, let $\ger m=\ger z_\ger k(\ger a)$, $M_0=Z_{K_0}(\ger a)$, and $M=C(M_0,\ger m)$. 

We will liberally use the notions introduced in Appendix \ref{s:app-c}. In particular, recall the notation $\Gamma_c$ for compactly supported global sections, the absolute Berezinian sheaf, the integral defined for its sections, the fibre integral (shriek map), invariant absolute Berezinians, and the concepts of geometric and analytic unimodularity.

\subsection{Integral formul\ae{} for the Iwasawa decomposition}

\begin{Par}
	Let $AN=NA=C(AN_0,\ger a\oplus\ger n)=C(N_0A,\ger a\oplus\ger n)$. This is a closed \emph{cs} subsupergroup of $G$. Similarly, define $MA$. 
\end{Par}

\begin{Lem}[gkunimod]
	The \emph{cs} Lie supergroups $G$, $G'$, $K$, $K'$, $M$, $A$, $N$ are geometrically and analytically unimodular. The \emph{cs} manifolds $G/K$, $G'/K'$, $G/MA$, and $G/A$ are geometrically and analytically unimodular as $G$-spaces. 
\end{Lem}

\begin{proof}
	All of the above statements follow by successive applications of \thmref{Prop}{unimodconds}. 
\end{proof}

\begin{Par}
	Let $\lambda\in\ger a^*$. We may define a function $e^\lambda\in\Ct[^\infty]0{\ger a_\reals,\cplxs}$ by setting $e^\lambda(x)=e^{\lambda(x)}$ \fa $x\in\ger a_\reals$. Define the linear form $\vrho=\tfrac12\str_{\ger n}\ad|_{\ger a}$. We have $\vrho=\vrho_0-\vrho_1$ where $2\vrho_j=\sum_{\lambda\in\Sigma_j^+}m_{\lambda,j}\cdot\lambda$, $m_{\lambda,j}=\dim_\cplxs\ger g_{j,\ger a}^\lambda$. Observe that $\vrho_0$ is real on $\ger a_\reals$. 
\end{Par}

\begin{Prop}[iwasawa-int]
	The invariant absolute Berezinians on $G$, $K$, $A$, and $N$ can be normalised \scth the following equations hold simultaneously: 
	\begin{align}
		&\int_{NA}f\,\Abs0{D(na)}=\int_{N\times A}m^*f\cdot p_2^*\log^*(e^{-2\vrho})\,\Abs0{Dn}\,da\ ,\quad f\in\Gamma_c(\mathcal O_{NA})\ ,\label{eq:na-fmla}\\
		&\int_{NA}f\Abs0{D(na)}=\int_{A\times N}m^*f\,da\,\Abs0{Dn}\ ,\quad f\in\Gamma_c(\mathcal O_{NA})\ ,\label{eq:an-fmla}\\
		&\int_Gf\,\Abs0{Dg}=\int_{N\times A\times K}\tilde\phi^*f\cdot p_2^*\log^*(e^{-2\vrho})\,\Abs0{Dn}\,da\,\Abs0{Dk}\ ,\quad f\in\Gamma_c(\mathcal O_G)\ ,\label{eq:nak-fmla}\\
		&\int_Gf\,\Abs0{Dg}=\int_{K\times A\times N}\phi^*f\cdot p_2^*\log^*(e^{2\vrho})\,\Abs0{Dk}\,da\,\Abs0{Dn}\ ,\quad f\in\Gamma_c(\mathcal O_G)\ ,\label{eq:kan-fmla}\\
		&\int_Gf\,\Abs0{Dg}=\int_{K\times N\times A}m^{(2)^*}f\,\Abs0{Dk}\,\Abs0{Dn}\,da\ ,\quad f\in\Gamma_c(\mathcal O_G)\ ,\label{eq:kna-fmla}
	\end{align}
\end{Prop}

\begin{proof}
	Equation \eqref{eq:na-fmla} follows directly from \thmref{Prop}{iwasawa-global}, \eqref{eq:prodsubgroup}, and \thmref{Lem}{berfactor}; \eqref{eq:an-fmla} follows in the same way, using the nilpotency of $\ger n$. Now \eqref{eq:nak-fmla}, and \eqref{eq:kna-fmla} follow in the same vein from \eqref{eq:na-fmla} (resp.~\eqref{eq:an-fmla}), \eqref{eq:unimodcond}, and the analytic unimodularity of $G/K$ in \thmref{Lem}{gkunimod}. The invariant absolute Berezinians can be normalised to give all equations simultaneously by the remark following \eqref{eq:unimodcond}. Finally, \eqref{eq:kan-fmla} follows from \eqref{eq:nak-fmla} by applying the invariance of $\Abs0{Dg}$ under $i^*$ (where $i$ is the inversion of $G$).
\end{proof}

\begin{Cor}
	For any $a\in A$ and $F\in\Gamma_c(\mathcal O_{NA})$, we have
	\begin{equation}\label{eq:aconjfactor}
		\int_N (R_a^*F)|_N\,\Abs0{Dn}=e^{2\vrho(a)}\cdot\int_N(L_a^*F)|_N\,\Abs0{Dn}\ .
	\end{equation}
\end{Cor}

\begin{proof}
	Considering $(R_a^*F)|_N\otimes\chi$ where $\chi\in\Cc0A$ with $\int_A\chi(a)\,da=1$, the result follows by comparing Equations \eqref{eq:an-fmla} and \eqref{eq:na-fmla}.  
\end{proof}

In the proof of the following lemma, we use the language of generalised points, \emph{cf}.~Appendix~\ref{app:s-points}. This method will be used repeatedly. 

\begin{Lem}[kgtransl]
	\Fa $f\in\Gamma_c(\mathcal O_K)$, we have
	\[
		p_{2!}\Parens1{(\id\times i)^*m^*u^*f\cdot(\Abs0{Dk}\otimes1)}=p_{2!}\Parens1{p_1^*fm^*A^*(e^{2\vrho})\cdot(\Abs0{Dk}\otimes1)}\ .
	\]
\end{Lem}

\begin{proof}
	Take $\chi\in\Gamma_c(\mathcal O_{N\times A})$ \scth $\int_{N\times A}\chi\cdot p_2^*\log^*\!e^{-2\vrho}\,\Abs0{Dn}\,da=1$, and define $h\in\Gamma_c(\mathcal O_G)$ by $\tilde\phi^*h=p_3^*\chi\cdot p_{12}^*f$. Then, since $G$ is analytically unimodular as a \emph{cs} Lie supergroup,
	\begin{align*}
		\int_K f\,\Abs0{Dk}&=\int_{N\times A\times K}\tilde\phi^*h\cdot p_2^*\log^*(e^{-2\vrho})\,\Abs0{Dn}\,da\,\Abs0{Dk}\\
		&=\int_Gh\,\Abs0{Dg}=p_{2!}\Parens1{m^*h\cdot(\Abs0{Dg}\otimes1)}=(*)
	\end{align*}
	
	Let $U$ be a \emph{cs} manifold. Take $U$-points $n\in_UN$, $a\in_UA$, $k\in_UK$, and $g\in_UG$. Then
	\[
		nakg=nan(kg)\exp A(kg)u(kg)=nan(kg)a^{-1}\cdot a\exp A(kg)\cdot u(kg)\ ,
	\]
	where the three factors are $U$-points of $N$, $A$, and $K$, respectively. Since $U$ and $n$, $a$, $k$, $g$ were arbitrary, by Yoneda's Lemma, there exist morphisms $f_1:A\times K\times G\to N$, $f_2:K\times G\to A$ \scth 
	\[
		m\cdot(\tilde\phi\times\id)=\tilde\phi\circ(m\circ(\id\times f_1),m\circ(p_2\times f_2),u\circ m\circ p_{34})\ .
	\]
	
	Therefore, $(*)$ equals 
	\begin{align*}
		p_{4!}\bigl(&(m\circ(\id\times f_1),m\circ(p_2\times f_2),u\circ m\circ p_{34}))^*\tilde\phi^*hp_2^*\log^*(e^{-2\vrho})(\Abs0{Dn}\,da\,\Abs0{Dk}\otimes1)\bigr)\\
		&=p_{4!}\Parens1{(p_{12},u\circ m\circ p_{34})^*\tilde\phi^*h\cdot p_2^*\log^*(e^{-2\vrho})\cdot p_{34}^*m^*A^*(e^{2\vrho})\cdot(\Abs0{Dn}\,da\,\Abs0{Dk}\otimes1)}\\
		&=p_{1!}\Parens1{m^*(u^*f\cdot A^*(e^{2\vrho}))\cdot(\Abs0{Dk}\otimes1)}\ .
	\end{align*}
	Here, we have used the fact that $N$ and $A$ are analytically unimodular \emph{cs} Lie supergroups. 
	
	To arrive at our claim, we need to `invert' $u^*$. Thus, let $U$ be a \emph{cs} manifold and $k\in_UK$, $g\in_UG$ be $U$-points. There exist unique $n\in_UN$, $a\in_UA$ \scth $nau(kg)=nakg$, which gives $u(u(kg)g^{-1})=k$. Hence, if $h(k)=f(u(kg))$ \fa $k\in_UK$, then $h(u(kg^{-1}))=f(k)$ \fa $k\in_UK$. Since $U$ and $g$ were arbitrary, this implies
	\[
		p_{2!}\Parens1{(\id\times i)^*m^*u^*h\cdot(\Abs0{Dk}\otimes1)}=p_{1!}\Parens1{p_1^*h\cdot m^*A^*(e^{2\vrho}))\cdot(\Abs0{Dk}\otimes1)}
	\]
	\fa $h\in\Gamma_c(\mathcal O_K)$. 
\end{proof}

\subsection{Joint eigen-superfunctions}

From now on, we will assume $\ger z(\ger g)=0$, so that $G=G'$, and $K=K'$ has a compact base.

In this section, we introduce a family of joint eigen-superfunctions for $\Uenv0{\ger g}^{\ger k}$, similar to the elementary spherical functions. It known \cite{zirnbauer-disk} that $\int_K1\,\Abs0{Dk}=0$ if $\dim\ger k_1\neq0$; thus, the generalisation of Harish-Chandra's Eisenstein integral yields superfunctions which are non-zero, but not obviously so.\footnote{In order to see that the Eisenstein integral is non-zero, one needs to study its asymptotic behaviour as $\lambda\to\infty$. This study will be the subject of a series of subsequent papers.}

If we are however willing to sacrifice $K$-biinvariance, we may introduce, into the Eisenstein integral, an auxiliary superfunction $\psi$ \scth $\int_K\psi\,\Abs0{Dk}=1$. As we will see presently, this defines a set of joint eigen-superfunctions which are obviously non-zero.

\begin{Par}
	For $D\in\Uenv0{\ger g}$, we define $D_{\ger a}\in\Uenv0{\ger a}$ by 
	\[
		D-D_{\ger a}\in\ger k\,\Uenv0{\ger g}+\Uenv0{\ger g}\ger n\ .
	\]
	Such a definition is possible due to the Poincar\'e--Birkhoff--Witt theorem. 
	
	The proof of the following lemma is standard, so we omit it. 
\end{Par}

\begin{Lem}[apart-hom]
	Let $D\in\Uenv0{\ger g}$, $D'\in\Uenv0{\ger g}^{\ger k}$. Then $(DD')_{\ger a}=D_{\ger a}D'_{\ger a}$. 
\end{Lem}

\begin{Par}
	For $D\in\Uenv0{\ger g}$, we define $\Gamma(D)=e^{-\vrho}D_{\ger a}e^{\vrho}\in S(\ger a)$. Recall that there is an algebra isomorphism $S(\ger a)\cong\cplxs[\ger a^*]:D\mapsto p$ defined by $p(\mu)=(De^\mu)(0)$. In these terms, $\Gamma(D)(\mu)=D_{\ger a}(\mu+\vrho)$. 
	
	By \thmref{Lem}{apart-hom}, have an algebra homomorphism $\Gamma:\Uenv0{\ger g}^{\ger k}\to S(\ger a)$, called the \emph{Harish-Chandra homomorphism}. (Not to be confused with the global sections functor $\Gamma$.) Obviously, we have $(\Uenv0{\ger g}\ger k)^\ger k\subset\ker\Gamma$. (The converse inclusion will be established below.)
\end{Par}

\begin{Par}
	In what follows, we fix $\psi\in\Gamma(\mathcal O_K)$ \scth $\int_K\psi\,\Abs0{Dk}=1$. Since $K$ is analytically unimodular as a \emph{cs} Lie supergroup, we also have $\int_Ki^*\psi\,\Abs0{Dk}=1$; hence, we may and will assume that $\psi=i^*\psi$. 
	
	Recall the notation from \thmref{Def}{iwasawa-proj}. Then $\tilde\psi:=k^*\psi\in\Gamma(\mathcal O_G)$ extends $\psi$. Slightly abusing notation, we will write $\psi=\tilde\psi$. Due to the Iwasawa decomposition, 
	\begin{equation}\label{eq:inttranspsi}
	p_{2!}\Parens1{m^*\psi\cdot(\Abs0{Dk}\otimes1)}=p_{2!}\Parens1{p_1^*\psi\cdot(\Abs0{Dk}\otimes1}=\int_K\psi\,\Abs0{Dk}=1\in\Gamma(\mathcal O_G)\ .
	\end{equation}
	
	For $f\in\Gamma(\mathcal O_G)$, we define 
	\[
	p_K(f)=f_K=p_{1!}\Parens1{m^*(\psi f)\cdot(1\otimes\Abs0{Dk})}\ .
	\]
	Then $f_K\in\Gamma(\mathcal O_G)^K=\Gamma(\mathcal O_{G/K})$ (\emph{cf.}~\thmref{Par}{quotient}). We let
	\begin{equation}\label{eq:psiphidefn}
		\phi^\psi_\mu=(H^*e^{\mu-\vrho})_K\mathfa \mu\in\ger a^*\ .
	\end{equation}
	Observe that $e^*\phi^\psi_\mu=1$, so that $\phi^\psi_\mu\neq0$. (Here, $e:*\to G$ is the unit of $G$.)

	In the following, we denote the $r$-action of $\Uenv0{\ger g}$ on $\Gamma(\mathcal O_G)$ (\emph{cf.}~Appendix \ref{s:app-b}) by juxtaposition, \emph{i.e.}~$Df$ instead of $r_Df$. 
\end{Par}

\begin{Prop}[phimu-jointeigenfn]
	For all $D\in\Uenv0{\ger g}^\ger k$, we have 
	\begin{equation}\label{eq:sph-eigen}
		D\phi^\psi_\lambda=D_\ger a(\lambda+\vrho)\cdot\phi^\psi_\lambda=\Gamma(D)(\lambda)\cdot\phi^\psi_\lambda\ .
	\end{equation}
	In particular, $\phi^\psi_\lambda$ is a joint eigenfunction of all $D\in\Uenv0{\ger g}^\ger k$. 
\end{Prop}

For the \emph{proof}, we first note the following lemma. 

\begin{Lem}[sphfn-afmla]
	We have
	\[
		\phi^\psi_\lambda=p_{2!}\Parens1{p_1^*\psi\cdot m^*A^*(e^{\lambda+\vrho})\cdot(\Abs0{Dk}\otimes1)}\ .
	\]
\end{Lem}

\begin{proof}
	First, we establish some identities for the morphisms $A$, $H$, $u$, using higher points. Thus, let $U$ be a \emph{cs} manifold and $k\in_UK$, $g\in_UG$ be $U$-points. There exists a unique $n\in_UN$ \scth $k(g)\exp H(g)n=g$. Then 
	\[
	g^{-1}=n^{-1}\exp(-H(g))k(g)^{-1}=n^{-1}\exp A(g)u(g)\ .
	\]
	Thus, 
	\begin{equation}\label{eq:Au-inv}
		H\circ i=-A\nd i\circ k=u\circ i\ .
	\end{equation}
	
	Next, let the unique $n\in_UN$ \scth $kg^{-1}=n\exp(A(kg^{-1}))u(kg^{-1})$. Then 
	\begin{align*}
		k&=n\exp(A(kg^{-1}))u(kg^{-1})g\\
		&=n\exp(A(kg^{-1}))n'\exp(A(u(kg^{-1})g))u(u(kg^{-1})g)\\
		&=n''\exp\Parens1{A(kg^{-1})+A(u(kg^{-1})g}k
	\end{align*}
	\fs $n'\in_UN$, where we write $n''=n\exp(A(kg^{-1}))n'\exp(-A(kg^{-1}))\in_UN$ and recall that $u(u(kg^{-1})g)=k$ (an identity which we have, in fact, rederived). 
	
	But this implies $A(kg^{-1})+A(u(kg^{-1})g)=0$. Since $U$, $k$, $g$ were arbitrary, we conclude that
	\begin{equation}\label{eq:Au}
		A\circ m\circ(u\circ m\circ (\id\times i),p_2)=-A\circ m\circ(\id\times i)
	\end{equation}
	as morphisms $K\times G\to\ger a$. 
	
	Now we can compute, using \eqref{eq:Au-inv}, $i^*\Abs0{Dk}=\Abs0{Dk}$, \eqref{eq:Au}, \thmref{Lem}{kgtransl} and the `inversion' of $u$ in its proof, 
	\begin{align*}
		\phi_\lambda^\psi&=p_{1!}\Parens1{m^*(\psi\cdot H^*e^{\lambda-\vrho})\cdot(1\otimes\Abs0{Dk})}\\
		&=p_{2!}\Parens1{\sigma^*m^*\psi\cdot (i\times i)^*m^*A^*e^{-\lambda+\vrho}\cdot(\Abs0{Dk}\otimes 1)}\\
		&=p_{2!}\Parens1{(i\times\id)^*\sigma^*m^*\psi\cdot(\id\times i)^*m^*A^*e^{-\lambda+\vrho}\cdot(\Abs0{Dk}\otimes 1)}\\
		&=p_{2!}\Parens1{(i\times\id)^*\sigma^*m^*\psi\cdot (u\circ m\circ(\id\times i),p_2)^*m^*A^*e^{\lambda-\vrho})\cdot(\Abs0{Dk}\otimes1)}\\
		&=p_{2!}\Parens1{(p_2,i\circ u\circ m)^*m^*\psi\cdot m^*A^*e^{\lambda+\vrho})\cdot(\Abs0{Dk}\otimes1)}
	\end{align*}
	Here, $\sigma=(p_2,p_1):K\times G\to G\times K$ denotes the flip.
	
	To complete the proof, we remark that $\psi=k^*(\psi|_K)$ and $i^*(\psi|_K)=\psi|_K$ by construction. Since $k(gu(kg)^{-1})=u(u(kg)g^{-1})^{-1}=k^{-1}$ \fa \emph{cs} manifolds $U$ and $k\in_UK$, $g\in_UG$, this implies $(p_2,i\circ u\circ m)^*m^*\psi=p_1^*\psi$, and, hence, our claim. 
\end{proof}

\begin{proof}[\protect{Proof of \thmref{Prop}{phimu-jointeigenfn}}]
	In view of \thmref{Lem}{sphfn-afmla}, it suffices to prove the equation $DA^*e^{\lambda+\vrho}=\Gamma(D)(\lambda)\cdot A^*e^{\lambda+\vrho}$, but this is trivial. 
\end{proof}

\subsection{The Harish-Chandra orbital integral} 

\begin{Def}
	The elements of the dense open subset of $\ger a$ given by the equation $\ger a'=\ger a\setminus\bigcup_{\lambda\in\Sigma}\lambda^{-1}(0)$ are called \emph{algebraically super-regular}. We also consider the sets $\ger a''=\ger a\setminus\bigcup_{\lambda\in\Sigma}\lambda^{-1}(2\pi i\ints)$ and $A'=\exp(\ger a_{\reals}\cap\ger a'')$. The elements of the latter set are called \emph{analytically super-regular}. 
	
	For $\lambda\in\Sigma$, $j=0,1$, let $m_{\lambda,j}$ be the multiplicity of $\lambda$ in the $\ger a$-module $\ger g_j$. For $a=\exp(h)\in A'$, the function 
	\begin{equation}\label{eq:weyldenom}
		D(a)=\frac{\prod_{\lambda\in\Sigma^+_0}\Abs0{\sinh\tfrac12\lambda(h)}^{m_{\lambda,0}}}{\prod_{\lambda\in\Sigma^+_1}(\sinh\tfrac12\lambda(h))^{m_{\lambda,1}}}=e^{\vrho(h)}\cdot\frac{\prod_{\lambda\in\Sigma_0^+}2^{m_{\lambda,0}}\Abs0{1-e^{-\lambda(h)}}^{m_{\lambda,0}}}{\prod_{\lambda\in\Sigma_1^+}2^{m_{\lambda,1}}(1-e^{-\lambda(h)})^{m_{\lambda,1}}}
	\end{equation}
	where we set $\sinh\lambda=\tfrac12(e^\lambda-e^{-\lambda})$, is well-defined and non-zero.
\end{Def}

\begin{Par}[conjdef]
	Let $c:G\times G\to G$ denotes \emph{conjugation} on $G$, \emph{i.e.}~$c=m\circ(m,i\circ p_1)$.

	We compute $c$ explicitly in terms of the \emph{cs} supergroup pair. To that end, fix an open subset $U\subset G_0$, $f\in\mathcal O_G(U)$, $u,v\in\Uenv0{\ger g}$, $(g,h)\in m^{-1}(U)$. We write $\Delta(u)=\sum_ju_j\otimes v_j$ and $\Delta(v)=\sum_iw_i\otimes z_i$. Then (\emph{cf.}~Appendix \ref{s:app-b})
	\begin{align*}
		c^*f(u&\otimes v;g,h)\\
		&=\sum\nolimits_{ji}(-1)^{\Abs0{v_j}\Abs0{w_i}}m^*f\Parens1{\Ad(h^{-1})(u_j)w_i\otimes\Ad(g)(S(v_j))\eps(z_i);gh,g^{-1}}\\
		&=\sum\nolimits_j(-1)^{\Abs0{v_j}\Abs0v}m^*f\Parens1{\Ad(h^{-1})(u_j)v\otimes\Ad(g)(S(v_j));gh,g^{-1}}\\
		&=\sum\nolimits_j(-1)^{\Abs0{v_j}\Abs0v}f\Parens1{\Ad(g)(\Ad(h^{-1})(u_j)vS(v_j));ghg^{-1}}\ .
	\end{align*}
	One may define
	\begin{equation}\label{eq:conjdef-infin}
		\begin{split}
			c_u^*f(v;g)&=(-1)^{\Abs0v\Abs0f}c^*f(u\otimes v;1,g)\\
			&=(-1)^{\Abs0v\Abs0f}\sum\nolimits_j(-1)^{\Abs0{v_j}\Abs0v}m^*f\Parens1{\Ad(g^{-1})(u_j)vS(v_j));g}\ .
		\end{split}
	\end{equation}
	Then $c_x^*=-\ell_x-r_x$ \fa $x\in\ger g$, by \eqref{eq:coind-infinaction} and \eqref{eq:gh-infinaction}. In particular, any subspace of $\Gamma(\mathcal O_G)$ invariant under $\ell_\ger g$ and $r_\ger g$ is invariant under $c_{\Uenv0{\ger g}}^*$, and \emph{vice versa}. 
	
	For $g\in G_0$, we set 
	\begin{equation}\label{eq:conjdef-grp}
		c_g^*f(v,h)=c^*f(1\otimes v;g,h)=f(\Ad(g)(v);ghg^{-1})\ .
	\end{equation}
	From \eqref{eq:coind-grpaction} and \eqref{eq:gh-grpaction}, one sees that $c_{g^{-1}}^*f=\ell_gr_gf$. Observe further that the pair $c_g=(h\mapsto ghg^{-1},c_g^*)$ is an automorphism of the \emph{cs} Lie supergroup $G$.
\end{Par}

\begin{Lem}[grpcomm-an]
	For $a\in A'$, $\xi_a=m\circ(c_{a^{-1}},i)$ is a \emph{cs} manifold automorphism of $N$.
\end{Lem}

\begin{proof}
	On the level of the underlying manifolds, $\xi$ is a diffeomorphism by \cite[Chapter I, \S~5, Lemma 6.4]{helgason-gaga}. We compute the tangent map $T_n\xi$. 
	
	To that end, observe that 
	\begin{gather*}
		(T_nc_{a^{-1}})(dL_n(x))=dL_{a^{-1}na}\Ad(a^{-1})(x)\ ,\\
		(T_ni)(dL_n(x))=dL_{n^{-1}}\Ad(n)(x)\ .
	\end{gather*}
	By \eqref{eq:tangentmult}, we find
	\begin{equation}\label{eq:grpcomm-an}
		(T_n\xi)(dL_n(x))=dL_{a^{-1}nan^{-1}}\Ad(n)\Parens1{\Ad(a^{-1})(x)-x}\ .
	\end{equation}
	By assumption, $\Ad(a^{-1})-\id:\ger n\to\ger n$ is an isomorphism, so the statement follows from \thmref{Prop}{invfnthm}. 
\end{proof}

\begin{Cor}[dfactor]
	\Fa $f\in\Gamma_c(\mathcal O_N)$, $a\in A'$, we have 
	\begin{equation}
		\int_Nf\,\Abs0{Dn}=\frac{\prod_{\lambda\in\Sigma_0^+}\Abs0{1-e^{-\lambda(\log a)}}^{m_{\lambda,0}}}{\prod_{\lambda\in\Sigma_1^+}(1-e^{-\lambda(\log a)})^{m_{\lambda,1}}}\cdot\int_N\xi_a^*f\,\Abs0{Dn}\ .
	\end{equation}
\end{Cor}

\begin{proof}
	By \thmref{Lem}{grpcomm-an}, $\xi_a$ is an isomorphism of \emph{cs} manifolds. By the invariance of absolute Berezin integrals under isomorphisms (\emph{cf.}~Appendix \ref{s:app-c}), the expression in \eqref{eq:grpcomm-an} of the tangent map of $\xi_a$, and the $N$-invariance of the absolute Berezinian $\Abs0{Dn}$ on $N$, everything comes down to the equation 
	\[
		\ABer[_{\ger n}]1{\Ad(n)\Parens1{\Ad(a^{-1})-1}}=\frac{\prod_{\lambda\in\Sigma_0^+}\Abs0{1-e^{-\lambda(\log a)}}^{m_{\lambda,0}}}{\prod_{\lambda\in\Sigma_1^+}(1-e^{-\lambda(\log a)})^{m_{\lambda,1}}}
	\]
	which is immediate.
\end{proof}

\begin{Par}[orbintprep]
	Let $\pi:G\to G/A$ be the canonical projection. Let $f\in\mathcal O_G(G_0)$ and $a\in A'$. Define $\tilde f_a\in\Gamma(\mathcal O_G)$ by $\tilde f_a=(\id\times a)^*c^*f$ where we consider $a$ as a morphism $*\to A$. This function is right $A$-invariant. In particular, there exists a unique $f_a\in\Gamma(\mathcal O_{G/A})$ \scth $\pi^*f_a=\tilde f_a$. 
\end{Par}

\begin{Prop}[orbintfmla]
	Let $a\in A'$. For any $f\in\Gamma_c(\mathcal O_G)$, we have $f_a\in\Gamma_c(\mathcal O_{G/A})$. Moreover, there is a normalisation of the invariant absolute Berezinian on $G/A$ (independent of $a$ and $f$), \scth
	\[
		D(a)\cdot\int_{G/A}f_a\,\Abs0{D\dot g}=e^{\vrho(a)}\cdot\int_{K\times N}L_{(1,a)}^*c^*f\,\Abs0{Dk}\,\Abs0{Dn}\ .
	\]
\end{Prop}

\begin{proof}
	That $f_a$ is compactly supported when considered as a superfunction on $G/A$ follows from \cite[Chapter I, \S 5, Proposition 5.6]{helgason-gaga}. Since $G/A$ is analytically unimodular, we deduce from \eqref{eq:kna-fmla} and \eqref{eq:fubfmla} that the invariant absolute Berezinian can be normalised \scth 
	\[
		\int_{G/A} h\,\Abs0{D\dot g}=2^{\dim\ger n_1-\dim\ger n_0}\cdot\int_{K\times N}m^*\pi^*h\,\Abs0{Dk}\,\Abs0{Dn}
	\]
	\fa $h\in\Gamma_c(\mathcal O_{G/A})$. Then, setting $C=2^{\dim\ger n_1-\dim\ger n_0}$, 
	\begin{align*}
		\int_{G/A}f_a\,\Abs0{D\dot g}&=C\cdot\int_{K\times N}m^*(\id,a)^*c^*f\,\Abs0{Dk}\,\Abs0{Dn}\\
		&=C\cdot\int_{K\times N}(\id\times\xi_a)^*L_{(1,a)}^*c^*f\,\Abs0{Dk}\,\Abs0{Dn}\\
		&=C\cdot\frac{\prod_{\lambda\in\Sigma_1^+}(1-e^{-\lambda(a)})^{m_{\lambda,1}}}{\prod_{\lambda\in\Sigma_0^+}\Abs0{1-e^{-\lambda(a)}}^{m_{\lambda,0}}}\cdot\int_{K\times N}L_{(1,a)}^*c^*f\,\Abs0{Dk}\,\Abs0{Dn}\ ,
	\end{align*}
	by \thmref{Cor}{dfactor}. The equation follows from \eqref{eq:weyldenom}. 
\end{proof}

\begin{Def}[orbint-def]
	Let $f\in\Gamma_c(\mathcal O_G)$. We define, \fa $a\in A$, the \emph{Harish-Chandra weighted orbital integral}
	\[
		F_f(a)=e^{\vrho(\log a)}\cdot\int_{K\times N}L_{(1,a)}^*c^*f\,\Abs0{Dk}\,\Abs0{Dn}\ .
	\]
	By \thmref{Prop}{orbintfmla}, we have \fa $a\in A'$, 
	\begin{equation}\label{eq:orbintfmla}
		F_f(a)=D(a)\cdot\int_{G/A}f_a\,\Abs0{D\dot g}\ .
	\end{equation}
\end{Def}

\section{The Harish-Chandra isomorphism}\label{s:hc}

Let us retain the assumptions and notation from Section \ref{s:integral}. We will at first generally assume that $\ger z(\ger g)=0$ (exceptions to this rule will be expressly stated). 

\subsection{Even Weyl group invariance}

\begin{Par}
	Recall the definition of the Harish-Chandra homomorphism $\Gamma$. Denote by $W_0=W(\ger g_0:\ger a)$ the even Weyl group.
\end{Par}

\begin{Prop}[hc-weylinv]
	We have $\Gamma(\Uenv0{\ger g}^{\ger k})\subset S(\ger a)^{W_0}$.
\end{Prop}

	The \emph{proof} requires a little preparation. The idea comes from Harish-Chandra's proof of the corresponding fact: One integrates the joint eigen-superfunctions $\phi_\mu^\psi$ against an arbitrary function, expresses this as the Abel transform of a weighted orbital integral, and uses the invariance of the latter. Compared to the even case, an additional complication is the lack of a well-behaved `invariant integral over $K$', and the ensuing occurrence of the auxiliary superfunction $\psi$. 

\begin{Prop}[orbint-weylinv]
	Let $f\in\Gamma_c(\mathcal O_G)$. Then 
	\[
		F_f(a^w)=F_f(a)\mathfa w\in W_0\,,\,a\in A\ .
	\]
	Here, $W_0=N_{K_0}(A)/Z_{K_0}(A)$, and $a^w=kak^{-1}$ for any $k\in N_{K_0}(A)$ \scth $w=kZ_{K_0}(A)$, and any $a\in A$. 
\end{Prop}

\begin{proof}
	By continuity, it suffices to check the equality on $A'$. Fix $a\in A'$ and $w=kZ_K(A)\in W_0$ where $k\in N_{K_0}(\ger a)$. Since $m_{\lambda,1}=\dim\ger g_{1,\ger a}^\lambda$ is even \fa $\lambda\in\Sigma$ \cite[Proposition 2.10 (v)]{ahz-chevalley}, $D(a^w)=D(a)$. We claim that 
	\[
		\int_{G/A}f_{a^w}\,\Abs0{D\dot g}=\int_{G/A}(c_k^*f)_a\,\Abs0{D\dot g}\ .\tag{$*$}
	\]
	It will then follow that $F_f(a^w)=F_{c_k^*f}(a)$, by \eqref{eq:orbintfmla}. Let us prove ($*$). 
	
	Clearly, $c_k$ induces a morphism $G/A\to G/A$. \Fa $h\in\Gamma_c(\mathcal O_G)$, define $h_A$ by the vector-valued integral $h_A=\int_AR_a^*h\,da$. Then 
	\[
	\int_{G/A}h_A\,\Abs0{D\dot g}=\int_Gh\,\Abs0{Dg}=\int_Gc_k^*h\,\Abs0{Dg}=\int_{G/A}(c_k^*h)_{MA}\,\Abs0{D\dot g}\ .
	\]
	On the other hand, the measure $da$ is $Z_{K_0}(A)$-conjugation invariant. This implies $(c_k^*h)_A=c_k^*(h_A)$, and we conclude that $c_k^*\Abs0{D\dot g}=\Abs0{D\dot g}$. Hence,
	\[
		\int_{G/A}f_{a^w}\,\Abs0{D\dot g}=\int_{G/A}c_k^*(f_{kak^{-1}})\,\Abs0{D\dot g}=\int_{G/A}(c_k^*f)_a\,\Abs0{D\dot g}\ ,
	\]
	so ($*$) holds, and we find that $F_f(a^w)=F_{c_k^*f}(a)$. On the other hand,
	\begin{align*}
		F_{c_k^*f}(a)&=e^{\vrho(\log a)}\cdot\int_{K\times N}L_{(k,a)}^*c^*f\,\Abs0{Dk}\,\Abs0{Dn}\\
		&=e^{\vrho(\log a)}\cdot\int_{K\times N}L_{(1,a)}^*c^*f\,\Abs0{Dk}\,\Abs0{Dn}=F_f(a)\ ,
	\end{align*}
	since $\Abs0{Dk}$ is $K$-invariant. This proves the claim.
\end{proof}

\begin{Lem}[abeltrans]
	For any $\lambda\in\ger a^*$ and $f\in\Gamma_c(\mathcal O_G)$, we have 
	\[
		\int_G\phi_\lambda^\psi\cdot f\,\Abs0{Dg}=\int_Ae^{\lambda(\log a)}F_f(a)\,da\ .
	\]
\end{Lem}

\begin{proof}
	We observe $\ABer[_\ger n]0{\Ad(a)}=e^{2\vrho(\log a)}$ (because the linear form $\vrho_0|_{\ger a_\reals}$ is real). Thus, we compute, using the analytic unimodularity of $G$, $\psi=k^*(\psi|_K)$, \eqref{eq:kan-fmla}, the analytic unimodularity of $K$, and $i^*(\psi|_K)=\psi|_K$, 
	\begin{align*}
		\int_G&\phi_\lambda^\psi f\,\Abs0{Dg}=\int_{G\times K}m^*(H^*(e^{\lambda-\vrho})\cdot\psi)\cdot p_1^*f\,\Abs0{Dg}\,\Abs0{Dk}\\
		&=\int_{G\times K}p_1^*H^*(e^{\lambda-\vrho})\cdot p_1^*k^*\psi\cdot (\id\times i)^*m^*f\,\Abs0{Dg}\,\Abs0{Dk}\\
		&=\int_Ae^{(\lambda+\vrho)(\log a)}\int_{K\times N\times K}p_1^*\psi\cdot L_{(1,a,1)}^*(\id\times\id\times i)^*m^{(2)*}f\,\Abs0{Dk_1}\,\Abs0{Dn}\,\Abs0{Dk_2}\,da\\
		&=\int_Ae^{(\lambda+\vrho)(\log a)}\int_{K\times N\times K}(i\circ p_3,p_1)^*m^*\psi\cdot L_{1,a}^*c^*f\,\Abs0{Dk_1}\,\Abs0{Dn}\,\Abs0{Dk_2}\\
		&=\int_Ae^{(\lambda+\vrho)(\log a)}\int_{K\times N}p_{1!}((i\times\id)^*m^*\psi\cdot(1\otimes\Abs0{Dk}))\cdot L_{1,a}^*c^*f\,\Abs0{Dk}\,\Abs0{Dn}\,da\ .
	\end{align*}
	
	Now, the invariance of $\Abs0{Dk}$ implies 
	\[
		p_{1!}((i\times\id)^*m^*\psi\cdot(1\otimes\Abs0{Dk}))=\int_K\psi\,\Abs0{Dk}=1\ .
	\]
	This proves the assertion.
\end{proof}

\begin{proof}[\protect{Proof of \thmref{Prop}{hc-weylinv}}]
	Let $w\in W_0$ and $\mu\in\ger a^*$. We have  
	\begin{align*}
		\int_G\phi_{w\lambda}^\psi\cdot f\,\Abs0{Dg}&=\int_Ae^{w\lambda(\log a)}\cdot F_f(a)\,da=\int_Ae^{\lambda(\log a)}F_f(a^w)\,da\\
		&=\int_Ae^{\lambda(\log a)}F_f(a)\,da=\int_G\phi_\lambda^\psi\cdot f\,\Abs0{Dg}
	\end{align*}
	\fa $f\in\Gamma_c(\mathcal O_G)$, by \thmref{Lem}{abeltrans} and \thmref{Prop}{orbint-weylinv}. 
	
	By \thmref{Lem}{berint-nondegenpair}, we conclude that $\phi_{w\lambda}^\psi=\phi_\lambda^\psi$. Then 
	\[
		\Gamma(D)(w\lambda)\cdot\phi^\psi_\lambda=\Gamma(D)(w\lambda)\cdot\phi^\psi_{w\lambda}=D\phi^\psi_{w\lambda}=D\phi^\psi_\lambda=\Gamma(D)(\lambda)\cdot\phi^\psi_\lambda\ ,
	\]
	by \eqref{eq:sph-eigen}. This proves the proposition, since $\phi_\lambda^\psi\neq0$. 
\end{proof}

\subsection{Odd Weyl group invariance}

\begin{Par}
	Recall the notation from \thmref{Par}{roots}. In addition, we let $\bar\Sigma_1^0=\Set1{\lambda\in\Sigma_1}{\Dual0\lambda\lambda=0}$ and $\bar\Sigma_1^1=\bar\Sigma_1\setminus\bar\Sigma_1^0$, where $\Dual0\cdot\cdot$ denotes the dual form of $b$ on $\ger a^*$. 
	
	Fix $\lambda\in\bar\Sigma_1$. We define $b^\theta(x,y)=b(x,\theta y)$ \fa $x,y\in\ger g$. The restriction of $b^\theta$ is a symplectic form on $\ger g_{1,\ger a}^{\pm\lambda}$ \cite[Proposition 2.10]{ahz-chevalley}, and this defines on $\ger g_{1,\ger a}^{\pm\lambda}$ the structure of symplectic $\ger m_0$-modules (where we recall that $\ger m=\ger z_\ger k(\ger a)$). Moreover, $\theta:\ger g_{1,\lambda}^\lambda\to\ger g_{1,\lambda}^{-\lambda}$ is an $\ger m_0$-equivariant symplectomorphism. Let $x_j,\tilde x_j\in\ger g_{1,\ger a}^\lambda$ be a $b^\theta$-symplectic basis, \emph{i.e.}
	\[
		b^\theta(x_i,x_j)=b^\theta(\tilde x_i,\tilde x_j)=0\ ,\ b^\theta(x_i,\tilde x_j)=2\delta_{ij}\ .
	\]
	Let $\ger g_{1,\ger a}^\lambda\oplus\ger g_{1,\ger a}^{-\lambda}=\ger k_1^\lambda\oplus\ger p_1^\lambda$ where $\ger k_1^\lambda\subset\ger k_1$ and $\ger p_1^\lambda\subset\ger p_1$. Set $x_j=y_j+z_j$ and $\tilde x_j=\tilde y_j+\tilde z_j$, according to this decomposition. Then 
	\begin{gather*}
		b(y_i,y_j)=b(\tilde y_i,\tilde y_j)=b(z_i,z_j)=b(\tilde z_i,\tilde z_j)=0\ ,\\
		b(y_i,z_j)=b(\tilde y_i,\tilde z_j)=b(y_i,\tilde z_j)=b(\tilde y_i,z_j)=0\ ,\\
		b(y_i,\tilde y_j)=b(\tilde z_j,z_i)=\delta_{ij}\ ,\\
		[h,y_i]=\lambda(h)z_i\ ,\ [h,\tilde y_i]=\lambda(h)\tilde z_i\ ,\ [h,z_i]=\lambda(h)y_i\ ,\ [h,\tilde z_i]=\lambda(h)\tilde y_i
	\end{gather*}
	\fa $i,j$ and $h\in\ger a$. Let $\ger m_0^\lambda=[\ger k_1^\lambda,\ger k_1^\lambda]$.
\end{Par}

\begin{Lem}[pureoddroot-relations]
	Retain the above notation, and define $A_\lambda\in\ger a$ by $b(A_\lambda,h)=\lambda(h)$ \fa $h\in\ger a$. We have the following equations.
	\begin{enumerate}
		\item \Fa $i,j$, 
		\[
			[y_i,z_j]=[\tilde y_i,\tilde z_j]=0\ ,\ [\tilde y_i,z_j]=-[y_i,\tilde z_j]=\delta_{ij}A_\lambda\ .
		\]
		\item \Fa $i,j$, 
			\[
			[y_i,y_j]=-[z_i,z_j]\ ,\ [\tilde y_i,\tilde y_j]=-[\tilde z_i,\tilde z_j]\ ,\ [y_i,\tilde y_j]=-[z_i,\tilde z_j]\ .
			\]
		\item\Fa $i,j,k$,\begin{gather*}
			[[y_i,y_j],y_k]=[[\tilde y_i,\tilde y_j],\tilde y_k]=0\ ,\\
			[[y_i,y_j],\tilde y_k]=\Dual0\lambda\lambda(\delta_{jk}y_i+\delta_{ik}y_j)\ ,\ 
			[[\tilde y_i,\tilde y_j],y_k]=-\Dual0\lambda\lambda(\delta_{jk}\tilde y_i+\delta_{ik}\tilde y_j)\ ,\\
			[[y_i,\tilde y_j],y_k]=-\Dual0\lambda\lambda\delta_{jk}y_i\ ,\ 
			[[y_i,\tilde y_j],\tilde y_k]=\Dual0\lambda\lambda\delta_{ik}\tilde y_j\ .
		\end{gather*}
		\item \Fa $i,j,k$, \begin{gather*}
			[[y_i,y_j],z_k]=[[\tilde y_i,\tilde y_j],\tilde z_k]=0\ ,\\
			[[y_i,y_j],\tilde z_k]=\Dual0\lambda\lambda(\delta_{jk}z_i+\delta_{ik}z_j)\ ,\ 
			[[\tilde y_i,\tilde y_j],z_k]=-\Dual0\lambda\lambda(\delta_{jk}\tilde z_i+\delta_{ik}\tilde z_j)\ ,\\
			[[y_i,\tilde y_j],z_k]=-\Dual0\lambda\lambda\delta_{jk}z_i\ ,\ 
			[[y_i,\tilde y_j],\tilde z_k]=\Dual0\lambda\lambda\delta_{ik}\tilde z_j\ .
		\end{gather*}
	\end{enumerate}
\end{Lem}

\begin{proof}
In the following, to simplify notation, we let
	\[
	(y,z),(y',z'),(y'',z''),(y''',z''')\in\Set1{(y_\ell,z_\ell),(\tilde y_\ell,\tilde z_\ell)}{\ell=1,\dotsc,\tfrac12m_{1,\lambda}}\ . 
	\]

	We have $[\ger g_{1,\ger a}^\lambda,\ger g_{1,\ger a}^\lambda]=0$. Applying this to $x_i,\tilde x_j$ and taking $\ger k$- and $\ger p$-projections gives the equations
	\[
		[y,y']=-[z,z']\ ,\ [y,z']=-[y',z]\ .
	\]
	This proves (ii), and reduces (i) to the cases $(y,y')\in\{(y_i,y_j),(y_i,\tilde y_j)\}$. Since $2\lambda\not\in\Sigma$, we have $[y,z']\in\ger p_0\cap(\ger a\oplus\ger g_\ger a^{\pm2\lambda})=\ger a$. On the other hand, for $h\in\ger a$, 
	\[
		b([y,z'],h)=-\lambda(h)b(y,y')=\begin{cases}0&(y,y')\in\{(y_i,y_j),(\tilde y_i,\tilde y_j)\}\ ,\\-\delta_{ij}b(A_\lambda,h)&(y,y')=(y_i,\tilde y_j)\ .\end{cases}
	\]
	The assertions in (i) follow by supersymmetry. 
	
	To prove (iii), we observe
	\[
		b([[y,y'],y''],y'''])=b([y,y'],[y'',y'''])=-b([y,y'],[z'',z'''])=-b([[y,y'],z''],z''')\ .
	\]
	Since $b(z''',z'')=b(y'',y''')$, (iii) follows from (iv). 
	
	Finally, 
	\[
		[[y,y'],z'']=[y,[y',z'']]+[y',[y,z'']]=:(*)\ .
	\]
	This is zero if $(y,y',z'')\in\{(y_i,y_j,z_k),(\tilde y_i,\tilde y_j,\tilde z_k)\}$, as follows from (i). Next, if $(y,y',z'')\in\{(y_i,y_j,\tilde z_k),(\tilde y_i,\tilde y_j,z_k)\}$, then, again by (i), 
	\[
		(*)=\mp(\delta_{jk}[y,A_\lambda]+\delta_{ik}[y',A_\lambda])=\pm\Dual0\lambda\lambda(\delta_{jk}z+\delta_{ik}z')
	\]
	where the sign $+$ occurs in the first case, and $-$ in the second. Finally, if we have $(y,y',z'')\in\{(y_i,\tilde y_j,z_k),(\tilde y_i,y_j,\tilde z_k)\}$, then 
	\[
		(*)=[y,[y',z'']]=\pm\delta_{jk}[y,A_\lambda]=\mp\Dual0\lambda\lambda\delta_{jk}z
	\]
	where the sign $-$ occurs in the first case, and $+$ in the second. This proves (iv), and as remarked, Assertion (iii) follows. 
\end{proof}

\begin{Cor}
	The graded subalgebra of $\ger g$ generated by $\ger g_1(\lambda):=\ger g_{1,\ger a}^\lambda\oplus\ger g_{1,\ger a}^{-\lambda}$ is exactly $\ger m_0^\lambda\oplus\cplxs A_\lambda\oplus\ger g_1(\lambda)$. It is invariant under $\theta$ and $\ger a$. In particular, $\ger m_0^\lambda\oplus\ger k_1^\lambda$ is a graded subalgebra of $\ger k$, and it leaves $\ger a\oplus\ger p_1^\lambda$ invariant. The radical of $b$ on $\ger m_0^\lambda\oplus\cplxs A_\lambda\oplus\ger g_1(\lambda)$ is $\ger m_0^\lambda\oplus\cplxs A_\lambda$ if $\lambda\in\bar\Sigma_1^0$, and $0$ if $\lambda\in\bar\Sigma_1^1$. 
\end{Cor}

\begin{proof}
	The last statement follows from part (iii) in \thmref{Lem}{pureoddroot-relations} and its proof. The other statements are immediate from the lemma. 
\end{proof}

\begin{Par}[ilambda-jlambda]
	To show that the image of $\Gamma$ is satisfies `odd Weyl group invariance', we adopt and adapt a technique due to Lepowsky, \emph{cf.}~\cite{gvbook}.
	
	Let $\ger n_\lambda=\ger g_{1,\ger a}^\lambda$, $\ger n_\lambda^\perp=\bigoplus_{\Sigma^+\ni\mu\neq\lambda}\ger g_\ger a^\mu$ and $\ger k_\lambda=\ger m_0^\lambda\oplus\ger k_1^\lambda$. If $\lambda\in\bar\Sigma_1^0$, we choose $h_0\in\ger a$ \scth $\lambda(h_0)=1$ and $b(h_0,h_0)=0$; in this case, we set $\ger a_\lambda=\Span0{h_0,A_\lambda}_\cplxs$. If, on the other hand, $\lambda\in\bar\Sigma_1^1$, then we let $\ger a_\lambda=\cplxs A_\lambda$. In any case, $\ger a_\lambda$ is $b$-non-degenerate, and $\ger m_\lambda:=\ger k_\lambda\oplus\ger a_\lambda\oplus\ger n_\lambda$ is a graded subalgebra of $\ger g$. Moreover, $\ger m_\lambda$ is $\theta$-invariant with $\ger p_\lambda:=\ger m_\lambda\cap\ger p=\ger a_\lambda\oplus\ger p_1^\lambda$. We write $\ger a_\lambda^\perp=\ger a\cap(\ger a_\lambda)^\perp$.
	
	We denote by $I_\lambda\subset S(\ger a)$ and $I_{\lambda,\ger m_\lambda}\subset S(\ger a_\lambda)$ the image of $S(\ger a\oplus\ger p_1^\lambda)^{\ger k_\lambda}$ and $S(\ger p_\lambda)^{\ger k_\lambda}$, respectively, under the projection onto $S(\ger a)$ along $\ger a^\perp S(\ger p)$. We observe that since $[\ger a_\lambda^\perp,\ger k_1^\lambda]=[\ger a_\lambda^\perp,\ger p_1^\lambda]=0$, we have $S(\ger p_\lambda)^{\ger k_\lambda}S(\ger a_\lambda^\perp)=S(\ger a\oplus\ger p_1^\lambda)^{\ger k_\lambda}$, so $I_\lambda=I_{\lambda,\ger m_\lambda}S(\ger a_\lambda^\perp)$. 
	
	Let $\Gamma_{\ger m_\lambda}$ denote the Harish-Chandra homomorphism for $\ger m_\lambda$. Whenever $\lambda\in\bar\Sigma_1^0$, then we let $J_\lambda=I_\lambda$ and $J_{\lambda,\ger m_\lambda}=I_{\lambda,\ger m_\lambda}$. For $\lambda\in\bar\Sigma_1^1$, we denote by $J_{\lambda,\ger m_\lambda}$ and $J_\lambda$ the sets $\Gamma_{\ger m_\lambda}(\Uenv0{\ger m_\lambda}^{\ger k_\lambda})$ and $J_{\lambda,\ger m_\lambda}S(\ger a_\lambda^\perp)$, respectively. Then we have  $J_\lambda=J_{\lambda,\ger m_\lambda}S(\ger a_\lambda^\perp)$ in any case. 
\end{Par}

\begin{Lem}[lepowsky-trick]
	Assume that $\Gamma_{\ger m_\lambda}(\Uenv0{\ger m_\lambda}^{\ger k_\lambda})\subset J_{\lambda,\ger m_\lambda}$. Then $\Gamma(\Uenv0{\ger g}^\ger k)\subset J_\lambda$. 
\end{Lem}

\begin{proof}
	For any $u\in\Uenv0{\ger g}$, there exists a unique $u_0\in\beta(S(\ger a\oplus\ger p_1^\lambda))$ \scth $u\equiv u_0\pmod{\ger n_\lambda^\perp\Uenv0{\ger g}+\Uenv0{\ger g}\ger k}$. (This follows from the Poincar\'e--Birkhoff--Witt theorem, applied to the vector space decomposition $\ger g=\ger a\oplus\ger p_1^\lambda\oplus\ger n_\lambda^\perp\oplus\ger k$.) In particular, if $u$ is $\ger k_\lambda$-invariant, then so is $u_0$. Moreover, 
	\[
	u_\ger a-(u_0)_\ger a=u_\ger a-u+u-u_0+u_0-(u_0)_\ger a\in\ger n\,\Uenv0{\ger g}+\Uenv0{\ger g}\ger k\ ,
	\]
	so $u_\ger a=(u_0)_\ger a$. 
	
	Let $u\in\Uenv0{\ger g}^\ger k$.  We have $S(\ger p_\lambda)^{\ger k_\lambda}S(\ger a_\lambda^\perp)=S(\ger a\oplus\ger p_1^\lambda)^{\ger k_\lambda}$ and $\beta(S(\ger p_\lambda)^{\ger k_\lambda})\subset\Uenv0{\ger m_\lambda}^{\ger k_\lambda}$. Hence, there exist $v_j\in\Uenv0{\ger m_\lambda}^{\ger k_\lambda}$ and $w_j\in\Uenv0{\ger a_\lambda^\perp}$ \scth $u_0=\sum_jw_jv_j$. For $\mu\in\ger a^*$, we have
	\[
		\Gamma(u)(\mu)=u_\ger a(\mu+\vrho)=(u_0)_{\ger a}\Parens1{\mu+\vrho}=\sum\nolimits_jw_j\Parens1{\mu+\vrho}(v_j)_{\ger a}\Parens1{\mu+\tfrac12m_{1,\lambda}\lambda}\ ,
	\]
	so we have proved our claim. 
\end{proof}

\begin{Lem}[ilambda-char]
	The set $I_\lambda$ is exactly the common domain of the differential operators $D$ with rational coefficients, whose local expression at any super-regular $\mu\in\ger a^*$ is exactly $\gamma_\mu(d)$, \fs fixed but arbitrary $d\in S(\ger p_1^\lambda)$. 
\end{Lem}

\begin{proof}
This is the main content of the discussion in \cite[Section 3.2]{ahz-chevalley}. 
\end{proof}

\begin{Lem}[rhoshift-invar]
	Let $\Dual0\lambda\lambda=0$. For any $z\in\cplxs$, the automorphism $p\mapsto p(\cdot+z\lambda)$ of $S(\ger a_\lambda)$ leaves $I_{\lambda,\ger m_\lambda}$ invariant. 
\end{Lem}

\begin{proof}
	For any fixed $p\in I_{\lambda,\ger m_\lambda}$, the statement is that for any $k$, the polynomial map $z\mapsto\partial_{A_\lambda}^k[p(\cdot+z\lambda)]+(A_\lambda^k):\cplxs\to S(\ger a_\lambda)/(A_\lambda^k)$ is zero. By degree considerations, it takes values in a finite-dimensional vector space. Any non-zero univariate vector-valued polynomial has but a finite set of zeros, so it suffices to prove that it vanishes at integer $z$, and by a trivial induction, for $z=1$. Thus, all we need to show is that $p\mapsto p(\cdot+\lambda)$ leaves $I_{\lambda,\ger m_\lambda}$ invariant. 
	
	We have $\ger a_\lambda=\Span0{h_0,A_\lambda}_\cplxs$, and $\partial_{A_\lambda}(h_0^kA_\lambda^\ell)=kh_0^{k-1}A_\lambda^\ell$. By \thmref{Lem}{ilambda-char} and \cite[Theorem 3.25]{ahz-chevalley}, $I_{\lambda,m_\lambda}$ is the common domain of the operators $A_\lambda^{-j}\partial_{A_\lambda}^j$ where $j=1,\dotsc,q:=\tfrac12m_{1,\lambda}$. Clearly, the latter is spanned by $h_0^kA_\lambda^\ell$ where $k\in\nats$ and $\ell\sge\min(k,q)$. The automorphism $p\mapsto p(\cdot+\lambda)$ maps $h_0^kA_\lambda^\ell$ to the element 
	\[
	(h_0+1)^kA_\lambda^\ell=\sum_{j=0}^k\binom kjh_0^jA_\lambda^\ell\ ,
	\]
	and therefore leaves $I_{\lambda,\ger m_\lambda}$ invariant. 
\end{proof}

\begin{Lem}[localinvar-isotropic]
	We have $\Gamma_{\ger m_\lambda}(\Uenv0{\ger m_\lambda}^{\ger k_\lambda})\subset I_{\lambda,\ger m_\lambda}$ for $\lambda\in\bar\Sigma_1^0$.
\end{Lem}

\begin{proof}
	We have $\Uenv0{\ger m_\lambda}=\beta(S(\ger p_\lambda))\oplus\Uenv0{\ger g}\ger k$, so in view of \thmref{Lem}{rhoshift-invar}, it suffices to prove that $\beta(p)_\ger a\in I_{\lambda,\ger m_\lambda}$ \fa $p\in S(\ger p_\lambda)^{\ger k_\lambda}$. 
	
	Let $k\in\nats$, $\ell\sge\min(k,q)$ where $q=\tfrac12m_{1,\lambda}$, and set 
	\[
		p_{k\ell}=\sum_{j=0}^{\min(k,q)}\binom kj\cdot h_0^{k-j}A_\lambda^{\ell-j}Z^j
	\]
	where $Z=z_1\tilde z_1+\dotsm+z_q\tilde z_q$. We have $\ad(y_n)(Z)=A_\lambda z_n$ and $\ad(\tilde y_n)(Z)=A_\lambda\tilde z_n$, so 
	\[
		\ad(y_n)(p_{k\ell})=\sum_{j=0}^{\min(k,q)}\Bracks2{-(k-j)\binom kjh_0^{k-j-1}A_\lambda^{\ell-j}z_nZ^j+j\binom kjh_0^{k-j}A_\lambda^{\ell-j+1}z_nZ^{j-1}}
	\]
	Since $z_nZ^q=0$, we have $\ad(y_n)(p_{k\ell})=0$. Similarly, $\ad(\tilde y_n)(p_{k\ell})=0$. By definition, $\ger k_\lambda$ is generated by $y_j$, $\tilde y_j$, so $p_{k\ell}\in S(\ger p_\lambda)^{\ger k_\lambda}$. Moreover, $p_{k\ell}(\mu)=(h_0^kA_\lambda^\ell)(\mu)$ \fa $\mu\in\ger a^*$, and this uniquely determines $p_{k\ell}$. It follows that $S(\ger p_\lambda)^{\ger k_\lambda}$ is spanned by the $p_{k\ell}$ where $k\sge0$ and $\ell\sge\min(k,q)$. 
	
	Next, we remark that $[\ger m_\lambda,\ger m_\lambda]\cap\ger p=\cplxs A_\lambda$. Hence, applying $\beta$ to $p_{k\ell}$ does not increase the $h_0$-degree. Since $A_\lambda$ is central in $\ger m_\lambda$, $\beta$ does also not decrease the $A_\lambda$-degree. It follows that $\beta(p_{k\ell})_\ger a\in I_{\lambda,\ger m_\lambda}$, and hence the assertion. 
\end{proof}

\begin{Par}
	For $\lambda\in\bar\Sigma_1^1$, we have defined $J_{\lambda,\ger m_\lambda}=\Gamma(\Uenv0{\ger m_\lambda}^{\ger k_\lambda})$. To prove our main result, we have to determine this set explicitly.  
	
	To that end, we will need a basic understanding of the relation between $\Gamma$ and the `restriction' map on $S(\ger p)$, and this will also be useful below. Thus, consider $\Uenv0{\ger g}$ with its standard filtration, which we denote $F_p\Uenv0{\ger g}$. The associated graded algebra $\gr\Uenv0{\ger g}$ has graded pieces $\gr_p\Uenv0{\ger g}=F_p\Uenv0{\ger g}/F_{p-1}\Uenv0{\ger g}$ and is supercommutative. The canonical map $\ger g\to F_1\Uenv0{\ger g}$ induces a canonical map $\ger g\to\gr_1\Uenv0{\ger g}\subset\gr\Uenv0{\ger g}$, which in turn gives a canonical map $S(\ger g)\to\gr\Uenv0{\ger g}$. 
	
	This map is an isomorphism, and we henceforth identify the algebras $S(\ger g)$ and $\gr\Uenv0{\ger g}$. Under this identification, it is known that $\beta:S(\ger g)\to\Uenv0{\ger g}$ equals the inverse of the canonical map $\Uenv0{\ger g}\to\gr\Uenv0{\ger g}$. Since $\Uenv0{\ger g}\ger k$ is a filtered subspace of $\Uenv0{\ger g}$, $S(\ger p)=\gr\Uenv0{\ger g}/\Uenv0{\ger g}\ger k$ as filtered super-vector spaces, and the isomorphism $\beta:S(\ger p)\to\Uenv0{\ger g}/\Uenv0{\ger g}\ger k$ is inverse to the canonical map $\Uenv0{\ger g}/\Uenv0{\ger g}\ger k\to\gr \Uenv0{\ger g}/\Uenv0{\ger g}\ger k$. 
\end{Par}

\begin{Lem}[grgamma]
	The map $\gr\Gamma:S(\ger p)\to S(\ger a)$ induced by $\Gamma:\Uenv0{\ger g}/\Uenv0{\ger g}\ger k\to\Uenv0{\ger a}$, considered as a map of filtered super-vector spaces, is the projection along $\ger a^\perp S(\ger p)$. 
\end{Lem}

\begin{proof}
	Denote the image of $p\in S(\ger p)$ in $S(\ger a)$ by $\bar p$. Since $1-\theta:\ger n\to\ger a^\perp$ is surjective, we have $\ger a^\perp\subset\ger k\oplus\ger n$, so that $p-\bar p\in(\ger n\oplus\ger k)S^{d-1}(\ger p)$ \fa $p\in S^d(\ger p)$.  
	
	The supersymmetrisation map $\beta:S(\ger a)\to\Uenv0{\ger a}$ is the identity, and on the level of $\ger g$, \fa $p_j\in S^{d_j}(\ger p)$, $j=1,2$, $\beta(p_1p_2)-\beta(p_1)\beta(p_2)\in\sum_{k<d_1+d_2}\beta(S^k(\ger p))$. It follows that 
	\[
	\beta(p)-\bar p\in\ger n\,\Uenv0{\ger g}+\Uenv0{\ger g}\ger k+\sum\nolimits_{k<d}\beta(S^k(\ger p))\mathfa p\in S^d(\ger p)\ .
	\]
	Since $\Gamma(\ger n\,\Uenv0{\ger g}+\Uenv0{\ger g}\ger k)=0$, we conclude 
	\begin{equation}\label{eq:hc-degfmla}
		\Gamma(\beta(p))-\bar p\in\bigoplus_{k<d}S^k(\ger a)\mathfa p\in S^d(p)\ .
	\end{equation}
	This proves the claim. 
\end{proof}

\begin{Par}
	Let $\lambda\in\bar\Sigma_1$, $c=\Dual0\lambda\lambda\neq0$, $2q=m_{1,\lambda}$. Choose some square root of $c$ and set
	\[
		a=c^{-1}A_\lambda\ ,\ w_j=c^{-1/2}z_j\ ,\ \tilde w_j=c^{-1/2}\tilde z_j\ ,\ v_j=c^{-1/2}y_j\ ,\ \tilde v_j=c^{-1/2}\tilde y_j\ .
	\]
	Then 
	\begin{gather*}
		[a,w_j]=v_j\ ,\ [a,v_j]=w_j\ ,\ [a,\tilde w_j]=\tilde v_j\ ,\ [a,\tilde v_j]=\tilde w_j\ ,\\
		[v_j,w_k]=[\tilde v_j,\tilde w_j]=0\ ,\ -[v_j,\tilde w_k]=[\tilde v_j,w_k]=\delta_{jk}a\ .
	\end{gather*}
	Let $\tilde J_{\lambda,\ger m_\lambda}$ be the subalgebra of $\cplxs[a]$ generated by $a^2-q^2$ and $a(a^2-q^2)^q$.
\end{Par}

\begin{Lem}
	Let $\lambda\in\bar\Sigma_1$, $\Dual0\lambda\lambda\neq0$. Then $\Gamma_{\ger m_\lambda}(\Uenv0{\ger m_\lambda}^{\ger k_\lambda})=\tilde J_{\lambda,\ger m_\lambda}$.
\end{Lem}

\begin{proof}	 
	We abbreviate $\Gamma_\lambda=\Gamma_{\ger m_\lambda}$. Since $\Uenv0{\ger m_\lambda}^{\ger k_\lambda}=\beta(S(\ger p_\lambda)^{\ger k_\lambda})\oplus(\Uenv0{\ger m_\lambda}\ger k_\lambda)^{\ger k_\lambda}$, it is sufficient to prove that $\Gamma_\lambda\circ\beta$ maps $S(\ger p_\lambda)^{\ger k_\lambda}$ onto $\tilde J_{\lambda,\ger m_\lambda}$. 
	
	To that end, let $W=w_1\tilde w_1+\dotsm+w_q\tilde w_q$. Then by \cite[4.7]{ahz-chevalley}, $S(\ger p^\lambda)^{\ger k_\lambda}$ is generated by the elements $P_2=a^2+2W$ and 
	\[
	P_{2q+1}=a^{2q+1}+\sum_{k=1}^q(-1)^{\frac12k(k+3)}a_{2q+1.k}a^{2(q-k)+1}W^k\ ,
	\]
	where
	\[
	a_{Nk}=\sum_{i=(k-N)_+}^{k-1}\Parens2{-\frac12}^iN\dotsm(N-k+i+1)\frac{(k-1+i)!}{(k-1-i)!i!}\ .
	\]
	In particular, 
	\[
		S^m(\ger p_\lambda)^{k_\lambda}=\begin{cases}
			\cplxs P_2^{2k}&m=2k,\\
			0&m=2k+1,k<q,\\ 
			\cplxs P_{2q+1}P_2^{k-q}&m=2k+1,k>q.
		\end{cases}
	\]
	
	In view of this statement, to see that $\Gamma_\lambda(\beta(S(\ger p_\lambda)^{\ger k_\lambda})\subset \tilde J_{\lambda,\ger m_\lambda}$, it suffices to prove the following: For all $k\in\nats$ and $\ell\sle k$, $\Gamma_\lambda\circ\beta$ maps $P_2^\ell$, and if $\ell\sge q$, then also $P_{2q+1}P_2^{\ell-q}$, to $\tilde J_{\lambda,\ger m_\lambda}$. We will prove this assertion by induction on $k$.
	
	First, we consider $L_2=\beta(P_2)$. Define $\tilde W=\tilde w_1w_1+\dotsm+\tilde w_qw_q$. Then we have $L_2=a^2+W-\tilde W$. Since 
	\[
	W+\tilde W=\sum\nolimits_{j=1}^q[w_j,\tilde w_j]\in\ger m_0^\lambda\subset\ger k_\lambda\ ,
	\]
	we have, modulo $\ger n\,\Uenv0{\ger g}+\Uenv0{\ger g}\ger k$, 
	\[
		W-\tilde W\equiv2W\equiv-\sum\nolimits_{j=1}^qv_j\tilde w_j\equiv-\sum\nolimits_{j=1}^q[v_j,\tilde w_j]=qa
	\]
	where we have used $u_j=v_j+w_j\in\ger n_\lambda$ and $[v_j,\tilde w_j]=-a$. Thus, 
	\[
	L_{2,\ger a}=a^2+2qa=(a+q)^2-q^2\ ,
	\]
	and since $\frac12\str_{\ger n_\lambda}\ad|_{\ger a_\lambda}=-\frac12 m_{1,\lambda}\lambda=-q\lambda$, 
	\[
	\Gamma_\lambda(L_2)=a^2-q^2\in \tilde J_{\lambda,\ger m_\lambda}\ .
	\]

	Now, assume we have proved the assertion for $k-1$, and we wish to prove it for $k$. We have $L_2^k\in\Uenv0{\ger m_\lambda}^{\ger k_\lambda}$, so $L_2^k\equiv\beta(p)\pmod{\Uenv0{\ger m_\lambda}\ger k_\lambda}$ \fs $p\in S(\ger p_\lambda)^{\ger k_\lambda}$. On the other hand, $\Gamma(L_2^k)=(a^2-q^2)^k$ and the $\ger a$-restriction of $P_2^k$ is $\bar P_2^k=a^{2k}$. 
	
	Now, the maps $\Gamma_\lambda\circ\beta$ and $p\mapsto\bar p$ coincide to leading order by \thmref{Lem}{grgamma}, so that $p'=P_2^q-p\in\bigoplus_{j<2k}S^j(\ger p_\lambda)$. The polynomial $p'$ being $\ger k_\lambda$-invariant, we obtain $\Gamma_\lambda(\beta(p'))\in \tilde J_{\lambda,\ger m_\lambda}$ by the inductive assumption. Since we have shown that $\Gamma_\lambda(\beta(p))=(a^2-q^2)^k\in \tilde J_{\lambda,\ger m_\lambda}$, it follows that $\Gamma_\lambda\circ\beta$ maps $P_2^k=p+p'$ to $\tilde J_{\lambda,\ger m_\lambda}$. 
	
	If $k\sge q$, then we need to show that $\Gamma_\lambda(P_{2q+1}P_2^{k-q})\in \tilde J_{\lambda,\ger m_\lambda}$, too. Similarly as above, we may take $p\in S(\ger p_\lambda)$ such that $aL_2^k\equiv\beta(p)\pmod{\Uenv0{\ger m_\lambda}^{\ger k_\lambda}}$, but we will need a little argument to prove that $p$ is $\ger k_\lambda$-invariant. 
	
	To that end, we compute for $(u,v,w)\in\{(u_j,v_j,w_j),(\tilde u_j,\tilde v_j,\tilde w_j)\}$, where we let $u_j=v_j+w_j$ and $\tilde u_j=\tilde v_j+\tilde w_j$,  
	\[
		[v,aL_2^q]=[v,a]L_2^q=-wL_2^q\equiv-uL_2^q+vL_2^q\equiv[v,L_2^q]=0\pmod{\ger n\,\Uenv0{\ger g}+\Uenv0{\ger g}\ger k}\ .
	\]
	
	Since $[v,aL_2^q]\equiv\beta([v,p])\pmod{\Uenv0{\ger g}\ger k}$, this implies $\beta([v,p])\in\ger n_\lambda\,\Uenv0{\ger a_\lambda\oplus\ger n_\lambda}$. Observe that $\theta(aL_2^k)=-aL_2^k$ since $\theta(L_2)=L_2$. Because $\theta(v)=v$, we find
	\[
		-\beta([v,p])]=\theta(\beta([v,p]))\in\ger n\,\Uenv0{\ger a\oplus\ger n}\cap\bar{\ger n}\,\Uenv0{\ger a\oplus\bar{\ger n}}=0\ . 
	\]
	Because $\beta$ is injective, this shows $[v,p]=0$, and it follows that $p$ is $\ger k_\lambda$-invariant. 
	
	Since $L_2^k$ is $\ger k_\lambda$-invariant, \thmref{Lem}{apart-hom} implies 
	\[
	\Gamma_\lambda(aL_2^k)=\Gamma_\lambda(a)\Gamma_\lambda(L_2^k)=(a-q)(a^2-q^2)^k\in \tilde J_{\lambda,\ger m_\lambda}\ .
	\]
	Of course, $\bar P_{2q+1}\bar P_2^{q-k}=a^{2k+1}$, and the same argument as above allows us to conclude by induction that $\Gamma_\lambda(\beta(P_{2q+1}P_2^{q-k}))\in \tilde J_{\lambda,\ger m_\lambda}$. 
	
	Moreover, if $k=q$, then $\Gamma_\lambda(\beta(P_{2q+1}))$ and $\Gamma_\lambda(aL_2^q)=(a-q)(a^2-q^2)^q$ coincide to leading order. Thus, if $p\in S(\ger p_\lambda)^{\ger k_\lambda}$ is such that $aL_2^k\equiv\beta(p)\pmod{\Uenv0{\ger g}\ger k}$, then the difference $P_{2q+1}-p$ is, by the inductive assumption, a polynomial in $P_2$. Since $\theta(P_{2q+1})=-P_{2q+1}$, $\theta(aL_2^q)=-aL_2^q$, $\theta(P_2)=P_2$, and $\beta$ is $\theta$-equivariant, we find that $p=P_{2q+1}$, so that $\Gamma(\beta(P_{2q+1}))=\Gamma_\lambda(aL_2^q)=(a-q)(a^2-q^2)^q$. 
	
	This finally proves our claim that $\Gamma_\lambda(\Uenv0{\ger m_\lambda}^{\ger k_\lambda})\subset \tilde J_{\lambda,\ger m_\lambda}$. Since we have the identities $\Gamma_\lambda(\beta(P_2))=a^2-q^2$ and $\Gamma_\lambda(\beta(P_{2q+1}))=(a-q)(a^2-q^2)^q$, we have also proved the equality $\Gamma_\lambda(\Uenv0{\ger m_\lambda}^{\ger k_\lambda})=\tilde J_{\lambda,\ger m_\lambda}$.
\end{proof}

\begin{Cor}[jlambda-explicit]
	Let $\lambda\in\bar \Sigma_1^1$. Then $J_{\lambda,\ger m_\lambda}=\cplxs[a^2-q^2,(a-q)(a^2-q^2)^q]$ where $\Dual0\lambda\lambda\cdot a=A_\lambda$ and $2q=m_{1,\lambda}$. 
\end{Cor}

\subsection{The Harish-Chandra isomorphism}

\begin{Par}
	We can now finally state and prove our main result. To that end, recall some notation. Let $(\ger g,\ger k,\theta)$ be a strongly reductive symmetric superpair of even type, and $\ger a$ an even Cartan subspace, giving rise to the set $\Sigma$ of restricted roots. Denote by $W_0$ the Weyl group of $\Sigma_0$ and fix a positive system $\Sigma^+$. Let 
	\[
		I(\ger a)=S(\ger a)\cap\bigcap_{\lambda\in\bar\Sigma_1}I_\lambda\nd J(\ger a)=S(\ger a)^{W_0}\cap\bigcap_{\lambda\in\bar\Sigma_1}J_\lambda
	\]
	where $I_\lambda$ and $J_\lambda$ are defined in \thmref{Par}{ilambda-jlambda}. Let $\Gamma:\Uenv0{\ger g}^\ger k\to S(\ger a)$ be the Harish-Chandra homomorphism.
\end{Par}

\begin{Th}[hc]
	The image of $\Gamma$ is $J(\ger a)$, and its kernel is $(\Uenv0{\ger g}\ger k)^\ger k$. Therefore, it induces an algebra isomorphism $D(X)^G\to J(\ger a)$ for any global \emph{cs} form $(G,K)$ of $(\ger g,\ger k)$, where $D(X)^G$ is the set $G$-invariant differential operators on $X=G/K$. 
\end{Th}

\begin{proof}
	Since $\Uenv0{\ger g}^\ger k=S(\ger z(\ger g))\otimes\Uenv0{\ger g'}^{\ger k'}$ where $\ger g'=[\ger g,\ger g]$ and $\ger k'=\ger g'\cap\ger k$, there is no restriction in assuming that $\ger z(\ger g)=0$. By the existence of a standard non-compact global \emph{cs} form (\thmref{Prop}{globalnccsformexists}), we may apply \thmref{Prop}{hc-weylinv}. By \thmref{Lem}{localinvar-isotropic} and  \thmref{Lem}{lepowsky-trick}, the image of $\Gamma$ is contained in $J(\ger a)$. (Note that the assumption of the latter lemma is trivially verified for $\lambda\in\bar\Sigma_1^1$.)
	
	Define $I_0=S(\ger a)^{W_0}\cap\bigcap_{\lambda\in\bar\Sigma_1^0}I_\lambda$. For $\lambda\in\bar\Sigma_1^1$, we have, by \thmref{Par}{ilambda-jlambda}, \cite[Proposition 4.5]{ahz-chevalley} and \thmref{Cor}{jlambda-explicit}, 
	\[
		I_\lambda=\cplxs[a^2,a^{2q+1}]S(A_\lambda^\perp)\nd J_\lambda=\cplxs[a^2-q^2,(a-q)(a^2-q^2)^q]S(A_\lambda^\perp)
	\]
	where $\Dual0\lambda\lambda\cdot a=A_\lambda$ and $2q=m_{1,\lambda}$. 
	
	Then $I_0$ and $I_\lambda$ are $\ints$-graded subspaces of $S(\ger a)$, whereas $J_\lambda$ is a filtered subspace with $\gr J_\lambda=I_\lambda$. We have $I(\ger a)=I_0\cap\bigcap_{\lambda\in\bar\Sigma_1^1}I_\lambda$ and $J(\ger a)=I_0\cap\bigcap_{\lambda\in\bar\Sigma_1^1}J_\lambda$. For any filtered subspace $V$ of $S(\ger a)$, $\gr V$ injects into $\gr S(\ger a)=S(\ger a)$. In this sense, one has 
	\[
	\gr J(\ger a)=I_0\cap\gr\bigcap_{\lambda\in\bar\Sigma_1^1}J_\lambda\subset I_0\cap\bigcap_{\lambda\in\bar\Sigma_1^1}\gr J_\lambda=I_0\cap\bigcap_{\lambda\in\bar\Sigma_1^1}I_\lambda=I(\ger a)\ .
	\]
	
	On the other hand, $\gr\Gamma$ is the `restriction' map $S(\ger p)\to S(\ger a)$, so \cite[Theorem 3.25]{ahz-chevalley} implies $I(\ger a)=\mathrm{im}\gr\Gamma\subset\gr J(\ger a)$, \emph{i.e.}~$I(\ger a)=\gr J(\ger a)$.
	
	Consider the filtered complex
	\begin{align*}
			C\colon&\qquad\xymatrix{0\ar[r]&(\Uenv0{\ger g}\ger k)^\ger k\ar[r]&\Uenv0{\ger g}^\ger k\ar[r]^-{\Gamma}&J(\ger a)\ar[r]&0\ .}
	\intertext{We wish to see that this is an exact sequence, \emph{i.e.}~$H(C)=0$. Since the filtration on $C$ is bounded below and exhaustive, the spectral sequence $(E^r)$ of the filtration converges to $H(C)$ \cite[Theorem 5.5.1]{weibel-homolalg}. By \thmref{Lem}{grgamma}, }
			E^0=\gr C\colon&\qquad\xymatrix{0\ar[r]&(S(\ger g)\ger k)^\ger k\ar[r]&S(\ger g)^\ger k\ar[r]^-{R}&I(\ger a)\ar[r]&0}
	\end{align*}
	where $R(S(\ger g)\ger k)=0$ and on $S(\ger p)$, $R$ is the `restriction' map $p\mapsto\bar p$. (Observe that $S(\ger g)\ger k\subset S(\ger g)=\gr \Uenv0{\ger g}$ identifies with $\gr(\Uenv0{\ger g}\ger k)$.)
	 
	By \cite[Theorem 3.25]{ahz-chevalley}, $R$ is injective on $S(\ger p)^\ger k$, with image $I(\ger a)$. Thus, $E^1=H(\gr C)=0$, and this proves the theorem. 
\end{proof}

\renewcommand\thesection{A}
\section{Appendix: The category of \emph{cs} manifolds}\label{s:app-a}

\begin{Par}
	Unlike their ungraded counterparts, simple complex Lie superalgebras do not in general possess real forms whose even part is compact \cite{serganova-symmetric}. At first sight, this seems to make the generalisation of many aspects of harmonic analysis to the graded setting unfeasible. On the other hand, it is common in the physics community to work with \emph{real} `variables' and \emph{complex} `Grassmann variables', \emph{i.e.}~to choose a real form only of the even part. 
	
	Taken seriously, this leads to the category of \emph{cs} manifolds introduced by J.~Bernstein \cite{bernstein-qft,deligne-morgan}.\footnote{The abbreviation \emph{cs} stands for `complex super', but \emph{cs} manifolds are \emph{not} the same thing as complex supermanifolds.} These are complex super-ringed spaces, locally isomorphic to \emph{real} superdomains with \emph{complexified} structure sheaves. 
	
	For ordinary real manifolds, complexifying the structure sheaves does not change anything. That is, one obtains a fully faithful embedding of real manifolds into \emph{cs} manifolds. However, for supermanifolds, matters do change, \emph{i.e.} the embedding does not extend fully faithfully to real supermanifolds. However, this should not be seen as a defect, since it is precisely this fact that resolves some of the issues related to real structures in the super world. The following two reasons make \emph{cs} manifolds a useful tool for our purposes.
	
	First, the usual theory of Berezin integration on real supermanifolds does not require the choice of an orientation in the odd coordinate directions, in particular, it does not require a real structure. This allows for the transposition of this theory to the category of \emph{cs} manifolds, as was observed by J.~Bernstein \cite{bernstein-qft,deligne-morgan}. 
	
	Second, the group objects in the category of \emph{cs} manifolds (which we call \emph{\emph{cs} Lie supergroups}) can be described by linear actions of \emph{real} Lie groups $G$ on \emph{complex} Lie superalgebras satisfying suitable compatibility assumptions. This allowed, in Section \ref{s:iwasawa}, for the definition of compact and non-compact `\emph{cs} forms'. 
	
	The theory of \emph{cs} manifolds seems to have been largely ignored by the mathematical community, perhaps because its utility has not been fully appreciated. For this reason, we find it appropriate to briefly review the foundations. Most arguments will carry over from the theory of real supermanifolds, but there are a few subtle points. Let us mention two: the theory of linear \emph{cs} manifolds, and the theory of \emph{cs} Lie supergroups. In both cases, real structures on the even parts intervene in a non-trivial fashion. 
\end{Par}

\subsection{Basic theory of \emph{cs} manifolds}

We will use notions of `super' mathe\-matics. We refer the reader to \cite{leites,kostant-supergeom,deligne-morgan}.

\begin{Def}
	A \emph{complex super-ringed space} is a pair $X=(X_0,\mathcal O)$ where $X_0$ is a topological space and $\mathcal O=\mathcal O_X$ is a sheaf of \emph{complex} super-com\-mu\-ta\-tive super\-al\-gebras. We define a subsheaf $\mathcal N\subset\mathcal O$ as $\mathcal N(U)=\mathcal O(U)\cdot\mathcal O(U)_1$, the ideal generated by the odd elements, \fa open $U\subset X_0$. We denote $\pi_X:\mathcal O\to\mathcal O/\mathcal N$ the canonical sheaf morphism.
	
	A \emph{morphism $\vphi:X\to Y$ of complex super-ringed spaces} is a tuple $(f,f^*)$ where $f:X_0\to Y_0$ is a continuous function and $f^*:\mathcal O_Y\to f_*\mathcal O_X$ is a morphism of sheaves of complex superalgebras (in particular, even and unital). One obtains a well-defined category of complex super-ringed spaces. Occasionally, when $f:X\to Y$ is a morphism of complex super-ringed spaces, we shall write $f:X_0\to Y_0$ for the underlying continuous map, and $f^*:\mathcal O_Y\to f_*\mathcal O_X$ for the map of sheaves, thus slightly abusing the notation. 
	
	There is a natural functor from real super-ringed spaces to complex super-ringed spaces, given by the complexification of the structure sheaves. We shall denote its application by the subscript $_\cplxs$. This functor is neither full nor faithful. By complexification, we associate with the real supermanifolds $\reals^{p|q}=(\reals^p,\mathcal C^\infty_{\reals^p}\otimes\bigwedge(\reals^q)^*)$ the complex super-ringed spaces $\reals^{p|q}_\cplxs=(\reals^p,\mathcal C^\infty_{\reals^p}\otimes\bigwedge(\cplxs^q)^*)$. 
	
	If $X=(X_0,\mathcal O)$ is a complex super-ringed space, then an \emph{open subspace} $U$ is a complex super-ringed space of the form $(U_0,\mathcal O|_{U_0})$, for some open subset $U_0\subset X_0$. A \emph{\emph{cs} domain} of graded dimension $p|q$ (where $p,q\in\nats$) is a complex super-ringed space isomorphic to an open subspace of $\reals^{p|q}_\cplxs$; the \emph{cs} domain whose underlying set is $\vvoid$ may be given any graded dimension. 
	
	Finally, a complex super-ringed space $X=(X_0,\mathcal O)$ is called a \emph{\emph{cs} manifold} if $X_0$ is Hausdorff and second countable, and $X$ possesses a cover by open subspaces which are \emph{cs} domains (whose graded dimensions may vary). Such open subspaces (and their underlying open sets) are called \emph{coordinate neighbourhoods}. If the open cover may be chosen such that the coordinate neighbourhoods all possess the same graded dimension, then we say that $X$ is of (pure) graded dimension $p|q$. One obtains a category of \emph{cs} manifolds as a full subcategory of complex super-ringed spaces. When no confusion is possible, we simply say that $X$ is a \emph{cs} manifold, and denote its structure sheaf by $\mathcal O_X$ (and similarly for other letters of the alphabet). We will the global sections functor by $\Gamma$.
	
	We apologise to the reader for the somewhat unfortunate parlance `\emph{cs} manifold'; we point out that besides being extant terminology, obvious replacements such as `semi-real supermanifold' would probably not constitute an improvement on the already rather laboured `super' nomenclature. Moreover, although this is clearly not a consistent use of the abbreviation `\emph{cs}', we will often explicitly add the prefix `super' in derived terminology (\emph{viz.}~`\emph{cs} Lie supergroup'), since contracting this prefix seems to hide the true nature of the objects considered. 
\end{Def}

\begin{Par}
	Let $X=(X_0,\mathcal O)$ be a \emph{cs} manifold. If $U\subset X_0$ is open, $f\in\mathcal O(U)$, and $x\in U$, then we define $f(x)$ to be the unique complex number $\alpha\in\cplxs$ \scth $(f-\alpha)_x$ is not invertible---here, $h_x$ denotes the \emph{germ} of $h$ at $x$. Such a complex number exists uniquely, because $\mathcal O_x$ is isomorphic to $\mathcal C_{V,y}^\infty\otimes\bigwedge(\cplxs^q)^*$ for some $p,q$, some open $V\subset\reals^p$, and some $y\in V$, and is hence a local algebra. We write $\tilde f:U\to\cplxs$ for the function $x\mapsto f(x)$. 
	
	We say that $f\in\mathcal O(U)$ takes \emph{values} in any given set $A\subset\cplxs$ whenever the function $\tilde f$ does. Hence, $f$ may take real values, \emph{etc}. Moreover, there exists a \emph{canonical} antilinear involution on $\mathcal C^\infty_{X_0,\cplxs}=\mathcal O/\mathcal N$. We denote by $\mathcal C^\infty_{X_0}$ the sheaf of fixed points for this involution, given by the elements that take real values. 
	
	Let $(f,f^*):X\to Y$ be a morphism of \emph{cs} manifolds, and $h\in\mathcal O(V)$ where $V\subset Y_0$ is open. Let $x\in f^{-1}(V)$. Since $(h-h(f(x)))_{f(x)}$ is not invertible, nor is $(f^*h-h(f(x)))_x=(f^*h)_x-h(f(x))$. By the definition of $f^*h(x)$, it follows that $f^*h(x)=h(f(x))$, \emph{i.e.}~$\widetilde{f^*h}=\tilde h\circ f$. In particular, any morphism of \emph{cs} manifolds is local, and the induced morphism  $\mathcal C^\infty_{Y_0,\cplxs}\to f_*\mathcal C^\infty_{X_0,\cplxs}$ respects the canonical antilinear involutions; so it gives a morphism $(X_0,\mathcal C^\infty_{X_0})\to(Y_0,\mathcal C^\infty_{Y_0})$ of \emph{real} ringed spaces. This has the following two consequences.
\end{Par}

\begin{Prop}
	For any \emph{cs} manifold $X=(X_0,\mathcal O)$, there exists on $X_0$ a unique structure of smooth manifold for which the morphism $\mathcal C^\infty_{X_0}\to\mathcal C^{\vphantom\infty}_{X_0}$ which sends $h\pmod{\mathcal N}$ to the continuous function $\tilde h$ defines an isomorphism of $\mathcal C^\infty_{X_0}$ with the sheaf of real-valued smooth functions on $X_0$. 
\end{Prop}

\begin{Prop}
	The complexification functor from the category of real super-ringed spaces to the category complex super-ringed spaces restricts to a fully faithful embedding of smooth real manifolds into \emph{cs} manifolds. 
\end{Prop}

\begin{Par}
	Let $X=(X_0,\mathcal O)$ be a \emph{cs} manifold. Although $\mathcal O$ does not have a canonical $\mathcal C^\infty_{X_0}$-module structure, we can use the sheaf $\mathcal C^\infty_{X_0}$ to show that $\mathcal O$ is a fine and, in particular, a soft sheaf (see below). Let us remark that \emph{a posteriori}, the choice of a partition of unity in $\mathcal O$ subordinate to some locally finite cover by coordinate neighbourhoods can be used to define a $\mathcal C^\infty_{X_0}$-module structure, but this construction is non-canonical. 
\end{Par}

\begin{Prop}[strucsheaffine]
	Let $X=(X_0,\mathcal O)$ be a \emph{cs} manifold. Then $\mathcal O$ is fine. 
\end{Prop}

In the \emph{proof} of this proposition, the following lemmata are crucial. 

\begin{Lem}[bumpfunc]
	Let $U,V\subset X_0$ be open, $K\subset V$ compact, \scth $\overline V\subset U$ and $U$ is a coordinate neighbourhood. There exists $f\in\mathcal O(U)$, $\supp f\subset V$, whose germ at $K$ is $1$, and which on $U$ takes values in the set $[0,1]\subset\reals$. 
\end{Lem}

\begin{proof}
	For any coordinate neighbourhood $U$, there exists an even algebra isomorphism $\mathcal O(U)\cong\mathcal C^\infty_{X_0}(U)\otimes\bigwedge(\cplxs^q)^*$ over the identity $U\to U$, and thus an embedding $\mathcal C^\infty_{X_0}(U)\subset\mathcal O(U)$ as an even real subalgebra. Hence the claim. 
\end{proof}

\begin{Lem}[cptext]
	Let $U\subset X_0$ be open and $f\in\mathcal O(U)$ with $\supp f$ compact. There exists $h\in\Gamma(\mathcal O)$ \scth $h|_U=f$ and $h_x=0$ \fa $x\not\in U$. 
\end{Lem}

\begin{proof}
	Let $V=X\setminus(\supp f)$. By definition of the support, we have \fa $x\in U\cap V$ that the germ $f_x=0$. Then $f$ on $U$ and $0$ on $V$ define the data of a global section of $\mathcal O$, by the sheaf axiom. 
\end{proof}

\begin{proof}[Proof of \thmref{Prop}{strucsheaffine}]
	Let $(U_\alpha)$ be a locally finite open cover of $X_0$. We need to construct $h_\alpha\in\Gamma(\mathcal O)$, \scth $\supp h_\alpha\subset U_\alpha$ form a locally finite family of closed sets, and the locally finite sum $\sum_\alpha h_\alpha=1$. Since $X_0$ is paracompact, passing to a locally finite refinement, we may assume that $U_\alpha$ are coordinate neighbourhoods. There are open subsets $V_\alpha,W_\alpha\subset X_0$ \scth $\overline{W_\alpha}\subset V_\alpha$ and $\overline{V_\alpha}\subset U_\alpha$ are compact, and $X_0=\bigcup_\alpha W_\alpha$. 
	
	By \thmref{Lem}{bumpfunc}, there exist $f_\alpha\in\mathcal O(U_\alpha)$ \scth $\supp f_\alpha\subset V_\alpha$, the germ of $f_\alpha$ at $\overline{W_\alpha}$ is $1$, and $f_\alpha$ takes values in $[0,1]$. Applying \thmref{Lem}{cptext} we extend $f_\alpha$ by zero to an element $f_\alpha\in\Gamma(\mathcal O)$ which takes values in $[0,1]$. 
		
	The sum $f=\sum_\alpha f_\alpha$ is locally finite and therefore exists in $\Gamma(\mathcal O)$. Since $(W_\alpha)$ is a cover of $X_0$, $f$ takes positive values. In particular, $f$ is invertible in $\Gamma(\mathcal O)$ (each germ being invertible). Define $h_\alpha=f^{-1}\cdot f_\alpha$.
\end{proof}

\begin{Rem}
	\thmref{Prop}{strucsheaffine}) can be used to prove an analogue of Batchelor's theorem for \emph{cs} manifolds: any \emph{cs} manifold $(X_0,\mathcal O)$ is (non-canonically) isomorphic to $(X_0,\bigwedge\mathcal E\otimes\cplxs)$ where $\mathcal E$ is a locally free $\mathcal C^\infty_{X_0}$-module. 
	
	We shall need fineness in the construction of a Berezin integral for \emph{cs} manifolds. Another useful consequence is the following corollary.
\end{Rem}

\begin{Cor}[globsectmor]
	Let $X$, $Y$ be \emph{cs} manifolds, and assume given an even superalgebra morphism $\vphi:\mathcal O_Y(Y_0)\to\mathcal O_X(X_0)$ which sends real-valued superfunctions to real-valued superfunctions. Then there exists a unique morphism $(f,f^*):X\to Y$ \scth $\vphi=f^*$ on $\mathcal O_Y(Y_0)$. 
\end{Cor}

\begin{Rem}[csalg]
	It is natural to introduce the following category $\mathrm{Alg}_{cs}$: The objects are pairs $(A,A_{0,\reals})$ where $A$ is a superalgebra over $\cplxs$, $A_{0,\reals}$ is a real subalgebra of $A_0$ containing $J_0=A_0\cap J$ where $J=A\cdot A_1$ is the graded ideal generated by $A_1$, \scth $A_{0,\reals}/J_0$ is a real form of $A_0/J_0=A/J$; the morphisms of such pairs $\phi:(A,A_{0,\reals})\to(B,B_{0,\reals})$ are even unital complex algebra morphisms $A\to B$ \scth $\phi(A_{0,\reals})\subset B_{0,\reals}$. 
	
	\thmref{Cor}{globsectmor} may then be rephrased as the statement that the global sections functor $\Gamma$ is fully faithful from \emph{cs} manifolds to $\mathrm{Alg}_{cs}$; to be precise, $\Gamma$ is defined on objects as $\Gamma(Z)=(\mathcal O_Z(Z_0),\mathcal O_Z(Z_0)_{0,\reals})$ where $\mathcal O_Z(Z_0)_{0,\reals}$ consists of the real-valued even elements of $\mathcal O_Z(Z_0)$. We call the objects of $\mathrm{Alg}_{cs}$ \emph{\emph{cs} algebras}. 
\end{Rem}

\begin{Def}
	Let $U$ be a \emph{cs} domain of graded dimension $p|q$, and suppose given a morphism $(f,f^*):U\to\reals^{p|q}_\cplxs$ which factors through an isomorphism onto an open subspace of $\reals^{p|q}_\cplxs$. 
	
	We define $x_j=f^*(\pr_j)\in\mathcal O_U(U_0)_0$, $j=1,\dotsc,p$, and $\xi_j=f^*(\theta_j)\in\mathcal O_U(U_0)_1$, $j=1,\dotsc,q$; here $\pr_j:\reals^p\to\reals$ are the coordinate projections, and $\theta_j=\pr_j\in(\cplxs^q)^*$ denote the standard generators of $\bigwedge(\cplxs^q)^*$. 
	
	The collection $(x,\xi)=(x_1,\dotsc,x_p,\xi_1,\dotsc,\xi_q)$ will be called a \emph{coordinate system} for $U$, and its entries will be called \emph{coordinates}. Note that the $x_j$ take real values; we denote the subset of $\reals^p$ consisting of all tuples $(x_1(u),\dotsc,x_p(u))$, where $u\in U_0$, by $x(U)=x(U_0)$. (Of course, the $\xi_j$ also takes real values, as does any odd superfunction.)	
\end{Def}

\begin{Rem}
	Hohnhold \cite{hohnhold-thesis} also considers complex coordinate systems. This is a useful concept when studying complex supermanifolds as \emph{cs} manifolds (by forgetting the holomorphic structure). The following mapping condition can, however, only be formulated with real coordinate systems. 
\end{Rem}

\begin{Prop}[coordmor]
	Let $X$ be a \emph{cs} manifold and $U$ a \emph{cs} domain of graded dimension $p|q$. Let $(x,\xi)$ be a coordinate system on $U$, where $y=(y_1,\dotsc,y_p)\in\mathcal O_X(X_0)^p_0$ and $\eta=(\eta_1,\dotsc,\eta_q)\in\mathcal O_X(X_0)^q_1$, \scth the map 
	\[
	X\to\cplxs^p:x\mapsto(y_1(x),\dotsc,y_p(x))
	\]
	takes values in $x(U)\subset\reals^p$. (In particular, the $y_j$ are real-valued.) Then there exists a unique morphism of \emph{cs} manifolds, $(f,f^*):X\to U$, \scth $f^*(x_j)=y_j$ and $f^*(\xi_j)=\eta_j$. 
\end{Prop}

\begin{proof}
	This follows as for real supermanifolds, \emph{cf.}~\cite[3.2]{schmitt-supergeom}. 
\end{proof}

\begin{Def}
	Let $V=V_0\oplus V_1$ be a finite-dimensional complex super-vector space \scth $V_0$ is equipped with a real form $V_{0,\reals}$. Then we call $(V,V_{0,\reals})$ a \emph{\emph{cs} vector space}. Moreover, $(V,V_{0,\reals})$ defines a \emph{cs} manifold $(V_{0,\reals},\mathcal C^\infty_{V_{0,\reals}}\otimes\bigwedge(V_1)^*)$ of graded dimension equal to the complex graded dimension of the super-vector space $V\!$. By abuse of notation, we shall denote it by $V\!$, and its structure sheaf by $\mathcal O_V$. We call $V$ a \emph{linear} \emph{cs} manifold. 
	
	The structure sheaf of $V$ comes with a natural $\ints$-grading, and its $0$th degree part is exactly $\mathcal C^\infty_{V_{0,\reals},\cplxs}$. Hence, in this special case, the given real form $V_{0,\reals}$ of $V_0$ is precisely what specifies the real-valued elements of $\mathcal O_V$. For this reason, we find it justifiable to use the subscript $_{0,\reals}$ both for the real form of $V_0$ and the sheaf of real-valued superfunctions (although only the former \emph{really} defines a real form). We forewarn the reader that we will systematically indulge in this abuse of notation. 
	
	The set $V^*$ of all complex linear forms on $V$ is embedded as a graded subspace in $\mathcal O_V(V_{0,\reals})$. The following statement follows from \thmref{Prop}{coordmor}. 
\end{Def}

\begin{Cor}[linearmor]
	Let $X$ be a \emph{cs} manifold and $V$ a linear \emph{cs} manifold. For any even linear map $\phi:V^*\to\mathcal O_X(X_0)$ \scth any $h\in\phi(V_{0,\reals}^*)$ takes real values, there exists a unique morphism $(f,f^*):X\to V$ \scth $f^*|_{V^*}=\phi$. 
	
	In other words, $V(X)=\Hom0{X,V}$, the set of $X$-points of $V\!$, is exactly 
	\[
		V(X)=(\mathcal O_X(X_0)\otimes V)_{0,\reals}=\mathcal O_X(X_0)_{0,\reals}\otimes_\reals V_{0,\reals}\oplus\mathcal O_X(X_0)_1\otimes_\cplxs V_1\ .
	\]
\end{Cor}

\begin{Par}[csmfprod]
	To finalise our disquisition on the fundamentals of \emph{cs} manifolds, we discuss the existence of finite products in this category. If $V$ and $W$ are \emph{linear} \emph{cs} manifolds, then we define the \emph{even product} $V\times W=(V_0\times W_0)\oplus(V_1\times W_1)$. The even part $V_0\times W_0$ of $V\times W$ is endowed with the real form $V_{0,\reals}\times W_{0,\reals}$. We thus obtain a linear \emph{cs} manifold which we denote again by $V\times W$. By \thmref{Cor}{linearmor}, there are morphisms $p_1:V\times W\to V$ and $p_2:V\times W\to W$ induced by the inclusions 
	\[
	V^*=(V_0^*\times0)\oplus(V_1^*\times0)\subset(V\times W)^*\supset(0\times W_0^*)\oplus(0\times W_1^*)=W^*
	\]
	and $(V\times W)^*\subset\mathcal O_{V\times W}(V_{0,\reals}\times W_{0,\reals})$. By the same token, $V\times W$ is the product of $V$ and $W$ in the category of \emph{cs} manifolds. 
	
	Next, consider open subspaces $A\subset V$ and $B\subset W$. We define $A\times B$ to be the (unique) open subspace of $V\times W$ whose underlying set is $A_0\times B_0\subset V_{0,\reals}\times W_{0,\reals}$. The morphisms $p_1$ and $p_2$ defined above restrict to morphisms $A\times B\to A$ and $A\times B\to B$, respectively. By \thmref{Prop}{coordmor}, $A\times B$ is the product of $A$ and $B$ in the category of \emph{cs} manifolds. 
	
	Given \emph{cs} manifolds $X$ and $Y$, we define $\mathcal O_{X\times Y}$ to be the (up to canonical isomorphism) unique sheaf on $X_0\times Y_0$ \scth $\mathcal O_{X\times Y}|_{U_0\times V_0}=\mathcal O_{U\times V}$ for all coordinate neighbourhoods $U\subset X$ and $V\subset Y$. This determines a \emph{cs} manifold $X\times Y$. The canonical morphisms $U\times V\to U$ and $U\times V\to V$ for all coordinate neighbourhoods $U$, $V$ determine morphisms $X\times Y\to X$ and $X\times Y\to Y$. As in the case of real supermanifolds, one shows that $X\times Y$ is the product of $X$ and $Y$ in the category of \emph{cs} manifolds. 
	
	At this point we mention that for any \emph{cs} manifold $X$, $\mathcal O_X$ carries a natural structure of a sheaf of nuclear Fr\'echet spaces. Indeed, for any coordinate neighbourhood $U$, $\mathcal O_X(U_0)$ is a nuclear Fr\'echet space for the standard topology induced from $\mathcal C^\infty_{X_0}(U_0)\otimes\bigwedge(\cplxs^q)^*$, and one may take locally convex projective limits with respect to the restriction maps (here the paracompactness of $X_0$ ensures that countable limits suffice). Using the nuclearity, one can show that for any \emph{cs} manifolds $X$ and $Y$, and any open subsets $U_0\subset X_0$ and $V_0\subset Y_0$, the inclusion $\mathcal O_X(U)\otimes\mathcal O_Y(V)\to\mathcal O_{X\times Y}(U\times V)$ given by $f\otimes h\mapsto p_1^*f\cdot p_2^*h$ induces an isomorphism on the completion of the tensor product w.r.t.~any locally convex tensor product topology (for instance, one may take the projective tensor product topology). 
\end{Par}

\subsection{An inverse function theorem for \emph{cs} manifolds}

We give an easy special case of the inverse function theorem that is useful in various situations.

\begin{Def}
	Let $X=(X_0,\mathcal O)$ be a \emph{cs} manifold. For $x\in X_0$, we define the \emph{tangent space} $T_xX=\Der0{\mathcal O_x,\cplxs}$. If $X=X_0$ is an ordinary manifold, then $T_xX=T_xX_0\otimes\cplxs$. In general, we have $\sdim_\cplxs T_xX=\sdim_\cplxs\ger m_x/\ger m_x^2$ by standard arguments, where $\ger m_x$ is the maximal ideal of the local algebra $\mathcal O_x$ and $\sdim$ denotes graded dimension.
	
	Given a morphism $f:X\to Y$ of \emph{cs} manifolds, we define \fa $x\in X_0$ even $\cplxs$-linear maps $T_xf:T_x X\to T_{f(x)}Y$, called \emph{tangent maps}, by 
	\begin{equation}\label{eq:tangentmapdef}
		(T_xf)(\xi)h=\xi(f^*h)\mathfa\xi\in T_xX\,,\,h\in\mathcal O_{Y,f(x)}\ .
	\end{equation}
	Clearly, one has, for any morphism $g:Y\to Z$, the chain rule
	\[
		T_x(g\circ f)=T_{f(x)}g\circ T_xf\ .
	\]
\end{Def}

Our main application of tangent maps is the following proposition.

\begin{Prop}[invfnthm]
	Let $f:X\to Y$ be a morphism of \emph{cs} manifolds. Then $f$ is an isomorphism if and only if $f:X_0\to Y_0$ is bijective, and \fa $x\in X_0$, $T_xf:T_x X\to T_{f(x)} Y$ is an isomorphism of vector spaces. 
\end{Prop}

\begin{proof}
	For real supermanifolds, this is \cite[Corollary to Theorem 2.16]{kostant-supergeom}. The case of \emph{cs} manifolds can be treated analogously. 
\end{proof}

\subsection{The functor of points}\label{app:s-points}

We recall the functor of points. We use it chiefly as a device to prove identities of morphisms, so we do not go very deep in our discussion. 

\begin{Par}
	Let $\SMan_{cs}$ denote the category of supermanifolds. For $U,X\in\SMan_{cs}$, let $X(U)$ denote the set of morphisms $U\to X$ in $\SMan_{cs}$. Elements $x\in X(U)$ are called $U$-points, and one writes $x\in_UX$. If $*$ is the $0|0$-dimensional \emph{cs} manifold whose underlying topological space is a point, then $X(*)=X_0$ is the real manifold underlying $X$. (Hence the terminology.)
\end{Par}

\begin{Par}
	The assignment $X\mapsto X(-)$ extends in a natural way to a functor from $\SMan_{cs}$ to the category $[\SMan_{cs}^{\mathrm{op}},\mathrm{Set}]$ of set-valued functors on \emph{cs} manifolds. By Yoneda's lemma \cite{maclane}, this functor is fully faithful, and we call it the \emph{Yoneda embedding}. In particular, to define morphisms of \emph{cs} manifolds, it suffices to define morphisms of their images under the Yoneda embedding. 
	
	The Yoneda embedding commutes with finite products. It follows that it induces a fully faithful embedding of the category of \emph{cs} Lie supergroups (see below) into the category $[\SMan_{cs}^{\mathrm{op}},\mathrm{Grp}]$ of group-valued functors on \emph{cs} manifolds.
\end{Par}

\begin{Rem}
	A more complete discussion of the functor of points, and its applications to the theory Lie supergroups, is given in \cite[\S\S~2.8--11]{deligne-morgan}.
\end{Rem}

\renewcommand\thesection{B}
\section{Appendix: The category of \emph{cs} Lie supergroups}\label{s:app-b}

\subsection{The category of \emph{cs} Lie supergroups}

In this section, we will use standard notions of Lie theory and of the theory of Lie superalgebras without explicit reference. Standard texts for the latter would be \cite{kac-liesuperalgs,scheunert-liesuperalgs}. Useful references for the former could be \cite{knapp-beyond,hnbook-new}, but most Lie theory texts should be sufficient.

\begin{Def}
	A \emph{\emph{cs} Lie supergroup} is a group object in the ca\-te\-gory of \emph{cs} manifolds. Thus, it is a tuple $(G,m,i,e)$ where $G$ is a \emph{cs} manifold and $m:G\times G\to G$, $i:G\to G$, and $e:*\to G$ are morphisms of \emph{cs} manifolds subject to the obvious axioms. (Here, $*$ is the terminal object in the category of \emph{cs} manifolds.) In particular, $G_0$ is a real Lie group.
	
	A \emph{morphism $f:G\to H$ of \emph{cs} Lie supergroups} is a morphism of group objects in the category of \emph{cs} manifolds. \emph{I.e.}, it is a morphism of \emph{cs} manifolds, and satisfies the equations 
	\[
	\phi\circ m_G=m_H\circ(\phi\times\phi)\ ,\ \phi\circ i_G=i_H\circ\phi\ ,\ \phi\circ e_G=e_H\ .
	\]
\end{Def}

There is an easy way to construct \emph{cs} Lie supergroups, due to J.-L.~Koszul \cite{koszul-superaction} in the real case, and it uses the following concept. 

\begin{Def}
	A \emph{\emph{cs} supergroup pair} (a.k.a.~\emph{super Harish-Chandra pair}) is actually a triple $(G_0,\ger g,\Ad)$ subject to the following assumptions: $G_0$ is a \emph{real} Lie group with Lie algebra $\ger g_{0,\reals}$; $\ger g=\ger g_0\oplus\ger g_1$ is a \emph{complex} Lie superalgebra \scth $\ger g_{0,\reals}\subset\ger g_0$ is a real form of $\ger g_0$; and $\Ad:G_0\times\ger g\to\ger g$ is a smooth linear action of $G_0$ by even Lie superalgebra automorphisms which extends the adjoint action of $G$ on $\ger g_{0,\reals}$, and whose differential $d\Ad:\ger g_{0,\reals}\times\ger g\to\ger g$ is the restriction of the bracket $[\cdot,\cdot]:\ger g\times\ger g\to\ger g$. 
	
	By abuse of notation, we will generally write $(G_0,\ger g)$ for a \emph{cs} supergroup pair (hence the parlance). The action $\Ad=\Ad_G$ is understood, although strictly speaking, it is not determined by $G_0$ and $\ger g$ unless $G_0$ is connected. 
	
	A \emph{morphism of \emph{cs} supergroup pairs} is a pair $(f,df):(G_0,\ger g)\to(H_0,\ger h)$ fulfilling the following conditions: $f:G_0\to H_0$ is a morphism of real Lie groups; $df:\ger g\to\ger h$ is an even morphism of complex Lie superalgebras extending the differential of $f$; and $df$ is $G_0$-equivariant for the $G_0$-action on $\ger h$ induced by $f$, \emph{i.e.}
	\[
		\Ad_H(f(g))(df(x))=df(\Ad_G(g)(x))\mathfa g\in G_0\,,\,x\in\ger g\ .
	\]
\end{Def}

\begin{Par}
	A simple but salient point about \emph{cs} supergroup pairs $(G_0,\ger g)$ is that $\mathcal C^\infty_{G_0}\otimes\cplxs$ is a sheaf of $\ger g_0$-modules. Indeed, if $U\subset G_0$ is open, $x\in\ger g_{0,\reals}$, and $f\in\Ct[^\infty]0{U,\cplxs}$, then we may define
	\[
		(r_xf)(g)=\frac d{dt}f\Parens1{g\exp(tx)}\Big|_{t=0}\mathfa g\in U\ .
	\]
	By complex linear extension, this defines a $\ger g_0$-module structure on $\Ct[^\infty]0{U,\cplxs}$ which is compatible with the restriction maps of $\mathcal C^\infty_{G_0}\otimes\cplxs$. 
\end{Par}

\begin{Par}	
	Let $\phi:\ger h\to\ger g$ be a morphism of complex Lie superalgebras and $\mathcal V$ a graded $\ger h$-module sheaf. Then $\Uenv0{\ger g}$ is an $\ger h$-module \emph{via} $x.v=\phi(x)v$ \fa $x\in\ger h$, $v\in\Uenv0{\ger g}$. One defines the \emph{coinduced module sheaf} $\Coind[_{\ger h}^{\ger g}]0{\mathcal V}$ by 
	\[
		\Coind[_{\ger h}^{\ger g}]0{\mathcal V}(U)=\underline{\mathrm{Hom}}_{\Uenv0{\ger h}}\Parens0{\Uenv0{\ger g},\mathcal V(U)}\ .
	\]
	Here, $\underline{\mathrm{Hom}}$ denotes \emph{inner} $\mathrm{Hom}$, \emph{i.e.}~linear maps are considered without a parity constraint. Then $\Coind[_{\ger h}^{\ger g}]0{\mathcal V}$ is a graded $\ger g$-module sheaf with $\Uenv0{\ger g}$-module structure 
	\begin{equation}\label{eq:coind-infinaction}
		(r_uf)(v)=(-1)^{\Abs0u(\Abs0F+\Abs0v)}f(vu)
	\end{equation}
	\fa $u,v\in\Uenv0{\ger g}$, $f\in\Coind[_{\ger h}^{\ger g}]0{\mathcal V(U)}$. If $\mathcal V$ is purely even, then we may replace the sign $(-1)^{\Abs0u(\Abs0F+\Abs0v)}$ by $(-1)^{\Abs0u}$. One is given a canonical $\ger h$-equivariant morphism $\Coind[_{\ger h}^{\ger g}]0{\mathcal V}\to\mathcal V$ by the assignment $f\mapsto F(1)$. 
	
	If $\mathcal V$ is a sheaf of superalgebras and $\ger h$ acts by graded derivations, then we obtain on $\Coind[_{\ger h}^{\ger g}]0{\mathcal V}$ the structure of a sheaf of $\ger g$-superalgebras. The multiplication is 
	 \begin{equation}\label{eq:coind-algstr}
	 	f\otimes f'\mapsto\mu\circ(f\otimes f')\circ\Delta
	 \end{equation}
	 where $\mu$ is the multiplication of $\mathcal V$, and $\Delta:\Uenv0{\ger g}\to\Uenv0{\ger g}\otimes\Uenv0{\ger g}=\Uenv0{\ger g\oplus\ger g}$ is the coproduct of $\Uenv0{\ger g}$ (\emph{i.e.}~the unique extension to an even unital algebra morphism of the map $\ger g\to\Uenv0{\ger g}\otimes\Uenv0{\ger g}:x\mapsto x\otimes1+1\otimes x$). Let $\eps:\Uenv0{\ger g}\to\cplxs$ be the unique extension to an even unital algebra morphism of the zero map $\ger g\to\cplxs$; then the unit of $\Coind[_{\ger h}^{\ger g}]0{\mathcal V}$ is $\eta\circ\eps:\Uenv0{\ger g}\to\cplxs\to\mathcal V$ where $\eta$ is the unit of $\mathcal V$. 
\end{Par}

\begin{Par}	 
	 Next, assume again that $(G_0,\ger g)$ is a \emph{cs} supergroup pair. Define a superalgebra sheaf by $\mathcal O_G=\Coind[_{\ger g_0}^{\ger g}]0{\mathcal C^\infty_{G_0}\otimes\cplxs}$.  If $\beta:S(\ger g)\to\Uenv0{\ger g}$ is the supersymmetrisation map, then the morphism $\mathcal O_G\to\mathcal C^\infty_G\otimes\bigwedge\ger g_1^*$ defined by $f\mapsto (f\circ\beta)|_{\bigwedge\ger g_1}$ is an isomorphism of superalgebra sheaves, since $\beta$ is a supercoalgebra isomorphism. In particular, $G=(G_0,\mathcal O_G)$ is a \emph{cs} manifold. 
	 
	 We note that the action of $\ger g_{0,\reals}$ on $\mathcal O_G$ integrates to an action of $G_0$ via 
	 \begin{equation}\label{eq:coind-grpaction}
	 	(r_hf)(u;g)=f(\Ad(h^{-1})(u);gh)
	 \end{equation}
	 \fa $f\in\mathcal O_G(U)$, $u\in\Uenv0{\ger g}$, $g\in G_0$, $h\in g^{-1}U$. Here and in the following, we use the notation $f(u;g)=f(u)(g)$ for $f\in\mathcal O_G(U)$, $u\in\Uenv0{\ger g}$, and $g\in U$. 
\end{Par}

\begin{Prop}[coindfunctor]
	Let $(G_0,\ger g)$ be a \emph{cs} supergroup pair. Let $m,i,e$ be the structure maps of the real Lie group $G_0$. We define $m^*:\mathcal O_G\to m_*\mathcal O_{G\times G}$, $i^*:\mathcal O_G\to i_*\mathcal O_G$, and $e^*:\mathcal O_G\to\cplxs$ by 
	\begin{equation}\label{eq:multdef}
		m^*f(u\otimes v;g,h)=f\Parens1{\Ad(h^{-1})(u)v;gh}
	\end{equation}
	\fa $f\in\mathcal O_G(U)$, $u,v\in\Uenv0{\ger g}$, $(g,h)\in m^{-1}(U)$; 
	\begin{equation}\label{eq:invdef}
		i^*f(u;g)=f\Parens1{\Ad(g)(S(u));g^{-1}}
	\end{equation}
	\fa $f\in\mathcal O_G(U)$, $u\in\Uenv0{\ger g}$, $g\in i^{-1}(U)$; and $e^*f=f(1;1)$. Here, $S:\Uenv0{\ger g}\to\Uenv0{\ger g}$ is determined by $S(1)=1$, $S(x)=-x$ \fa $x\in\ger g$, and $S(uv)=(-1)^{\Abs0u\Abs0v}S(v)S(u)$ \fa homogeneous $u,v\in\Uenv0{\ger g}$. Then $C(G_0,\ger g)=(G,(m,m^*),(i,i^*),(e,e^*))$ is a \emph{cs} Lie supergroup. 
	
	If $(f,df):(G_0,\ger g)\to(H_0,\ger h)$ is a morphism of \emph{cs} supergroup pairs, then the assignment $C(f,df)=(f,f^*)$, where $f^*:\mathcal O_H\to f_*\mathcal O_G$ is defined by 
	\[
		(f^*h)(u;g)=h(df(u);f(g))\mathfa h\in\mathcal O_H(U)\,,\,u\in\Uenv0{\ger g}\,,\,g\in f^{-1}(U)\ ,
	\]
	gives a morphism of \emph{cs} Lie supergroups. Here, $df$ is extended uniquely to an even unital algebra morphism $\Uenv0{\ger g}\to\Uenv0{\ger h}$. This defines a functor $C$ from \emph{cs} supergroup pairs to \emph{cs} Lie supergroups. 
\end{Prop}

\begin{proof}
	The details of the proof are somewhat tedious, and most of the calculations are straightforward, so we only indicate the salient points. 
	
	First, by \thmref{Cor}{globsectmor}, it is sufficient to work on the level of global sections; this tidies matters up a little. One needs to check that $m^*$ is well-defined, \emph{i.e.}~that $m^*f$ is $(\ger g_0\oplus\ger g_0)$-equivariant for any $f$; this follows from the fact that $G_0$ acts on $\ger g$ by Lie superalgebra automorphisms. From $(\Delta\otimes\id)\circ\Delta=(\id\otimes\Delta)\circ\Delta$ (the coassociativity relation), it follows that 
	\[
		(m,m^*)\circ((m,m^*)\times\id)=(m,m^*)\circ(\id\times(m,m^*))\ .
	\]
	
	Next, the multiplication of $\mathcal O_G(G)$ extends through the canonical embedding $\Gamma(\mathcal O_G)\otimes\Gamma(\mathcal O_G)\to\Gamma(\mathcal O_{G\times G})$ to a unital algebra morphism which coincides with $\delta^*:\Gamma(\mathcal O_{G\times G})\to\Gamma(\mathcal O_G)$, given by $\delta^*f(u;g)=f(\Delta(u);\delta(g))$ where $\delta(g)=(g,g)$. But then $(\delta,\delta^*):G\to G\times G$ is the diagonal morphism in \emph{cs} manifolds given by the universal property of the product. The equations
	\[
		(m,m^*)\circ((i,i^*)\times\id)\circ(\delta,\delta^*)=(e,e^*)\circ(*,1)=(m,m^*)\circ(\id\times(i,i^*))\circ(\delta,\delta^*)
	\]
	where $(*,1):G\to *$ is the canonical morphism to the terminal object, then follow from the fact that $S$ is an antipode for the bialgebra $\Uenv0{\ger g}$. The latter fact means that $\mu\circ(S\otimes\id)\circ\Delta=1=\mu\circ(\id\otimes S)\circ\Delta$ where $\mu$ is the algebra multiplication of $\Uenv0{\ger g}$, and can be extracted from standard texts on Hopf algebras such as \cite{sweedler}, or the excellent introduction \cite{cartier-hopf}.
	
	The remaining statements are easily substantiated. 
\end{proof}

\begin{Par}
	Let $G$ be a \emph{cs} Lie supergroup with structure morphisms $(m,m^*)$, $(i,i^*)$, and $(e,e^*)$. We denote $1=e(*)\in G_0$ the neutral element of the underlying Lie group. It is clear that $e^*f=f(1)$ \fa $f\in\Gamma(\mathcal O_G)$. Let $\ger g\subset\underline{\mathrm{Hom}}\Parens0{\Gamma(\mathcal O_G),\cplxs}$ consist of the graded derivations along $e^*$. \emph{I.e.}, $\ger g$ is spanned by the homogeneous linear maps $x:\Gamma(\mathcal O_G)\to\cplxs$ which satisfy
	\[
		x(f\cdot h)=x(f)h(1)+(-1)^{\Abs0x\Abs0f}f(1)x(h)
	\]
	\fa homogeneous $f,h\in\Gamma(\mathcal O_G)$. 
	
	For any $x\in\ger g$, one defines by continuous linear extension from the algebraic tensor product a linear map $\id\Hat\otimes\,x:\Gamma(\mathcal O_{G\times G})\to\Gamma(\mathcal O_G)$. It is a graded derivation along $\id\Hat\otimes\,e^*$. Then $\mathcal L_x=(\id\Hat\otimes\,x)\circ m^*\in\Der0{\Gamma(\mathcal O_G)}$ is a graded derivation (along the map $\id$). 
\end{Par}

\begin{Lem}[invarvf]
	For any \emph{cs} Lie supergroup $G$, the map $\mathcal L:\ger g\to\Der0{\Gamma(\mathcal O_G)}$ defines a bijection onto the set of graded derivations $d$ which satisfy 
	\[
	(p_1,m)^*\circ(\id\Hat\otimes\,d)=(\id\Hat\otimes\,d)\circ(p_1,m)^*
	\]
	where $(p_1,m):G\times G\to G\times G$. 
\end{Lem}

\begin{proof}
	The map $\mathcal L$ is well-defined because of the associativity equation for $(m,m^*)$ which translates in particular to $(\id\Hat\otimes\,m^*)\circ m^*=(m^*\Hat\otimes\,\id)\circ m^*$. Then one recovers $x$ from $x=e^*\circ \mathcal L_x$, so $\mathcal L$ is injective. One may also use this equation to determine the image of $\mathcal L$.
\end{proof}

\begin{Par}[adjdef]
	Let $G$ be a \emph{cs} Lie supergroup. We shall call the graded derivations in the image of $\mathcal L$ \emph{left-invariant vector fields}. The set of left-invariant vector fields is easily shown to be a Lie subsuperalgebra of $\Der0{\Gamma(\mathcal O_G)}$ under the super-commutator bracket. In particular, $\ger g$ has the structure of a finite-dimensional complex Lie superalgebra; we call $\ger g$ the \emph{Lie superalgebra of $G$}. 
	
	Let $g\in G_0$. Then $g$ defines a morphism $\Gamma(\mathcal O_G)\to\cplxs$, and there is an even unital algebra morphism $g\otimes\id\otimes\,g^{-1}:\Gamma(\mathcal O_{G\times G\times G})\to\Gamma(\mathcal O_G)$. We define morphisms $m^{(2)}=m\circ(m\times\id)$ and $c_g^*=(g\otimes\id\otimes\,g^{-1})\circ m^{(2)*}$, the \emph{conjugation} by $g$. Then it is natural to define $\Ad:G_0\times\ger g\to\ger g$ by $\Ad(g)(x)=x\circ c_g^*$. One checks that this is a finite-dimensional continuous linear representation, and hence smooth. Straightforward calculations prove the following lemma. 
\end{Par}

\begin{Lem}
	Let $G$ be a \emph{cs} Lie supergroup with Lie superalgebra $\ger g$ and underlying real Lie group $G_0$. With the adjoint action $\Ad$, $(G_0,\ger g)$ is a \emph{cs} supergroup pair. 
\end{Lem}

\begin{Par}
	Let $G$ be a \emph{cs} Lie supergroup given by $G=C(G_0,\ger g)$ for a \emph{cs} supergroup pair $(G_0,\ger g)$. Then there is a canonical isomorphism $T_1G\cong\ger g$ which can be derived from \thmref{Lem}{invarvf}. 
	
	For any $g\in G_0$, we define morphisms of \emph{cs} manifolds $L_g,R_g:G\to G$ by taking the left (resp.~right) translation by $g$ on the level of $G$, and setting 
	\begin{gather*}
		L_g^*f(u;h)=m^*f(1\otimes u;g,h)=f(u;gh)\ ,\\
		R_g^*f(u;h)=m^*f(u\otimes1;h,g)=f(\Ad(g^{-1})(u);hg)\ ,
	\end{gather*}
	\fa $f\in\mathcal O_G(U)$, $u\in\Uenv0{\ger g}$, and $h\in g^{-1}U$ (resp.~$h\in Ug^{-1}$). 

	Clearly, $L_g$, $R_g$ are isomorphisms with inverses $L_{g^{-1}}$ and $R_{g^{-1}}$, respectively. In particular, if we write $dL_g=T_1L_g$ and $dR_g=T_1R_g$, then $dL_g,dR_g:\ger g\to T_gG$ are isomorphisms, by \thmref{Prop}{invfnthm}.
	
	A direct calculation shows that
	\begin{equation}\label{eq:tangentmult}
		\Parens1{T_{(g,h)}m}(dL_g(x),dL_h(y))=dL_{gh}(\Ad(h^{-1})(x)+y)
	\end{equation}
	\fa $g,h\in G_0$, $x,y\in\ger g$. 
\end{Par}

\begin{Prop}[catequiv]
	The functor $C$ is an equivalence of categories from \emph{cs} supergroup pairs to \emph{cs} Lie supergroups. 
\end{Prop}

\begin{Rem}
	In the real case, \thmref{Prop}{catequiv} is due to B.~Kostant \cite{kostant-supergeom}. As remarked above, we use a construction due to J.-L.~Koszul \cite{koszul-superaction}.
\end{Rem}

\begin{proof}[\protect{Proof of \thmref{Prop}{catequiv}}]
	We need to prove that $C$ is fully faithful and essentially surjective. We begin with the essential surjectivity. 
	
	To that end, let $G$ be a \emph{cs} Lie supergroup with Lie superalgebra $\ger g$ and underlying real Lie group $G_0$. We define a morphism $\phi:C(G_0,\ger g)\to G$, $\phi=(\id,\phi^*)$, by specifying $\phi^*:\mathcal O_G\to\mathcal O'=\Coind[_{\ger g_0}^{\ger g}]0{\mathcal C^\infty_{G_0}\otimes\cplxs}$ \emph{via} 
	\[
		(\phi^*f)(u;g)=(-1)^{\Abs0f\Abs0u}(\mathcal L_uf)(g)\mathfa f\in\mathcal O_G(U)\,,\,u\in\Uenv0{\ger g}\,,\,g\in U\ .
	\]
	
	Here, $\mathcal L$ is the extension of $\mathcal L:\ger g\to\Der0{\mathcal O_G(U)}$ to an even unital algebra morphism $\Uenv0{\ger g}\to\underline{\mathrm{End}}\Parens0{\mathcal O_G(U)}$ (inner $\mathrm{End}$). To see that $(\id,\phi^*)$ is an isomorphism of \emph{cs} manifolds, we apply \thmref{Prop}{invfnthm}. 
	
	Indeed, for $g\in G_0$, denote the maximal ideal of $\mathcal O_{G,g}$ and $\mathcal O'$ by $\ger m_g$ and $\ger n_g$, respectively. Then 
	\[
		\ger n_g=\Set1{h\in\mathcal O'_g}{h(1;g)=0}\ ,\ \ger n_g^2=\Set1{h\in\mathcal O'_g}{h(x;g)=0\text{ \fa}x\in\ger g}\ ,
	\]
	as follows from the definition of the algebra structure. Then $\phi^*\ger m_g\subset\ger n_g$. On the other hand, if $h\in\mathcal O_g$ is not contained in $\ger m_g^2$, then $(\mathcal L_xh)(g)\neq0$ \fs $x\in\ger g$. Thus, $\phi^*$ induces an injection $\ger m_g/\ger m_g^2\to\ger n_g/\ger n_g^2$, and the tangent map $T_g(\id,\phi^*):T_g(G_0,\mathcal O')=(\ger n_g/\ger n_g^2)^*\to(\ger m_g/\ger m_g^2)^*=T_gG$ is surjective. Since $G$ and $(G_0,\mathcal O')$ have the same graded dimensions, \thmref{Prop}{invfnthm} implies that $(\id,\phi^*)$ is an isomorphism of \emph{cs} manifolds. 

	One can use the isomorphism $\phi^*$ and \thmref{Lem}{invarvf} to prove that $C$ is fully faithful. We leave the details to the reader. 
\end{proof}

\begin{Rem}
	The argument using the inverse function theorem in the proof of \thmref{Prop}{catequiv} is due to E. G.$\,$Vishnyakova \cite{vishnyakova-supergrp}. 

	By a result of H.~Hohnhold \cite[\S~4.4, Proposition 12]{hohnhold-thesis}, there is a forgetful functor from complex supermanifolds to \emph{cs} manifolds, and it is faithful. Moreover, there is a natural notion of complex supergroup pairs (modeled over complex Lie groups and complex Lie superalgebras), and the forgetful functor from complex to \emph{cs} supergroup pairs which maps $(G_0,\ger g)$ to $(G_0,\ger g\otimes_\reals\cplxs)$, forgetting complex structure on $G_0$, is also faithful. 
	
	If we apply the functor $C$ to a \emph{cs} supergroup pair which comes from a complex supergroup pair, then we obtain the \emph{cs} Lie supergroup of a complex Lie supergroup. In particular, we have reproved the following result which is also due to E.~G.$\,$Vishnyakova \cite{vishnyakova-supergrp}.
\end{Rem}

\begin{Cor}[csgrpcateq]
	The categories of complex Lie supergroups and of complex supergroup pairs are equivalent.
\end{Cor}

\begin{Par}[glcs]
	The general linear supergroup deserves a separate discussion in the framework of \emph{cs} manifolds. If $V$ is a \emph{real} super-vector space, then the real Lie supergroup $\mathcal{GL}(V)$ can be complexified to a \emph{cs} Lie supergroup $\mathcal{GL}(V)_\cplxs$; in particular, this is the case for $\mathcal{GL}(p|q)=\mathcal{GL}(\reals^{p|q})$. For any real super-vector space $V$, the \emph{cs} supergroups $\mathcal{GL}(V)_\cplxs$ and $C(\GL(V_0)\times \GL(V_1),\End0{V\otimes\cplxs})$ are canonically isomorphic. Let us describe the \emph{cs} supergroup $\mathcal{GL}(V)_\cplxs$ in the language of points. 

	To that end, assume more generally that $(A,A_{0,\reals})$ be a supercommutative \emph{cs} algebra (see \thmref{Rem}{csalg}), and $(V,V_{0,\reals})$ be a \emph{cs} vector space. The free $A$-module $A\otimes V=A\otimes_\cplxs V$, and $\underline{\mathrm{End}}_A\Parens0{A\otimes V}=A\otimes\underline{\mathrm{End}}(V)$. Let $\End[_{cs}]0{A\otimes V}_\reals$ denote the following set of block matrices:
	\[
		\begin{Matrix}1
			A_{0,\reals}\otimes_\reals\End[_\reals]0{V_{0,\reals}}&A_1\otimes_\cplxs\Hom[_\cplxs]0{V_1,V_0}\\
			A_1\otimes_\cplxs\Hom[_\cplxs]0{V_0,V_1}&A_0\otimes_\cplxs\End[_\cplxs]0{V_1}
		\end{Matrix}\ .
	\]
	
	Recall the notation $J$ for the graded ideal of $A$ generated by $A_1$. Then $J_0$ is a complex subspace of $A_0$, so that 
	\begin{equation}\label{eq:j-end}
	J_0\otimes_\reals\End[_\reals]0{V_{0,\reals}}=(J_0\otimes_\cplxs\cplxs)\otimes_\reals\End0{V_{0,\reals}}=J_0\otimes_\cplxs\End[_\cplxs]0{V_0}\ ,
	\end{equation}
	because $\cplxs\otimes_\reals\End[_\reals]0{V_{0,\reals}}=\End[_\cplxs]0{V_0}$. Since $J_0\subset A_{0,\reals}$, Equation \eqref{eq:j-end} shows that the composition of block matrices leaves the set $\End[_{cs}]0{A\otimes V}_\reals$ invariant, and turns it into an associative $\reals$-algebra. 
	
	Let $\GL_{cs}(A\otimes V)$ denote its group of units. This is a subgroup of $\GL(A\otimes V)$, the group of units of $\End[_A]0{A\otimes V}=\underline{\mathrm{End}}_A(A\otimes V)_0$. For $(A,A_{0,\reals})=\Gamma(\reals^{p|q}_\cplxs)$, we obtain $\GL_{cs}(p|q,A)\subset\GL(p|q,A)$. 
	
	If $Z$ is any \emph{cs} manifold and $V$ is now again a \emph{real} super-vector space, then the set $\mathcal{GL}(V)_\cplxs(Z)$ of $Z$-points of $\mathcal{GL}(V)_\cplxs$ is exactly $\GL_{cs}(\mathcal O_Z(Z_0)\otimes V_\cplxs)$.  
\end{Par}
		
	Any linear \emph{cs} manifold is isomorphic to the complexification of a real linear supermanifold. However, this isomorphism is non-canonical. To describe linear actions of \emph{cs} Lie supergroups, one therefore needs a generalisation of the general linear group which does not take a real structure on $V_1$ into account. The description of the generalised points of the complexification of the \emph{real} general linear supergroup indicates the correct definition. 

\begin{Par}
	Let $(A,A_{0,\reals})$ be a supercommutative \emph{cs} algebra, and $(V,V_{0,\reals})$ a \emph{cs} vector space. In general, there does not exist a real subsuperalgebra $B$ of $\underline{\mathrm{End}}\Parens0{A\otimes V}$, spanning over $\cplxs$, \scth $\End[_{cs}]0{A\otimes V}_\reals=B_0$. However, for $(A,A_{0,\reals})=(\cplxs,\reals)$, we can consider the set $\underline{\mathrm{End}}_{cs}(V)$ of all block matrices
	\[
		\begin{Matrix}1
			\End[_\reals]0{V_{0,\reals}}\otimes_\reals\cplxs&\Hom[_\cplxs]0{V_1,V_0}\\
			\Hom[_\cplxs]0{V_0,V_1}&\End[_\cplxs]0{V_1}\otimes_\reals\cplxs
		\end{Matrix}\ .
	\] 
	
	Let $\mu_{00}$ denote the inverse of the isomorphism $\End[_\reals]0{V_{0,\reals}}\otimes_\reals\cplxs\to\End[_\cplxs]0{V_0}$ of $\cplxs$-algebras. Furthermore, define $\mu_{ij}:\Hom[_\cplxs]0{V_i,V_j}\otimes_\reals\cplxs\to\Hom[_\cplxs]0{V_i,V_j}$, for $(i,j)\in{(0,1),(1,0)}$, by 
	\[
		\mu_{ij}(a\otimes z)=z\cdot a\mathfa a\in\Hom[_\cplxs]0{V_i,V_j}\,,\,z\in\cplxs\ .
	\]
	We obtain an algebra multiplication on $\underline{\mathrm{End}}_{cs}(V)$ by 
	\[
		\begin{Matrix}1
			a&b\\c&d
		\end{Matrix}\cdot
		\begin{Matrix}1
			e&f\\ g&h
		\end{Matrix}
		=
		\begin{Matrix}1
			ae+\mu_{00}(bg)&\mu_{10}(a(f\otimes1)+(b\otimes1)h)\\
			\mu_{01}((c\otimes1)e+d(g\otimes1))&cf\otimes1+dh
		\end{Matrix}
	\]
	for any $a,e\in\End[_\reals]0{V_{0,\reals}}\otimes_\reals\cplxs$, $b,f\in\Hom0{V_1,V_0}$, $c,g\in\Hom0{V_0,V_1}$, and any $d,h\in\End0{V_1}\otimes_\reals\cplxs$. One checks that (with the obvious grading), this turns $\underline{\mathrm{End}}_{cs}(V)$ into a complex superalgebra \scth a real form of the even part is given by $\underline{\mathrm{End}}_{cs}(V)_{0,\reals}=\End[_{cs}]0V_\reals=\End[_\reals]0{V_{0,\reals}}\oplus\End[_\cplxs]0{V_1}$. 
		
	We define 
	\[
		\mathcal{GL}(V,V_{0,\reals})=C\Parens1{\GL(V_{0,\reals})\times\GL(V_1),\underline{\mathrm{End}}_{cs}(V)}
	\] 
	where the adjoint action is defined by the conjugation of block matrices. As a \emph{cs} manifold, $\mathcal{GL}(V,V_{0,\reals})$ is the open subspace of the linear \emph{cs} manifold associated with $(\underline{\mathrm{End}}_{cs}(V),\End0{V_{0,\reals}}\oplus\End0{V_1})$ whose underlying subset of $\End0{V_{0,\reals}}\oplus\End0{V_1}$ is $\GL(V_{0,\reals})\times\GL(V_1)$. In particular, for $A=\mathcal O_Z(Z_0)$ ($Z$ being any \emph{cs} manifold), 
	\[
	\mathcal{GL}(V,V_{0,\reals})(Z)=\GL_{cs}(A\otimes V)\ ,
	\]
	as follows from \thmref{Cor}{linearmor}. 
\end{Par}

\begin{Def}
	Let $G$ be a \emph{cs} Lie supergroup and $X$ a \emph{cs} manifold. A morphism $\alpha:G\times X\to\mathcal X$ is called an \emph{action} if 
	\[
		\alpha\circ(m\times\id)=\alpha\circ(\id\times\alpha)\nd\alpha\circ(e\times\id)=\id\ .
	\]
	The action is called \emph{linear} if $X$ is the linear \emph{cs} manifold associated to a \emph{cs} vector space $(V,V_{0,\reals})$, and $\alpha^*(V^*)\subset\mathcal O_G(G_0)\otimes V^*$. 
\end{Def}

Combining \thmref{Cor}{csgrpcateq}, \thmref{Prop}{coordmor}, \thmref{Cor}{linearmor} and the \emph{formul\ae} from \thmref{Prop}{coindfunctor}, we obtain the following statement.

\begin{Prop}[linearaction]
	Let $(V,V_{0,\reals})$ be a \emph{cs} vector space, $(G_0,\ger g)$ a \emph{cs} supergroup pair, and $G=C(G_0,\ger g)$. The following data are in bijection.
	\begin{enumerate}
		\item Linear actions $\alpha:G\times V\to V$. 
		\item Even linear maps $f:V^*\to\mathcal O_G(G_0)\otimes V^*$ \scth
		\[
			(m^*\Hat\otimes\id_{V^*})\circ f=(\id_{V^*}\Hat\otimes f)\circ f\nd (e^*\Hat\otimes\id_{V^*})\circ f=\id_{V^*}\ .
		\]
		\item Elements $F\in\GL_{cs}(\mathcal O_G(G_0)\otimes V)\subset(\mathcal O_G(G_0)\otimes\underline{\mathrm{End}}(V))_0$ \scth
		\[
			(m^*\Hat\otimes\id_V)\circ F=(\id_V\Hat\otimes\,F)\circ F\nd (e^*\Hat\otimes\id_V)\circ F=\id_V\ .
		\]
		\item Morphisms $\vphi:G\to\mathcal{GL}(V,V_{0,\reals})$ of \emph{cs} Lie supergroups. 
		\item Morphisms $(\vphi,d\vphi):(G_0,\ger g)\to(\GL(V_{0,\reals})\times\GL(V_1),\underline{\mathrm{End}}_{cs}(V))$ of \emph{cs} supergroup pairs. 
	\end{enumerate}
	
	The data in (iii) and (v) are related by the equation 
	\begin{equation}\label{eq:linearaction}
		F(u;g)=\vphi(g)\circ d\vphi(u)\in\underline{\mathrm{End}}(V)\mathfa u\in\Uenv0{\ger g}\,,\,g\in G_0\ .
	\end{equation}
\end{Prop}

\begin{proof}
	Given an even linear map $f:V^*\to\mathcal O_G(G_0)\otimes V^*$, we obtain an element $F\in(\mathcal O_G(G_0)\otimes\underline{\mathrm{End}}(V))_0$. Then $f(V_{0,\reals}^*)$ consists of real-valued superfunctions if and only $F\in\End[_{cs}]0{\mathcal O_G(G_0)\otimes V}_\reals$. Assume now that $f=\alpha^*$ where $\alpha:G\times V\to V$ is a linear action.
	
	We recall now from \thmref{Par}{csmfprod} that for any \emph{cs} manifold, the global sections module of the structure sheaf is endowed with a natural nuclear Fr\'echet topology. Then, more explicitly, $F$, considered as an even element of the tensor product $\mathcal O_G(G_0)\otimes\underline{\mathrm{End}}(V)$, is given by 
	\[
		(\id\Hat\otimes\,\xi)(F(x))=f(\xi)(x)\mathfa\xi\in V^*\,,\,x\in V\ .
	\]
	In fact, this equation can easily be extended to hold for $v\in\mathcal O_G(G_0)\otimes V$ (by extending $F$ and $f$ $\mathcal O_G(G_0)$-linearly). Since $\alpha$ is an action, we have the relation $(m^*\Hat\otimes\id)\circ f=(\id\Hat\otimes f)\circ f$. Hence, we compute \fa $\mu\in V^*$, $x\in V$, 
	\begin{align*}
		(\id\Hat\otimes\id\Hat\otimes\,\xi)((\id\Hat\otimes F)(F(x)))&=(\id\Hat\otimes f(\xi))(F(v))=(\id\Hat\otimes f)(f(\xi))(x)\\
		&=(m^*\Hat\otimes\id)(f(\xi))(x)=m^*(f(\xi)(x))\\
		&=m^*((\id\Hat\otimes\,\xi)(F(x)))=(m^*\Hat\otimes\,\xi)(F(x))\ ,
	\end{align*}
	so $(\id\Hat\otimes F)\circ F=(m^*\Hat\otimes\id)\circ F$. 
	
	Moreover,
	\[
		\xi((e^*\Hat\otimes\id)(F(x)))=e^*(f(\xi)(x))=(e^*\Hat\otimes\id)(f(\xi))(x)=\xi(x)\ ,
	\]
	so $(e^*\Hat\otimes\id)\circ F=\id$. 
	
	Denote by $\delta:G\to G\times G$ the diagonal; $\delta^*$ is the algebra multiplication of $\mathcal O_G(G_0)$. Moreover, $F$, considered as an element of $\End[_{\mathcal O_G(G_0)}]0{\mathcal O_G(G_0)\otimes V}$, is $(\id\Hat\otimes F)$. The composite in the endomorphism ring of $F$ and $(i^*\Hat\otimes\id)\circ F$ is the left hand side of the following equation:
	\begin{align*}
		(\delta^*\Hat\otimes\id)\circ(\id\Hat\otimes\,i^*\Hat\otimes\id)\circ(\id\Hat\otimes F)\circ F&=\Parens1{(\delta^*\circ(\id\Hat\otimes\,i^*)\circ m^*)\Hat\otimes\id}\circ F\\
		&=(1\cdot e^*\Hat\otimes\id)\circ F=1\otimes\id\ .
	\end{align*}
	This shows that $F$ is left invertible; hence, it is invertible and thus an element 
	of $\GL_{cs}(\mathcal O_G(G_0)\otimes V)$. 
	
	Let $F\in\GL_{cs}(\mathcal O_G(G_0)\otimes V)=\mathcal{GL}(V,V_{0,\reals})(G)$ satisfy 
	\[
	(m^*\Hat\otimes\id)\circ F=(\id\Hat\otimes F)\circ F\nd(e^*\Hat\otimes\id)\circ F=\id\ .
	\]
	By the same computation as above, the inverse in $\GL_{cs}(\mathcal O_G(G_0)\otimes V)$ of $F$ is given by $F^{-1}=(i^*\Hat\otimes\id)\circ F$. 
	
	The element $F$ represents a morphism $\vphi:G\to\mathcal{GL}(V,V_{0,\reals})$ of \emph{cs} manifolds. By \thmref{Prop}{coordmor}, $\vphi\circ m$ is represented by $(m^*\Hat\otimes\id)\circ F$. If $m_0$ is the multiplication morphism of $\mathcal{GL}(V,V_{0,\reals})$, then $m_0\circ(\vphi\times\vphi)$ is represented by $(\id\Hat\otimes F)\circ F$. Similarly, $(i^*\Hat\otimes F)$ represents $\vphi\circ i$, and if $i_0$ denotes the inversion morphism of $\mathcal{GL}(V,V_{0,\reals})$, then $i_0\circ F$ is represented by $F^{-1}$. Finally, $\vphi\circ e$ is represented by $(e^*\Hat\otimes\id)\circ F$, and the unit morphism $e_0$ of $\mathcal{GL}(V,V_{0,\reals})$ is represented by $\id$. These considerations show that $\vphi$ is a morphism of \emph{cs} Lie supergroups.
	
	Morphisms of \emph{cs} Lie supergroups and of \emph{cs} supergroup pairs are in bijection. Given a morphism $(\vphi,d\vphi):(G_0,\ger g)\to(\GL(V_{0,\reals})\times\GL(V_1),\underline{\mathrm{End}}_{cs}(V))$, we may define a map $f:V^*\to\mathcal O_G(G_0)\otimes V^*\subset\mathcal O_{G\times V}(G_0\times V_{0,\reals})$ by 
	\[
		f(\xi)(u;g)(x)=(\xi\circ\vphi(g)\circ d\vphi(u))(x)
	\]
	\fa $\xi\in V^*$, $x\in V$, $u\in\Uenv0{\ger g}$, $g\in G_0$.

	Let $\alpha:G\times V\to V$ be the morphism of \emph{cs} manifolds which corresponds \emph{via} \thmref{Cor}{linearmor} to $f$. Then 
	\begin{align*}
		((m^*\Hat\otimes\id)\circ f)(\xi)(u\otimes v;g,h)(x)&=f(\mu)(\Ad(h^{-1})(u)v;gh)(x)\\
		&=(\xi\circ\vphi(gh)\circ d\vphi(\Ad(h^{-1})(u))(x)\\
		&=(\xi\circ\vphi(g)\circ d\vphi(u)\circ\vphi(h)\circ d\vphi(v))(x)\\
		&=f(\xi)(u;g)((\vphi(h)\circ d\vphi(v))(x))\\
		&=((\id\Hat\otimes f)\circ f)(\xi)(u\otimes v;g,h)(x)
	\end{align*}
	Hence, $(m^*\Hat\otimes\id)\circ f=(\id\Hat\otimes f)\circ f$. By the uniqueness statement in \thmref{Cor}{linearmor}, it follows that $\alpha\circ(m\times\alpha)=\alpha\circ(\id\times\alpha)$. Analogously, one proves that $\alpha\circ(e\times\id)=\id$. We have proved the claim. 
\end{proof}

\renewcommand\thesection{C}
\section{Appendix: Berezin integration on \emph{cs} manifolds}\label{s:app-c}

\subsection{Berezin integral and absolute Berezin integral}

\begin{Par}
	Let $(A,A_{0,\reals})$ be a supercommutative \emph{cs} algebra, and $(V,V_{0,\reals})$ a \emph{cs} vector space. The homomorphism $\mathrm{Ber}:\GL(A\otimes V)\to A^\times_0$ defined by 
	\[
		\mathrm{Ber}\begin{Matrix}1a&b\\c&d\end{Matrix}=\det(a-bd^{-1}c)(\det d)^{-1}=\det(d-ca^{-1}b)(\det a)^{-1}
	\]
	called the \emph{Berezinian}. Observe the following: If $\bar A=A/A\cdot A_1$, then a matrix as above is invertible if and only if $a$ and $d$ are, if and only if their images in $\End[_{\bar A_0}]0{\bar A_0\otimes V_j}$ are (where $j=0$ and $j=1$, respectively). 
	
	If $(A,A_{0,\reals})=\Gamma(Z)=(\mathcal O_Z(Z_0),\mathcal O_Z(Z_0)_{0,\reals})$ for some supermanifold $Z$, then one can define a further homomorphism $\Abs0{\mathrm{Ber}}:\GL_{cs}(A\otimes V)\to A^\times_0$ by 
	\[
		\Abs0{\mathrm{Ber}}\begin{Matrix}1a&b\\c&d\end{Matrix}=\sgn\widetilde{\det a}\cdot\mathrm{Ber}\begin{Matrix}1a&b\\ c&d\end{Matrix}\ .
	\]
	Here, recall that $\tilde f$ is the function underlying the superfunction $f$, $\sgn z=\Abs0z^{-1}z$ for $z\in\cplxs^\times$, and $\sgn 0=0$. The homomorphism $\Abs0{\mathrm{Ber}}$ is slightly non-standard, \emph{cf.}~\cite{voronov-geomint} where the notation $\mathrm{Ber}_{1,0}$ is used. Compare also \cite{shander-semiorient}. 
\end{Par}

\begin{Par}[ber-absber]
	Let $M$ be a free graded $A$-module of graded dimension $p|q$. We define an $A$-module $\Ber[_A]0M=\Ber0M$ of graded rank $1|0$ (if $q$ is even) resp.~$0|1$ (if $q$ is odd) as follows. With a basis $x_1,\dotsc,x_p,\xi_1,\dotsc,\xi_q$ where the $x_j$ are even and the $\xi_j$ are odd, we associate a distinguished basis $D(x,\xi)$ of $\Ber0M$, of parity $\equiv q\,(2)$. 
	
	If $y_1,\dotsc,y_p,\eta_1,\dotsc,\eta_q$ is another such basis, related to $x,\xi$ by 
	\[
	y_i=\textstyle\sum_ja_{ij}x_j+b_{ij}\xi_j\nd\eta_i=\sum_jc_{ij}x_j+d_{ij}\xi_j\ ,
	\]
	where $a_{ij},d_{ij}\in A_0$, $b_{ij},c_{ij}\in A_1$, then 
	\[
		D(y,\eta)=\mathrm{Ber}\begin{Matrix}0a&b\\c&d\end{Matrix}\cdot D(x,\xi)\ .
	\]
	
	Assume that we have $(A,A_{0,\reals})=\Gamma(Z)$ for some \emph{cs} manifold $Z$ (\emph{cf.}~\thmref{Rem}{csalg}); moreover, suppose that we are given a choice $M_{0,\reals}\subset M_0$ of a maximal proper $A_{0,\reals}$-submodule \scth $M_0=\Span0{M_{0,\reals}}_\cplxs$ and $A_1\cdot M_1\subset M_{0,\reals}$ where $\Span0\cdot_\cplxs$ denotes complex linear span. (Because $A_1\subset A_{0,\reals}$, such submodules manifestly exist whenever $A_0\neq0$ and $M_0\neq0$.) We define $A_\reals=A_{0,\reals}\oplus A_1$. This is a real graded subalgebra of $A$. 
		
	We now define a free $A$-module $\ABer[_{A,A_{0,\reals}}]0{M,M_{0,\reals}}$ of graded rank $1|0$ (if $q$ is even) resp.~$0|1$ (if $q$ is odd). With any graded basis $x,\xi$ of the $A_\reals$-module $M_\reals=M_{0,\reals}\oplus M_1$, we associate a distinguished basis $\Abs0{D(x,\xi)}$.
	
	If $y,\eta$ is related to $x,\xi$ as above, where now $a_{ij},d_{ij}\in A_{0,\reals}$, $b_{ij},c_{ij}\in A_1$, then we require 
	\[
		\Abs0{D(y,\eta)}=\Abs0{\mathrm{Ber}}\begin{Matrix}0a&b\\ c&d\end{Matrix}\cdot\Abs0{D(x,\xi)}\ .
	\]
\end{Par}

\begin{Par}[ber-absber-sheaf]
	Let $X=(X_0,\mathcal O)$ be a \emph{cs} manifold of graded dimension $p|q$. If $U$ is a coordinate neighbourhood and $(x,\xi)$ is a coordinate system, then any superfunction $f\in\mathcal O(U)$ may be written uniquely in the form 
	\[
	f=\sum\nolimits_If_I(x_1,\dotsc,x_p)\xi^I\mathtxt{where}f_I\in\Ct[^\infty]0{x(U),\cplxs}\ ,
	\]
	the sum extends over all $I=(1\sle i_1<\dotsm<i_k\sle q)$, and $\xi^I=\xi_{i_1}\dotsm\xi_{i_q}$. Here, $h(x_1,\dotsc,x_p)$, for $h\in\Ct[^\infty]0{x(U),\cplxs}$, is to be understood in the following way: Let $\vphi:U\to x(U)\subset\reals^p$ be the morphism determined by $\vphi^*(\pr_j)=x_j$; then $h(x_1,\dotsc,x_p)=\vphi^*h$.  
	
	Define even derivations $\frac\partial{\partial x_i}$ of $\mathcal O(U)$ by 
	\[
		\frac{\partial f}{\partial x_i}=\sum\nolimits_I(\partial_if_I)(x_1,\dotsc,x_p)\xi^I\ .
	\]
	In particular, $\frac{\partial x_\ell}{\partial x_i}=\delta_{i\ell}$. Moreover, define odd derivations $\frac\partial{\partial\xi_j}$ by 
	\[
		\frac{\partial\xi_k}{\partial\xi_j}=\delta_{jk}\nd\frac{\partial f}{\partial\xi_j}=\sum\nolimits_If_I(x_1,\dotsc,x_p)\frac{\partial\xi^I}{\partial\xi_j}\ .
	\]
	Then $\frac\partial{\partial x_i}$, $\frac\partial{\partial\xi_j}$ form an $\mathcal O(U)$-basis of $\Der0{\mathcal O(U)}$. The $\mathcal O(U)$-dual $\Form[^1_X]0U$ has the basis $dx_1,\dotsc,dx_p,d\xi_1,\dotsc,d\xi_q$ where
	\[
		\Dual1{dx_i}{\tfrac\partial{\partial x_j}}=\Dual1{d\xi_i}{\tfrac\partial{\partial\xi_j}}=\delta_{ij}\nd
		\Dual1{dx_i}{\tfrac\partial{\partial\xi_j}}=\Dual1{d\xi_i}{\tfrac\partial{\partial x_j}}=0\ .
	\]
	Let $\Form[^1_X]0U_{0,\reals}$ be the $\mathcal O(U)_{0,\reals}$ submodule spanned by $dx_1,\dotsc,dx_p$. One sees easily that this submodule is in fact independent of the choice of basis. (The point is that the $x_j$ are real-valued for any coordinate system.)
	
	We set $\mathcal Ber_X(U)=\Ber[_{\mathcal O(U)}]0{\Form[^1_X]0U}$. Denote the associated sheaf of $\mathcal O$-modules by $\mathcal Ber_X$; it is called the \emph{Berezinian sheaf} of $X$. Similarly, set 
	\[
	\Abs0{\mathcal Ber}_X(U)=\ABer[_{\mathcal O(U),\mathcal O(U)_{0,\reals}}]0{\Form[^1_X]0U,\Form[^1_X]0U_{0,\reals}}\ ,
	\]
	and let $\Abs0{\mathcal Ber}_X$ be the associated sheaf, called the \emph{absolute Berezinian sheaf}. 
	
	If $f:X\to Y$ is an isomorphism of supermanifolds, then we define a sheaf morphism $f^*:\mathcal Ber_Y\to f_*\mathcal Ber_X$ as follows. For any coordinate neighbourhood $V\subset Y_0$ and any coordinate system $(y,\eta)$ on $V$, we let 
	\[
	f^*\Parens1{h\cdot D(dy,d\eta)}=(f^*h)\cdot D(df^*y,df^*\eta)\mathfa h\in\mathcal O_Y(V)\ .
	\]
	If $(z,\zeta)$ is another coordinate system, then $D\Parens1{dz,d\zeta}=\Ber0J\cdot D\Parens1{dy,d\eta}$ where the Jacobian is given by 
	\[
		J^{st}=\begin{Matrix}1\frac{\partial z}{\partial y}&\frac{\partial\zeta}{\partial y}\\\frac{\partial z}{\partial\eta}&\frac{\partial\zeta}{\partial\eta}\end{Matrix}\ ,
	\]
	the superscript ${}^{st}$ denoting the super-transpose. It follows by the multiplicative property of the Berezinian that $f^*$ gives a well-defined sheaf morphism. 
	
	One can now proceed in exactly the same way to define a morphism of sheaves $f^*:\Abs0{\mathcal Ber}_{\mathcal Y}\to f_*\Abs0{\mathcal Ber}_{\mathcal X}$. What enters here crucially is that $\frac{\partial z}{\partial y}$ is real-valued.
\end{Par}

\begin{Def}[ber-absber-int]
	A \emph{cs} manifold $X$ is called \emph{evenly oriented} if the underlying manifold is oriented. For an evenly oriented supermanifold of graded dimenion $p|q$, we define as follows an even linear presheaf morphism $\int_X:\Gamma_c(\mathcal Ber_X)\to\cplxs_{X_0}$ called the \emph{Berezin integral}---here, $\Gamma_c$ denotes the set of compactly supported sections. 
	
	Let $(U_\alpha)$ be a locally finite cover of $X_0$ by coordinate neighbourhoods, $(x^\alpha,\xi^\alpha)$ coordinate systems where $x_1^\alpha,\dotsc,x_p^\alpha$ is are oriented coordinates systems of the underlying manifold $X_0$, and $(\chi_\alpha)$ an $\mathcal O$-partition of unity subordinate to $(U_\alpha)$. 
	
	For any $f\in\mathcal O_X(U_\alpha)$, we write $f=\sum_If_I\cdot\xi^{\alpha I}$ where $f_I=h_I(x_1^\alpha,\dotsc,x_p^\alpha)$ \fs $h_I\in\Ct[^\infty]0{x^\alpha(U_\alpha),\cplxs}$. If $f$ is compactly supported, let
	\[
		\int_\alpha D(dx^\alpha,d\xi^\alpha)\cdot f=(-1)^{pq}\cdot\int_{U_\alpha}d\tilde x^\alpha\cdot\tilde f_{1,\dotsc,1}\in\cplxs\ ,
	\]
	where $d\tilde x^\alpha=d\tilde x_1^\alpha\wedge\dotsm\wedge d\tilde x_p^\alpha$. Then, we define, \fa $\omega\in\Gamma_c(\mathcal Ber_X)$,
	\[
		\int_X\omega=\sum\nolimits_\alpha\int_\alpha\chi_\alpha\cdot\omega\ .
	\]
	
	The Berezin integral is well-defined, independent of all choices \cite[Theorem 2.4.5]{leites}. Moreover, if $f:X\to Y$ is an isomorphism of even oriented \emph{cs} manifolds preserving the orientations of the underlying manifolds, then 
	\[
		\int_Xf^*\omega=\int_Y\omega\mathfa\omega\in\Gamma_c(\mathcal Ber_Y)\ .
	\]

	In exactly the same way, we define a presheaf morphism $\int_X:\Gamma_c(\Abs0{\mathcal Ber}_X)\to\cplxs_{X_0}$ for \emph{any} \emph{cs} manifold $X$ (possibly without even orientation) by using the integral of densities $\int_{U_\alpha}\Abs0{d\tilde x^\alpha}\cdot\tilde f_{1,\dotsc,1}$. The resulting integral is invariant under \emph{all} isomorphisms, \emph{irrespective of preservation or even existence of orientations}. 
	
	We shall use the following lemma repeatedly.
\end{Def}

\begin{Lem}[berint-nondegenpair]
	Let $Z$ be a \emph{cs} manifold and $f\in\Abs0{Ber}_Z(Z_0)$. If $\int_Zf\cdot h=0$ \fa $h\in\Gamma_c(Z_0,\mathcal O_Z)$, then $f=0$. 
\end{Lem}

\begin{proof}
	We argue by contraposition. Then we may assume that $Z=\reals^{p|q}_\cplxs$ as a \emph{cs} manifold, and that $f=\vphi\cdot\Abs0{D(x,\xi)}$ where $0\in(\supp\vphi)^\circ$ and $x,\xi$ is the standard coordinate system. 
	
	Write $\vphi=\sum_I\vphi_I\xi^I$ with $\vphi_I\in\Ct[^\infty]0{\reals^p,\cplxs}$ (this is possible, since $x$ are the standard coordinates). There is some multi-index $I=(1\sle i_1<\dotsm i_m\sle q)$ \scth $\vphi_I(0)\neq0$. Let $J=(1-i_1,\dotsc,1-i_m)$ and $h=\xi^J\cdot\overline{\vphi_I}\cdot\chi$ where $\chi:\reals^p\to[0,1]$ is smooth, of compact support, and satisfies $\chi(0)=1$. Thus 
	\[
		\int_{\reals^{p|q}_\cplxs} f\cdot h=\pm\int_{\reals^{p|q}_\cplxs}\xi^{(1,\dotsc,1)}\cdot\Abs0{\vphi_I}^2\cdot\chi\cdot\Abs0{D(x,\xi)}=\pm\int_{\reals^p}\Abs0{\vphi_I}^2\cdot\chi\neq0\ .
	\]
	This proves the lemma.
\end{proof}

\begin{Par}
	For any morphism $\vphi:X\to Y$, and any open subset $U\subset Y$, we consider the set $\Gamma^\vphi_{cf}(U,\Abs0{\mathcal Ber}_X)$ of all local sections $\omega\in\Gamma(\vphi^{-1}(U),\Abs0{\mathcal Ber}_X)$ \scth $\vphi:\supp\omega\to Y$ is a proper map. This defines a presheaf $\Gamma_{cf}^\vphi$ on $Y$ , the presheaf of sections \emph{compactly supported along the fibres of $\vphi$}. 
	
	Assume that $\vphi$ is submersive and that its underlying map $\vphi_0$ is surjective. Then there is a well-defined even presheaf morphism $\vphi_!:\Gamma_{cf}^\vphi(\Abs0{\mathcal Ber}_X)\to\Abs0{\mathcal Ber}_Y$ of (Berezin) integration along the fibres.  Note that we give $\Abs0{\mathcal Ber}_X$ the parity $\eps^q$ where $\eps=0|1$ and $\dim X=p|q$. This makes $\vphi_!$ an even morphism. Two fundamental properties of $\vphi_!$ are 
	\[
		\supp\vphi_!(\omega)\subset\vphi(\supp\omega)\nd\vphi_!(\vphi^*(f)\cdot\omega)=f\cdot\vphi_!(\omega)\ .
	\]
	Compare \cite[Proposition 5.7]{ah-berezin} for the definition of $\vphi_!$.
	
	In the case of a projection $p_1:Y\times F\to Y$, one may define more generally a fibre integration map $p_{1!}:\Gamma_{cf}^{p_1}(\mathcal E\otimes\Abs0{\mathcal Ber}_F)\to\Gamma(\mathcal E)$, for any sheaf $\mathcal E$ on $Y$. 
\end{Par}

\subsection{Invariant Berezin integration}

In this section, we will review some results from \cite{ah-berezin} concerning invariant Berezin integration on homogenenous \emph{cs} manifolds, reformulating them in the language of \emph{cs} supergroup pairs. 

\begin{Par}[quotient]
	Let $G$ be a \emph{cs} Lie supergroup. As in \cite{ah-berezin}, one can define the concept of a quotient by an action of $G$ and show that quotients exist for free and proper $G$-actions. 
	
	In particular, assume that $H$ is a \emph{closed \emph{cs} subsupergroup}. \emph{I.e.}, we have a morphism $f:H\to G$ of \emph{cs} Lie supergroups which is, on the level of spaces, a closed embedding of $H_0$ in $G_0$, and $f^*:\mathcal O_G\to f_*\mathcal O_H$ is an epimorphism. (We will suppress $f$ from the notation.) Then the quotient \emph{cs} manifold $G/H=(G_0/H_0,\mathcal O_{G/H})$ exists. Here, \fa open subsets $U\subset G_0/H_0$, 
	\[
		\mathcal O_{G/H}(U)=\Set1{f\in\mathcal O_G(\pi^{-1}(U))}{m^*f=p_1^*f\in\mathcal O_{G\times H}(\pi^{-1}(U)\times H_0)}\ ,
	\]
	where $(m,m^*):G\times H\to G$ is the restriction of the multiplication morphism of $G$. On the level of sheaves, the canonical morphism $\pi:G\to G/H$ maps $f$ to $f=\pi^*f$ (\emph{i.e.}~$\pi^*$ is the inclusion). We will often view $\mathcal O_{G/H}$ as a subsheaf of $\mathcal O_G$. 
	
	If $G=C(G_0,\ger g)$ where $(G_0,\ger g)$ is a \emph{cs} supergroup pair, and $H=C(H_0,\ger h)$ where the \emph{cs} supergroup pair is given by a closed subgroup $H_0\subset G_0$ and a Lie subsuperalgebra $\ger h\subset\ger g$, then we have 
	\[
		\mathcal O_{G/H}(U)=\mathcal O_G(\pi^{-1}(U))^{H_0,\ger h}=\underline{\mathrm{Hom}}_{\Uenv0{\ger g_0},H}\Parens1{\Uenv0{\ger g}/\Uenv0{\ger g}\ger h,\Ct[^\infty]0{\pi^{-1}(U)}}
	\]
	\fa open $U\subset G_0/H_0$. Here, the superscript ${}^{H_0,\ger h}$ denotes simultaneous $H_0$- and $\ger h$-invariants, and we recall the definition of the actions from the equations \eqref{eq:coind-infinaction} and \eqref{eq:coind-grpaction}. 
	
	Moreover, $\mathcal O_{G/H}$ inherits actions by $\ger g$ and $G_0$ from $\mathcal O_G$, namely
	\begin{equation}\label{eq:gh-infinaction}
		(\ell_uf)(v;g)=(-1)^{\Abs0f\Abs0u}f\Parens1{\Ad(g^{-1})(S(u))v;g}
	\end{equation}
	\fa $f\in\mathcal O_{G/H}(U)$, $u,v\in\Uenv0{\ger g}$, $g\in\pi^{-1}(U)$, and 
	\begin{equation}\label{eq:gh-grpaction}
		(\ell_gf)(u;h)=f(u;g^{-1}h)
	\end{equation}
	\fa $f\in\mathcal O_{G/H}(U)$, $u\in\Uenv0{\ger g}$, $g\in G_0$, $h\in g\cdot\pi^{-1}(U)$, respectively.  
\end{Par}

\begin{Par}
	If $G$ is of graded dimension $p|q$, then \fa open $U\subset G_0/H_0$, there are canonical isomorphisms 
	\begin{gather*}
		\mathcal Ber_{G/H}(U)\cong\Parens1{\underline{\mathrm{Hom}}_{\Uenv0{\ger g_0}}\Parens0{\Uenv0{\ger g},\Ct[^\infty]0{\pi^{-1}(U)}}\otimes\Ber0{(\ger g/\ger h)^*}}^{H_0,\ger h}\ ,\\
		\Abs0{\mathcal Ber}_{G/H}(U)\cong\Parens1{\underline{\mathrm{Hom}}_{\Uenv0{\ger g_0}}\Parens0{\Uenv0{\ger g},\Ct[^\infty]0{\pi^{-1}(U)}}\otimes\ABer0{(\ger g/\ger h)^*}}^{H_0,\ger h}\ ,
	\end{gather*}
	where we abbreviate 
	\[
		\ABer0{(\ger g/\ger h)^*}=\ABer[_{\cplxs,\reals}]0{(\ger g/\ger h)^*,(\ger g_{0,\reals}/\ger h_{0,\reals})^*}\ .
	\]
	Here, we use the canonical isomorphism $\Ber[_R]0{R\otimes V}=R\otimes\Ber0V$, and the characterisation of $\mathcal Ber_{G/H}$ from \cite[Corollary 4.12]{ah-berezin}. (In this reference, we consider real supermanifolds, but as should be clear by now, everything goes through for \emph{cs} manifolds with only minor changes.) For the absolute Berezinians, essentially the same argument goes through. 
	
	If, according to the above isomorphism, we consider the elements of $\mathcal Ber_{G/H}(U)$ as maps $\Uenv0{\ger g}\times\pi^{-1}(U)\to\Ber0{(\ger g/\ger h)^*}$, then the actions of $\ger g$ and $G_0$ induced from $\mathcal O_{G/H}$ may be expressed simply by 
	\begin{equation}\label{eq:ber-infinaction}
		(\ell_uf)(v;g)=(-1)^{\Abs0f\Abs0u}f\Parens1{\Ad(g^{-1})(S(u))v;g}
	\end{equation}
	\fa $f\in\mathcal Ber_{G/H}(U)$, $u,v\in\Uenv0{\ger g}$, $g\in\pi^{-1}(U)$, and
	\begin{equation}\label{eq:ber-grpaction}
		(\ell_gf)(u;h)=f(u;g^{-1}h)
	\end{equation}
	\fa $f\in\mathcal Ber_{G/H}(U)$, $u\in\Uenv0{\ger g}$, $g\in G_0$, $h\in g\cdot\pi^{-1}(U)$, respectively. The same holds for the absolute Berezinians. 
\end{Par}
	
\begin{Par}
	With the aid of \thmref{Prop}{linearaction}, one sees that $\Ad_G:G\to\mathcal{GL}(\ger g)$, the adjoint morphism, is represented by $\Ad_{G_0,\ger g}\in\Gamma(\mathcal O_G)\otimes\underline{\mathrm{End}}\Parens0{\ger g}$, defined by 
	\begin{equation}\label{eq:ad-pairexpr}
		\Ad_{G_0,\ger g}(u;g)=\Ad(g)\circ\ad(u)\in\underline{\mathrm{End}}\Parens0{\ger g}\mathfa g\in G_0\,,\,u\in\Uenv0{\ger g}\ .
	\end{equation}
	Similarly, $\Ad_G\!|_H:H\to\underline{\mathrm{End}}\Parens0{\ger g}$ is represented by $\Ad_{G_0,\ger g}\!|_H\in\Gamma(\mathcal O_H)\otimes\underline{\mathrm{End}}\Parens0{\ger g}$, given by 
	\[
		\Ad_{G_0,\ger g}|_H(u;h)=\Ad(h)\circ\ad(u)\in\underline{\mathrm{End}}\Parens0{\ger g}\mathfa h\in H_0\,,\,u\in\Uenv0{\ger h}\ .
	\]
	Analogously, one has $\Ad_{H_0,\ger h}\in\Gamma(\mathcal O_H)\otimes\underline{\mathrm{End}}\Parens0{\ger h}$. 
	
	By \cite[Theorem 4.13]{ah-berezin}, $\Gamma(\mathcal Ber_{G/H})$ contains a non-zero $(G_0,\ger g)$-invariant element if and only if 
	\begin{equation}\label{eq:unimodcond}
		\Ber0{\Ad_{H_0,\ger h}}=\Ber0{\Ad_{G_0,\ger g}\!|_H}\in\Gamma(\mathcal O_H)\ .
	\end{equation}
	Here, $\Ber0{\Ad_{H_0,\ger h}}\in\Gamma(\mathcal O_H)_0^\times$ is given by 
	\[
		\Ber0{\Ad_{H_0,\ger h}}(u;h)=\Ber[_{\ger h}]0{\Ad(h)}\cdot\str_{\ger h}\ad(u)\mathfa u\in\Uenv0{\ger h}\,,\,h\in H_0\ .
	\]
	Here, $\str_{\ger g}\circ\ad:\Uenv0{\ger g}\to\cplxs$ is the unique extension to an even unital algebra morphism of the map $\ger g\to\Uenv0{\ger g}:x\mapsto\str_{\ger g}\ad(x)$. A similar equation defines $\Ber0{\Ad_{G_0,\ger g}\!|_{\mathcal H}}$. 
	
	Such a non-zero $(G_0,\ger g)$-invariant element (if it exists), is unique up to scalar multiples. In this case, we say that the $G$-space $G/H$ is \emph{geometrically unimodular}. We say that the \emph{cs} Lie supergroup $G$ is geometrically uni\-mo\-du\-lar if it is so as a $G\times G$-space (it always is as a $G$-space).
\end{Par}

\begin{Par}
	Similarly, $\Abs0{\mathcal Ber}_{G/H}$ possesses a non-zero $(G_0,\ger g)$-invariant global section if and only if 
	\begin{equation}\label{eq:anunimodcond}
		\ABer0{\Ad_{H_0,\ger h}}=\ABer0{\Ad_{G_0,\ger g}\!|_{\mathcal H}}\in\Gamma(\mathcal O_H)\ .
	\end{equation}
	Here, $\ABer0{\Ad_{H_0,\ger h}}\in\Gamma(\mathcal O_H)_0^\times$ is given by  
	\[
		\ABer0{\Ad_{H_0,\ger h}}(u;h)=\ABer[_{\ger h}]0{\Ad(h)}\cdot\str_{\ger h}\ad(u)\mathfa u\in\Uenv0{\ger h}\,,\,h\in H_0\ .
	\]
	Whenever \eqref{eq:anunimodcond} is satisfied, we say that the $G$-space $G/H$ is \emph{analytically unimodu\-lar}; we say that the \emph{cs} Lie supergroup $G$ is analytically unimodular if it is as a $G\times G$-space (it always is as a $G$-space).
\end{Par}

\begin{Prop}[unimodconds]
	Let $(G_0,\ger g)$ and $(H_0,\ger h)$ be \emph{cs} supergroup pairs where $H_0\subset G_0$ is a closed subgroup and $\ger h\subset\ger g$ is a Lie subsuperalgebra. Let $G=C(G_0,\ger g)$ and $H=C(H_0,\ger h)$. 
	\begin{enumerate}
		\item If $G$ and $H$ are geometrically (analytically) unimodular, then so is the $G$-space $G/H$.
		\item If $G$ and $H$ are geometrically (analytically) unimodular, then so is the direct product $G\times H$ of \emph{cs} Lie supergroups. 
		\item If $\ger g$ is nilpotent (in particular, if it is Abelian) and $G_0$ is connected, then $G$ is geometrically and analytically unimodular.
		\item If $\ger g$ is strongly reductive (\emph{cf.}~\thmref{Def}{red}) and $G_0$ is connected, then $G$ is geometrically and analytically unimodular. 
		\item If $\Ad_\ger g(G_0)\subset\GL(\ger g)$ is compact, then $G$ is analytically and geometrically unimodular. 
	\end{enumerate}
\end{Prop}

\begin{proof}
	Statements (i) and (ii) are immediate from \eqref{eq:unimodcond} and \eqref{eq:anunimodcond}. To prove that $G$ is geometrically unimodular, it wil be sufficient to prove that $\Ber[_\ger g]0{\Ad(g)}=1$ \fa $g\in G_0$ and $\str_\ger g\ad(x)=0$ \fa $x\in\ger g$. The latter condition is always verified for $x\in\ger g_1$, for reasons of parity. Thus, it will suffice to prove the former, since the latter then follows by differentiation; conversely, when $G_0$ is connected, the latter condition for $x\in\ger g_0$ implies the former. A similar reasoning holds true for the case of analytic unimodularity. 
	
	If $\ger g$ is nilpotent, then $\tr_{\ger g_0}\ad(x)=\tr_{\ger g_1}\ad(x)=0$ \fa $x\in\ger g_0$; thus, $G$ is geometrically and analytically unimodular. If $\ger g$ is strongly reductive, then one has $\ger g=\ger z(\ger g)\oplus\ger g'$. We have $\ad(\ger z(\ger g))=0$. But $\str_\ger g([a,b])=0$ for any $a,b\in\End0{\ger g}$. In particular, $\str_\ger g\ad(\ger g')=0$. Finally, if $\Ad_\ger g(G_0)$ is compact, then $\ABer0{\Ad_\ger g(G_0)}$ and $\Ber0{\Ad_\ger g(G_0)}$ are compact subgroups of $(\cplxs^\times,\cdot)$, and hence trivial.
\end{proof}

\begin{Par}
	If $G/H$ is a geometrically unimodular $G$-space, then $\Ber0{(\ger g/\ger h)^*}$ is a trivial $H_0$- and $\ger h$-module. Thus, in this case, for all open $U\subset G_0/H_0$, 
	\begin{equation}\label{eq:berunimod}
		\mathcal Ber_{G/H}(U)\cong\mathcal O_{G/H}(U)\otimes\Ber0{(\ger g/\ger h)^*}\ ;
	\end{equation}
	similarly for the absolute Berezinians if $G/H$ is analytically unimodular. 
	
	By Equations \eqref{eq:ber-infinaction} and \eqref{eq:ber-grpaction} for the $(G_0,\ger g)$-action, a non-zero invariant element of $\Gamma(\mathcal Ber_{G/H})$ is necessarily of the form $1\otimes\omega$ for a non-zero $\omega\in\Ber0{(\ger g/\ger h)^*}$. 
\end{Par}

\begin{Par}
	If $f\in\Gamma_c(\mathcal O_G)$ and an invariant $0\neq\Abs0{Dg}\in\Gamma(\Abs0{\mathcal Ber}_G)$ is fixed, we define $\int_Gf\,\Abs0{Dg}$ as the integral of the absolute Berezinian $f\cdot\Abs0{Dg}$. If $G$ is unimodular, then $i^*\Abs0{Dg}$ is $G\times G$-invariant again, and hence proportional to $\Abs0{Dg}$. Since there exist compactly supported superfunctions $f$ on $G$ \scth $\int_Gf\,\Abs0{Dg}\neq0$ and $i^*f=f$, we find that $i^*\Abs0{Dg}=\Abs0{Dg}$ in this case. 
	
	If $G/H$ is an analytically unimodular $G$-space, and $\Abs0{D\dot g}$ is a non-zero and invariant absolute Berezinian, then we define, for $f\in\Gamma_c(\mathcal O_G)$, 
	$\int_{G/H}f\,\Abs0{D\dot g}$ as the integral of the absolute Berezinian $f\cdot\Abs0{D\dot g}$. 
	
	Moreover, by \cite[Corollary 5.12]{ah-berezin} (which holds analogously for \emph{cs} Lie supergroups and absolute Berezinians), invariant absolute Berezinians can always be normalised \scth 
	\begin{equation}\label{eq:fubfmla}
		\int_Gf\,\Abs0{Dg}=\int_{G/H}p_{1!}\Parens1{m^*f\cdot(1\otimes\Abs0{Dh})}\,\Abs0{D\dot g}\mathfa f\in\Gamma_c(\mathcal O_G)
	\end{equation}
	where $m:G\times H\to G$ is multiplication, and $p_1:G\times H\to G$ the first projection.
	
	Let $\alpha:G\times G/H\to G/H$ denote the action $G$ on $G/H$ induced by left multiplication, and $i:G\to G$ the inversion. Then the invariance of the absolute Berezinian implies 
	\begin{equation}\label{eq:invarfibint}
		p_{1!}\Parens1{\alpha^*f\cdot p_2^*h\cdot(1\otimes\Abs0{D\dot g})}=p_{1!}\Parens1{p_2^*f\cdot(i\times\id)^*\alpha^*h\cdot(1\otimes\Abs0{D\dot g})}
	\end{equation}
	\fa $f,h\in\Gamma_c(\mathcal O_{G/H})$. In particular, if $H'$ is a closed \emph{cs} Lie subsupergroup, then 
	\begin{equation}\label{eq:invarint}
		\int_{H'\times G/H}\alpha^*f\cdot p_2^*f'\,\Abs0{Dh'}\,\Abs0{D\dot g}=\int_{H'\times G/H}p_2^*f\cdot(i\times\id)^*\alpha^*f'\,\Abs0{Dh'}\,\Abs0{D\dot g}
	\end{equation}
	\fa $f,f'\in\Gamma_c(\mathcal O_{G/H})$.
\end{Par}

\begin{Par}
	Let $U=C(U_0,\ger u)$ be a \emph{cs} Lie supergroup, and let $M=C(M_0,\ger m)$ and $H=C(H_0,\ger h)$ be closed subsupergroups \scth the restriction of the multiplication morphism $m:M\times H\to U$ defines an isomorphism onto an open subspace $V$ of $U$. For instance, this is clearly the case if $\ger u=\ger m\oplus\ger h$ as super-vector spaces. Similarly as for \cite[Proposition 5.16]{ah-berezin}, we have, for a suitable normalisation of the invariant absolute Berezinians,
	\begin{equation}\label{eq:prodsubgroup}
		\int_Uf\,\Abs0{Du}=\int_{M\times H}m^*f\cdot\frac{\pr_2^*\ABer0{\Ad_{H_0,\ger h}}}{\pr_2^*\ABer0{\Ad_{U_0,\ger u}\!|_H}}\,\Abs0{Dm}\,\Abs0{Dh}
	\end{equation}
	\fa $f\in\Gamma_c(U_0,\mathcal O_U)$ \scth $\supp f\subset V$. In fact, by \cite[Lemma 5.15]{ah-berezin}, the correct normalisation is given by $\Abs0{Du}=1\otimes\omega_1\otimes\omega_2$, $\Abs0{Dm}=1\otimes\omega_1$, $\Abs0{Dh}=1\otimes\omega_2$ in \eqref{eq:berunimod} where $\omega_1\otimes\omega_2\in\ABer0{\ger m^*}\otimes\ABer0{\ger h^*}=\ABer0{\ger u^*}$.
		
	Equation \eqref{eq:prodsubgroup} is formally quite similar to the Lie group case. But care is to be taken here: The product and inverse of functions are to be computed in the algebra $\mathcal O_{M\times H}(M_0\times H_0)$, and not in the pointwise sense. We give an explicit formula in the following lemma.
\end{Par}

\begin{Lem}[berfactor]
	Let $U=C(U_0,\ger u)$ be a \emph{cs} Lie supergroup and $H=C(H_0,\ger h)$ a closed \emph{cs} subsupergroup. Then 
	\[
		\Parens1{\ABer0{\Ad_{H_0,\ger h}}\ABer0{\Ad_{U_0,\ger u}\!|_H}^{-1}}(u;h)=\ABer[_{\ger u/\ger h}]0{\Ad(h^{-1})}\str_{\ger u/\ger h}\ad(S(u))
	\]
	\fa $u\in\Uenv0{\ger h}$, $h\in H_0$. 
\end{Lem}

\begin{proof}
	We compute explicitly. Define superfunctions $f$ and $g$ on $H$ by 
	\[
		f(u;h)=\str_{\ger h}\ad(u)\nd g(u;h)=\str_{\ger u}\ad(S(u))\ .
	\]
	For $x\in\ger h$, we have $\str_{\ger u}\ad(x)=\str_{\ger u/\ger h}\ad(x)+\str_{\ger h}\ad(x)$, so 
	\[
		(f\cdot g)(x;h)=\mu\circ((\str_{\ger h}\circ\ad)\otimes(\str_{\ger u}\circ\ad\circ S))(\Delta(x))=\str_{\ger u/\ger h}\ad(S(x))\ .
	\]

	Next, let $u'=xu$ where $u'\in\Uenv0{\ger h}$ and $x\in\ger h$. We write $s_{\ger h}=\str_{\ger h}\circ\ad$ and $s_{\ger u}=\str_{\ger u}\circ\ad$. Then, writing $\Delta(u)=\sum_ju_j\otimes v_j$, we deduce by induction  
	\begin{align*}
		(f&\cdot g)(u';h)\\
		&=\mu\Parens1{(s_{\ger h}(x)\otimes 1+1\otimes s_{\ger u}(S(x)))\cdot((s_{\ger h}\otimes s_{\ger u}\circ S)\circ\Delta)(u)}\\
		&=\sum\nolimits_js_{\ger h}(x)s_{\ger h}(u_j)s_{\ger u}(S(v_j))+(-1)^{\Abs0x\Abs0{u_j}}s_{\ger h}(u_j)s_{\ger u}(S(x))s_{\ger u}(S(v_j))\\
		&=s_{\ger h}(x)\sum\nolimits_js_{\ger h}(u_j)s_{\ger u}(S(v_j))+(-1)^{\Abs0x\Abs0u}\sum\nolimits_js_{\ger h}(u_j)s_{\ger u}(S(v_j))s_{\ger u}(S(x))\\
		&=s_{\ger h}(x)s_{\ger u/\ger h}(S(u))+(-1)^{\Abs0x\Abs0u}s_{\ger u/\ger h}(S(u))s_{\ger u}(S(x))\\
		&=s_{\ger u/\ger h}(S(xu))+(1-(-1)^{\Abs0x\Abs0u})\cdot s_{\ger h}(x)s_{\ger u/\ger h}(S(u))
	\end{align*}
	If $\Abs0x=\Abs0u=1$, then $s_{\ger h}(x)=s_{\ger u/\ger h}(S(u))=0$. Finally, if $\Abs0x\neq\Abs0u$ or $\Abs0x=\Abs0u=0$, then $(-1)^{\Abs0x\Abs0u}=1$. Thus, in any case, 
	\[
	(f\cdot g)(u';h)=s_{\ger u/\ger h}(S(xu))=\str_{\ger u/\ger h}\ad S(u')\ .
	\]
	The $H_0$-dependent parts of the superfunctions occurring in the assertion of the lemma are easily treated, and the claim follows. 
\end{proof}

\bibliographystyle{alpha}%
\bibliography{a-hchom}%

\end{document}